\documentclass[aap,MSNbibl,seceqn,citesort,dvips]{arximspdf}
\usepackage{mathbh}

\doi{10.1214/09-AAP666}
\volume{20}
\issue{5}
\pubyear{2010}
\firstpage{1907}
\lastpage{1965}

\makeatletter
\newtheorem{Theorem}{Theorem}[section]
\newtheorem{Lemma}[Theorem]{Lemma}
\newproclaim{Construction}[Theorem]{Construction}

\newtheorem{Proposition}[Theorem]{Proposition}
\newproclaim{Remark}[Theorem]{Remark}

\DeclareMathAlphabet\mathcaligr{OMS}{cmsy}{m}{n}
\renewcommand{\mathcal}{\mathcaligr}
\newcommand{\cal}{\mathcal}
\newcommand{\eqref}[1]{(\ref{#1})}

\newcommand{\prob}{\mathbb{P}}
\newcommand{\expec}{\mathbb{E}}
\newcommand{\bfd}{\mathbf{d}}
\newcommand{\Var}{\operatorname{Var}}
\newcommand{\Qprob}{{\mathbb Q}}
\newcommand{\SWG}{\mathrm{SWG}}
\newcommand{\AS}{\mathrm{AS}}
\newcommand{\CE}{\mathrm{C}}
\newcommand{\Art}{\mathrm{Art}}
\newcommand{\FS}{\mathrm{FS}}
\newcommand{\Op}{O_{\prob}}
\newcommand{\ulmn}{\underline{m}_n}
\newcommand{\olmn}{\overline{m}_n}
\newcommand{\hatI}{\hat{I}}
\newcommand{\hatG}{\hat{G}}
\newcommand{\hatS}{\hat{S}}
\newcommand{\Bindep}{B^{(\mathrm{ind})}}
\newcommand{\Sindep}{S^{(\mathrm{ind})}}
\newcommand{\gammaeuler}{\gamma^{(\mathrm{e})}}
\newcommand{\vep}{\varepsilon}
\newcommand{\Wn}{W_n}
\newcommand{\Hn}{H_n}
\def\1{\mathbh{1}}
\newcommand{\indic}[1]{\1_{\{#1\}}}
\newcommand{\indicwo}[1]{\1_{#1}}
\newcommand{\convd}{\stackrel{d}{\longrightarrow}}
\newcommand{\convp}{\stackrel{{\mathbb P}}{\longrightarrow}}
\newcommand{\convas}{\stackrel{\mathrm{a.s.}}{\longrightarrow}}
\makeatother

\begin{document}
\begin{frontmatter}

\title{First passage percolation on random graphs with~finite mean degrees}
\runtitle{First passage percolation on sparse random graphs}

\begin{aug}
\author[A]{\fnms{Shankar} \snm{Bhamidi}\thanksref{t1}\ead[label=e1]{bhamidi@email.unc.edu}},
\author[B]{\fnms{Remco} \snm{van der Hofstad}\corref{}\thanksref{t2}\ead[label=e2]{rhofstad@win.tue.nl}}\\
\and
\author[C]{\fnms{Gerard} \snm{Hooghiemstra}\ead[label=e3]{g.hooghiemstra@tudelft.nl}}
\runauthor{S. Bhamidi, R. van der Hofstad and G. Hooghiemstra}
\affiliation{University of North Carolina,
Eindhoven University of Technology\\ and
Delft University of Technology}
\address[A]{S. Bhamidi\\
Department of Statistics\\
\quad and Operations Research\\
304 Hanes Hall\\
University of North Carolina\\
Chapel Hill,
North Carolina 27599\\
USA\\
\printead{e1}} %adresu isvedimo komanda gale!
\address[B]{R. van der Hofstad\\
Department of Mathematics\\
\quad and
Computer Science\\
Eindhoven University of Technology\\
P.O. Box 513\\
5600 MB Eindhoven\\
The Netherlands\\
\printead{e2}}
\address[C]{G. Hooghiemstra\\
EEMCS\\
Delft University of Technology\\
Mekelweg 4\\
2628 CD Delft\\
The Netherlands\\
\printead{e3}}
\end{aug}
\thankstext{t1}{Supported by NSF Grant DMS-07-04159.}
\thankstext{t2}{Supported in part by Netherlands
Organisation for Scientific Research (NWO).}

% HISTORY:
\received{\smonth{3} \syear{2009}}
\revised{\smonth{10} \syear{2009}}

% ABSTRACT
%
\begin{abstract}
We study first passage percolation on the configuration model.
Assuming that each edge has an independent exponentially distributed
edge weight,
we derive explicit distributional asymptotics for the minimum weight between
two randomly chosen connected vertices in the network, as well as for
the number of edges on the least weight path, the so-called \textit{hopcount}.

We analyze the configuration model with degree power-law exponent
$\tau>2$, in which the degrees are assumed to be i.i.d. with a tail
distribution which is either of power-law form with exponent
$\tau-1>1$, or has even thinner tails ($\tau=\infty$).
In this model, the degrees have a finite first moment,
while the variance is finite for $\tau>3$, but infinite for
$\tau\in(2,3)$.

We prove a central limit theorem for the hopcount, with
asymptotically equal means and variances equal to $\alpha\log{n}$,
where $\alpha\in(0,1)$ for $\tau\in(2,3)$, while $\alpha>1$ for
$\tau>3$. Here $n$ denotes the size of the graph. For $\tau\in
(2,3)$, it is known that the graph distance
between two randomly chosen connected vertices is proportional to $\log
\log{n}$ [\textit{Electron. J. Probab.}
\textbf{12}
(2007)
703--766],
that is, distances are \textit{ultra small}. Thus,
the addition of edge weights causes a marked change in the
geometry of the network. We further study the weight
of the least weight path and prove convergence in
distribution of an appropriately centered version.

This study continues
the program initiated in [\textit{J. Math. Phys.}
\textbf{49} (2008) 125218] of showing that
$\log{n}$ is the correct scaling for the
hopcount under i.i.d. edge disorder, even if the graph distance
between two randomly chosen vertices is of much smaller order. The case
of infinite mean
degrees ($\tau\in[1,2)$) is studied in [Extreme value theory,
{P}oisson--{D}irichlet distributions and
first passage percolation on random networks (2009) Preprint] where it
is proved that
the hopcount remains uniformly bounded and converges in distribution.
\end{abstract}

% KEYWORDS
%
\begin{keyword}[class=AMS]
\kwd{60C05}
\kwd{05C80}
\kwd{90B15}.
\end{keyword}
\begin{keyword}
\kwd{Flows}
\kwd{random graph}
\kwd{first passage percolation}
\kwd{hopcount}
\kwd{central limit theorem}
\kwd{coupling to continuous-time branching processes}
\kwd{universality}.
\end{keyword}

\end{frontmatter}

%s1 ###
\section{Introduction}
\label{sec-int}
The general study of \textit{real-world} networks has seen a
tremendous growth
in the last few years. This growth occurred both at an empirical level of
obtaining data on networks such as the Internet, transportation
networks, such as rail and road networks, and biochemical networks,
such as gene regulatory networks, as well as at a theoretical level in
the understanding of the properties of various mathematical models for these
networks.

We are interested in one specific theoretical aspect of the above vast
and expanding field. The setting
is as follows: Consider a transportation network whose main aim is to transport
flow between various vertices in the network via the available edges.
At the very
basic level there are two crucial elements which affect the flow carrying
capabilities and delays experienced by vertices in the network:

\begin{longlist}[(b)]
\item[(a)] The actual graph topology, such as the density of edges and
existence of
short paths between vertices in the graph distance. In this context
there has been an enormous amount of interest in the concept of
\textit{small-world} networks where the typical graph distance between
vertices in the network is of order $\log{n}$ or even smaller. Indeed,
for many of the mathematical models used to model real-world transmission
networks, such as the Internet, the graph distance can be of order much
smaller than order $\log{n}$. See, for example, \cite
{CohHav03,hofs1}, where for the
configuration model with degree exponent $\tau\in(2,3)$, the
remarkable result that the graph distance between typical vertices is of
order $\log{\log{n}}$ is proved. In this case, we say
that the graph is \textit{ultra small}, a phrase invented in
\cite{CohHav03}. Similar results have appeared
for related models in \cite{ChuLu03,DomHofHoo08,NorRei06}.
The configuration model is described in more detail in Section~\ref{not}.
For introductions to scale-free random graphs, we refer to the monographs
\cite{ChuLu06c,Durr06}, for surveys of classical random graphs
focussing on the Erd\H{o}s--R\'enyi random graph (see \cite
{Boll01,JanLucRuc00}).

\item[(b)] The second factor which plays a crucial role is the edge
weight or
cost structure of the graph which can be thought of as representing
actual economic costs or congestion costs across edges. Edge
weights being identically equal to $1$ gives us back the graph
geometry. What can be said when the edge costs have some other behavior?
The main aim of this study is to understand what happens when each edge
is given an independent edge cost with mean $1$. For simplicity, we have
assumed that the distribution of edge costs is exponentially with mean
$1$ $[\operatorname{Exp}(1)]$, leading to first passage percolation on
the graph involved.
First passage percolation with exponential weights has received
substantial attention (see
\cite
{vcg-random-shanky,hofs-erdos-fpp,hofs-flood,douglas-fpp,smythe-fpp,janson123,wastlund}),
in particular on the complete graph, and, more recently, also on
Erd\H{o}s--R\'enyi random graphs. However, particularly the relation
to the
scale-free nature of the underlying random graph and the behavior of
first passage percolation on it has not yet been investigated.
\end{longlist}

In this paper, we envisage a situation where the edge weights represent
actual economic costs, so that all flow is routed through minimal
weight paths.
The actual time delay experienced by vertices in the network is given by
the number of edges on this least cost path or hopcount $\Hn$.
Thus, for two typical vertices
$1$ and $2$ in the network, it is important to understand both the minimum
weight $\Wn$ of transporting flow between two vertices as well as the
hopcount $\Hn$ or the number\vadjust{\goodbreak} of edges on this minimal weight path.
What we shall
see is the following universal behavior:

\begin{quote}

\textit{Even if the graph topology is of ultra-small nature, the
addition of
random edge weights causes a complete change in the geometry and, in particular,
the number of edges on the minimal weight path between two vertices
increases to
$\Theta(\log{n})$.}

\end{quote}

Here we write $a_n=\Theta(b_n)$ if there exist positive constants
$c$ and $C$, such that, for all $n$, we have $cb_n\le a_n\le Cb_n$.
For the precise mathematical results we refer to Section~\ref{results}.
We shall see that a remarkably universal picture emerges, in the sense
that for each $\tau>2$, the hopcount satisfies a central limit theorem (CLT)
with asymptotically equal mean and variance equal to $\alpha\log{n}$,
where $\alpha\in(0,1)$ for $\tau\in(2,3)$, while $\alpha>1$ for
$\tau>3$. The parameter $\alpha$ is the only feature which is left from
the randomness of the underlying random graph, and $\alpha$ is a
simple function
of~$\tau$ for $\tau\in(2,3)$, and of the average forward degree for
$\tau>3$.
This type of universality is reminiscent of that of simple random
walk, which, appropriately scaled, converges to Brownian motion,
and the parameters needed for the Brownian limit are only the mean and variance
of the step-size. Interestingly, for the Internet hopcount,
measurements show that
the hopcount is close to a normal distribution with equal
mean and variance (see, e.g., \cite{VanHooHof00}), and it would
be of interest to investigate whether first passage percolation
on a random graph can be used as a model for the Internet hopcount.

This paper is part of the program initiated in \cite{vcg-random-shanky}
to rigorously analyze the asymptotics of distances and weights of
shortest-weigh paths
in random graph models under the addition of edge weights. In this
paper, we rigorously
analyze the case of the configuration model with degree exponent $\tau
>2$, the
conceptually important case in practice, since the degree exponent
of a wide variety of real-world networks is conjectured to be in this interval.
In \cite{BhaHofHoo09a}, we investigate the case $\tau\in[1,2)$, where
the first moment of the degrees is infinite and we observe entirely
different behavior of the hopcount $\Hn$.

%s2 ###
\section{Notation and definitions}
\label{not}

We are interested in constructing a random graph on $n$ vertices. Given a
\textit{degree sequence}, namely a sequence of $n$ positive integers
$\bfd= (d_1,d_2,\ldots, d_n)$
with $\sum_{i=1}^n d_i$ assumed to be even, the configuration model
(CM) on $n$
vertices with degree sequence $\bfd$ is constructed as follows:

Start with $n$ vertices and $d_i$ stubs or half-edges
adjacent to vertex $i$. The graph is constructed by randomly pairing
each stub to some other stub to form edges. Let
%
%e2.1 ###
%
\begin{equation} l_n = \sum_{i=1}^n d_i
\end{equation}
denote the total degree. Number the stubs from $1$ to $l_n$ in some
arbitrary order. Then, at each step,
two stubs which are not already paired are chosen uniformly at random among
all the unpaired or \textit{free} stubs and are paired to form a
single edge
in the graph. These stubs are no longer free and removed from the list of
free stubs. We continue with this procedure of choosing and pairing two stubs
until all the stubs are paired. Observe that the order in which we
choose the\vadjust{\goodbreak}
stubs does not matter. Although self-loops may occur,
these become rare as $n\to\infty$ (see, e.g., \cite{Boll01} or
\cite{Jans06b}
for more precise results in this direction).

Above, we have described the construction of
the CM when the degree sequence is given. Here
we shall specify how we construct the actual degree sequence $\bfd$
which shall be \textit{random}. In general, we shall let a capital letter
(such as $D_i$) denote a random variable, while a lower case letter
(such as $d_i$) denote a deterministic object. We shall
assume that the random variables $D_1,D_2,\ldots, D_n$ are
independent and identically distributed (i.i.d.) with a certain
distribution function $F$.
(When the sum of stubs $L_n=\sum_{i=1}^n D_i$ is not even
then we shall use the degree sequence $D_1,D_2,\ldots,D_n$, with $D_n$ replaced
by $D_n+1$. This does not
effect our calculations.)

We shall assume that the degrees of all vertices are at least $2$ and that
the degree distribution $F$ is regularly varying. More precisely, we assume
%
%e2.2 ###
%
\begin{equation} \label{Fcond} \prob(D\geq2) =1\quad\mbox
{and}\quad
1-F(x)=x^{-(\tau-1)} L(x),
\end{equation}
with $\tau>2$, and where $x\mapsto L(x)$ is a slowly varying function
for $x\to\infty$.
In the case $\tau>3$, we shall replace~(\ref{Fcond})
by the less stringent condition~\eqref{distribution>3}. Furthermore, each
edge is given a random edge weight, which in this
study will always be assumed to be independent and identically distributed
(i.i.d.) exponential random variables with mean 1. Because in our
setting the
vertices are exchangeable, we let~$1$ and $2$ be
the two random vertices picked \textit{uniformly at random} in the network.

As stated earlier, the parameter $\tau$ is assumed to satisfy $\tau>2$,
so that the degree distribution has finite mean.
In some cases, we shall distinguish between $\tau>3$ and $\tau\in(2,3)$;
in the former case, the variance of the degrees is finite, while in the
latter, it is infinite.
It follows from the condition $D_i\geq2$, almost surely,
that the probability that the vertices 1 and 2 are connected converges
to 1.

Let $f=\{f_j\}_{j=1}^{\infty}$ denote the probability mass function
corresponding to
the distribution function $F$, so that $f_j=F(j)-F(j-1)$.
Let $\{g_j\}_{j=1}^{\infty}$ denote\vspace*{1pt} the \textit{size-biased} probability
mass function corresponding to $f$, defined by
%
%e2.3 ###
%
\begin{equation} \label{eqn:size-bias} g_j = \frac{(j+1)f_{j+1}}{\mu
},\qquad j\geq0,
\end{equation}
where $\mu$ is the expected size of the degree, that is,
%
%e2.4 ###
%
\begin{equation}
\mu=\expec[D]=\sum_{j=1}^\infty j f_j.
\end{equation}

%s3 ###
\section{Results}
\label{results}

In this section, we state the main results for $\tau>2$.
We treat the case where $\tau>3$ in
Section~\ref{sec-tau>3} and the case where $\tau\in(2,3)$ in
Section~\ref{sec-tau(2,3)}. The case where $\tau\in[1,2)$ is
deferred to \cite{BhaHofHoo09a}.

Throughout the paper, we shall denote by
%
%e3.1 ###
%
\begin{equation} \label{HnWn-def}
(\Hn, \Wn),
\end{equation}
the number of edges and total weight of the shortest-weight
path between vertices~1 and 2 in the CM with i.i.d.\
degrees with distribution function $F$, where we
condition the vertices 1 and 2 to be connected, and we assume that
each edge in the CM has an i.i.d. exponential weight with mean 1.

%s3.1 ###
\subsection{Shortest-weight paths for $\tau>3$}
\label{sec-tau>3}
In this section, we shall assume that the distribution function $F$
of the degrees in the CM is nondegenerate and satisfies
$F(x)=0, x<2$, so that the random variable $D$
is nondegenerate and satisfies $D\ge2$, a.s., and that there exist
$c>0$ and $\tau>3$ such that
%
%e3.2 ###
%
\begin{equation} \label{distribution>3} 1-F(x)\leq c x^{-(\tau
-1)},\qquad x\geq0.
\end{equation}
Also, we let
%
%e3.3 ###
%
\begin{equation} \label{nu-def} \nu=\frac{\expec[D(D-1)]}{\expec
[D]}.
\end{equation}
As a consequence of the conditions we have that $\nu>1$.
The condition $\nu>1$ is equivalent to the
existence of a \textit{giant component} in the CM,
the size of which is proportional to $n$
(see, e.g., \cite{hofs3,MolRee95,MolRee98}; for the
most recent and general result, see \cite{JanLuc07}).
Moreover, the proportionality constant is the survival probability
of the branching process with offspring distribution $\{g_j\}_{j\ge1}$.
As a consequence of the conditions on the distribution function $F$, in
our case, the
survival probability equals $1$, so that for $n\to\infty$ the graph
becomes asymptotically
connected in the sense that the giant component has $n(1-o(1))$ vertices.
Also,\vspace*{-1pt} when~(\ref{distribution>3}) holds, we have
that $\nu<\infty$. Throughout the paper, we shall let $\convd$ denote
convergence in distribution and $\convp$ convergence in probability.

\begin{Theorem}[(Precise asymptotics for $\tau>3$)]
\label{main>3}
Let the degree distribution~$F$ of the
CM on $n$ vertices be nondegenerate, satisfy $F(x)=0, x<2$ and
satisfy~(\ref{distribution>3})
for some $\tau>3$. Then:

\begin{longlist}[(b)]
\item[(a)] the hopcount $\Hn$ satisfies the CLT
%
%e3.4 ###
%
\begin{equation} \label{CLT-hopcount>3} \frac{\Hn-\alpha\log
{n}}{\sqrt{\alpha\log{n}}} \convd Z,
\end{equation}
where $Z$ has a standard normal distribution, and
%
%e3.5 ###
%
\begin{equation} \label{alpha>3} \alpha=\frac{\nu}{\nu-1}\in
(1,\infty);
\end{equation}
\item[(b)] there exists a random variable $V$ such that
%
%e3.6 ###
%
\begin{equation} \Wn-\frac{\log{n}}{\nu-1} \convd V.
\end{equation}
\end{longlist}
\end{Theorem}

In Appendix \hyperref[sec-weak-conv-weight-tau3]{C}, we
shall identify the
limiting random variable $V$ as
%
%e3.7 ###
%
\begin{equation} V=-\frac{\log{W_1}}{\nu-1}- \frac{\log{W_2}}{\nu
-1}+\frac{\Lambda}{\nu-1} +\frac{\log{\mu(\nu-1)}}{\nu-1},
\end{equation}
where $W_1, W_2$ are two independent copies of the limiting random
variable of
a certain supercritical continuous-time branching process, and $\Lambda
$ has a Gumbel
distribution.

%s3.2 ###
\subsection{Analysis of shortest-weight paths for $\tau\in(2,3)$}
\label{sec-tau(2,3)}
In this section, we shall assume that~(\ref{Fcond})
holds for some $\tau\in(2,3)$ and some slowly varying
function $x\mapsto L(x)$. When this is the case, the variance
of the degrees is infinite, while the mean degree is finite.
As a result, we have that $\nu$ in~(\ref{nu-def}) equals
$\nu=\infty$, so that the CM is always supercritical
(see \cite{hofs1,JanLuc07,MolRee95,MolRee98}).
In fact, for $\tau\in(2,3)$, we shall make a stronger assumption on
$F$ than~(\ref{Fcond}), namely, that there exists a $\tau\in(2,3)$ and
$0<c_1\leq c_2<\infty$ such that, for all $x\geq0$,
%
%e3.8 ###
%
\begin{equation} \label{Fcond(2,3)} c_1x^{-(\tau-1)}\leq1-F(x)\leq
c_2x^{-(\tau-1)}.
\end{equation}

\begin{Theorem}[{[Precise asymptotics for $\tau\in(2,3)$]}]
\label{main(2,3)}
Let the degree distribution $F$ of the
CM on $n$ vertices be nondegenerate, satisfy $F(x)=0, x<2$ and
satisfy~\textup{(\ref{Fcond(2,3)})}
for some $\tau\in(2,3)$. Then:

\begin{longlist}[(b)]
\item[(a)] the hopcount $\Hn$ satisfies the CLT
%
%e3.9 ###
%
\begin{equation} \label{CLT-hopcount(2,3)} \frac{\Hn-\alpha\log
{n}}{\sqrt{\alpha\log{n}}} \convd Z,
\end{equation}
where $Z$ has a standard normal distribution and where
%
%e3.10 ###
%
\begin{equation} \label{alpha(2,3)} \alpha=\frac{2(\tau-2)}{\tau
-1}\in(0,1);
\end{equation}
\item[(b)] there exists a limiting random variable $V$ such that
%
%e3.11 ###
%
\begin{equation} \Wn\convd V.
\end{equation}
\end{longlist}
\end{Theorem}

In Section~\ref{sec-lemma:CLT-sum2}, we shall identify the
limiting random variable $V$ precisely as
%
%e3.12 ###
%
\begin{equation} V=V_1+V_2,
\end{equation}
where $V_1, V_2$ are two independent copies of a random variable which
is the explosion time of a certain infinite-mean continuous-time
branching process.

%We now remark on what happens when \refeq{Fcond} holds, but
%In this case, we shall show that the CLT with equal means and
%variances persists to be true, but now
%the asymptotic means and variances are given by $2\log{a_n}$ for some
%$a_n$ which is regularly varying
%with exponent $(\tau-2)/(\tau-1)$.

%s3.3 ###
\subsection{Discussion and related literature}
\label{rlt}

\subsubsection*{Motivation}
The basic motivation of this work was to show that even though the
underlying graph topology might imply that the distance between\vadjust
{\goodbreak} two
vertices is very small, if there are edge weights representing capacities,
say, then the hopcount could drastically increase. Of course, the
assumption of i.i.d. edge weights is
not very realistic; however, it allows us to almost completely analyze
the minimum weight path. The assumption of exponentially distributed
edge weights is probably not necessary \cite{AdaBroLug09,janson123}
but helps in considerably
simplifying the analysis. Interestingly, hopcounts which are
close to normal with asymptotically equal means and variances
are observed in Internet (see, e.g., \cite{VanHooHof00}).
The results presented here might shed some light on the origin of
this observation.

\subsubsection*{Universality for first passage percolation on the CM}
Comparing Theorems~\ref{main>3} and~\ref{main(2,3)} we see that
a remarkably universal picture emerges. Indeed, the hopcount \textit{in
both cases}
satisfies a CLT with equal mean and variance proportional
to $\log{n}$, and the proportionality constant $\alpha$ satisfies
$\alpha\in(0,1)$
for $\tau\in(2,3)$, while $\alpha>1$ for $\tau>3$. We shall see that
the proofs of Theorems~\ref{main>3} and~\ref{main(2,3)} run, to a
large extent,
parallel, and we shall only need to distinguish when dealing with the
related branching process problem to which the neighborhoods can be coupled.

\subsubsection*{The case $\tau\in[1,2)$ and critical cases $\tau=2$ and
$\tau=3$}
In \cite{BhaHofHoo09a}, we study first passage percolation on the CM
when $\tau\in[1,2)$,
that is, the degrees have infinite mean. We show that a remarkably
different picture emerges,
in the sense that $\Hn$ remains uniformly bounded and converges in
distribution.
This is due to the fact that we can think of the CM, when $\tau\in[1,2)$,
as a union of an (essentially) finite number of stars.
%Moreover, the limiting distribution depends on whether we
%multiple edges or not,
%i.e., whether each multiple edge connecting two vertices receives an
%independent weight,
%or the same weight.
Together with the results in Theorems~\ref{main>3}--\ref{main(2,3)},
we see that only the critical cases $\tau=2$ and $\tau=3$ remain
open. We conjecture that the
CLT, with asymptotically equal means and variances, remains valid when
$\tau=3$, but
that the proportionality constant $\alpha$ can take any value in
$[1,\infty)$,
depending on, for example, whether $\nu$ in~(\ref{nu-def}) is finite
or not.
What happens for $\tau=2$ is less clear to~us.

\subsubsection*{Graph distances in the CM}
Expanding neighborhood techniques for random graphs have been
used extensively to explore shortest path structures and other
properties of
locally tree-like graphs. See the closely related papers
\cite{hofs2,hofs3,hofs1,ReiNor04} where an extensive study of the CM
has been carried out. Relevant to our context is
\cite{hofs1}, Corollary 1.4(i), where it has been shown that when $2<
\tau< 3$,
the graph distance ${\widetilde H}_n$ between two typical
vertices,which are conditioned to be
connected, satisfies the asymptotics\vspace*{-0.8pt}
%
%e3.13 ###
%
\begin{equation} \frac{{\widetilde H}_n}{\log{\log{n}}} \convp\frac
{2}{|\log{(\tau-2})|}
\end{equation}
as $n\to\infty$, and furthermore that the fluctuations of
${\tilde H}_n-\log\log{n}$
 remain uniformly bounded as $n\rightarrow\infty$.
For $\tau>3$, it is shown in \cite{hofs3}, Corollary\vspace*{1pt} 1.3(i),
and that ${\tilde H}_n-\log{n}$ has bounded fluctuations\vspace*{-0.8pt}
%
%e3.14 ###
%
\begin{equation} \frac{{\widetilde H}_n}{\log{n}} \convp\frac
{1}{\log{\nu}},
\end{equation}
again with bounded fluctuations. Comparing these results with
Theorems~\ref{main>3}--\ref{main(2,3)}, we see the drastic
effect that the addition of edge weights has on the geometry of the graph.

\subsubsection*{The degree structure} In this paper, as in
\cite{hofs2,hofs3,hofs1,ReiNor04}, we assume that the degrees are
i.i.d. with a certain
degree distribution function $F$. In the literature, also the setting where
the degrees $\{d_i\}_{i=1}^n$ are deterministic, and converge in an appropriate
sense to an asymptotic degree distribution is studied
(see, e.g., \cite{ChuLu03,FerRam04,JanLuc07,MolRee95,MolRee98}).
We expect that our
results can be adapted to this situation. Also, we assume that the degrees
are at least 2 a.s., which ensures that two uniform vertices lie, with
high probability
(w.h.p.) in the giant component. We have chosen for this setting to
keep the
proofs as simple as possible, and we conjecture that Theorems~\ref
{main>3}--\ref{main(2,3)},
when instead we condition the vertices~1 and~2 to be connected, remain
true verbatim
in the more general case of the supercritical CM.

\subsubsection*{Annealed vs. quenched asymptotics} The problem studied in
this paper,
first passage percolation on a random graph, fits in the more general framework
of stochastic processes in random environments, such as random walk in random
environment. In such problems, there are two interesting settings,
namely, when
we study results when averaging out over the environment and when we freeze
the environment (the so-called annealed and quenched asymptotics). In
this paper,
we study the \textit{annealed} setting, and it would be of interest to
extend our
results to the \textit{quenched} setting, that is, study the
first-passage percolation
problem \textit{conditionally on the random graph}. We expect the
results to change
in this case, primarily due to the fact that we know the exact
neighborhood of each point.
However, when we consider the shortest-weight problem between two
\textit
{uniform}
vertices, we conjecture Theorems~\ref{main>3}--\ref{main(2,3)} to
remain valid
verbatim, due to the fact that the neighborhoods of uniform vertices
converge to the
same limit as in the annealed setting (see, e.g., \cite{BerSid09,hofs3}).

\subsubsection*{First passage percolation on the Erd\H{o}s--R\'enyi
random graph}
We recall that the Erd\H{o}s--R\'enyi random graph
$G(n,p)$ is obtained by taking the vertex set $[n]=\{1,\ldots,n\}$
and letting each edge $ij$ be present, independently
of all other edges, with probability $p$.
The study closest in spirit to our study is \cite{vcg-random-shanky} where
similar ideas were explored for dense Erd\H{o}s--R\'enyi random graphs.
The Erd\H{o}s--R\'enyi random graph
$G(n,p)$ can be viewed as a close brother of the CM, with Poisson
degrees, hence with $\tau=\infty$.
Consider the case where $p=\mu/n$ and $\mu>1$.
In a future paper we plan to show, parallel to
the above analysis, that $\Hn$ satisfies a CLT
with asymptotically equal mean and variance given by
$\frac{\mu}{\mu-1}\log{n}$. This connects up
nicely with \cite{vcg-random-shanky}
where related results were shown for $\mu=\mu_n\rightarrow\infty,$
and $\Hn/\log{n}$ was proved to converge to 1 in probability.
See also \cite{hofs-erdos-fpp} where related statements
were proved under stronger assumptions on $\mu_n$.
Interestingly, in a recent paper, Ding et al. \cite{DinKimLubPer09}
use first passage percolation to study the diameter of the largest
component of
the Erd\H{o}s--R\'enyi random graph with
edge probability $p=(1+\vep)/n$ for $\vep=o(1)$ and $\vep^3n\to
\infty$.

\subsubsection*{The weight distribution}
It would be of interest to study the effect of
weights even further, for example, by studying the case where
the weights are i.i.d. random variables with distribution equal to $E^{s}$
where $E$ is an exponential random variable with mean 1 and
$s\in[0,\infty)$. The case $s=0$ corresponds to the graph
distance~${\tilde H}_n$ as
studied in \cite{hofs2,hofs3,hofs1}
while the case $s=1$ corresponds to the case
with i.i.d. exponential weights as studied here. Even the problem on
the complete
graph seems to be open in this case, and we intend to return to this problem
in a future paper. We conjecture that the CLT remains valid for
first passage perolation on the CM when the weights are given by
independent copies of
$E^s$ with asymptotic mean and variance proportional to $\log{n}$,
but, when $s\neq1$, we predict that the asymptotic means and variances
have \textit{different} constants.

We became interested in random graphs with edge weights from \cite
{weak-strong-diso}
where, via empirical simulations, a wide variety of behavior
was predicted for the shortest-weight paths in various random graph models.
The setup that we analyze is the \textit{weak disorder} case. In \cite
{weak-strong-diso}, also
a number of interesting conjectures regarding the \textit{strong
disorder case} were made, which would correspond to analyzing the
minimal spanning
tree of these random graph models, and which is a highly interesting problem.

\subsubsection*{Related literature on shortest-weight problems}
First passage percolation, especially on the integer lattice, has been
extensively
studied in the last fifty years (see, e.g., \cite{ham-wel,smythe-fpp} and the
more recent survey
\cite{douglas-fpp}). In these papers, of course, the emphasis is
completely different, in the sense
that geometry plays an intrinsic role and often the goal of the study
is to show that
there is a limiting ``shape'' to first passage percolation from the origin.

Janson \cite{janson123} studies first passage percolation on the
complete graph with
exponential weights. His main results are
%
%e3.15 ###
%
\begin{equation}
\label{janson}\qquad\quad
\frac{W_n^{(ij)}}{\log n/n}
\convp1,\qquad
\frac{\max_{j\le n}W_n^{(ij)}}{\log n/n}
\convp2,\qquad
\frac{\max_{i,j\le n}W_n^{(ij)}}{\log n/n}
\convp3,
\end{equation}
where $W_n^{(ij)}$ denotes the weight of the shortest path between
the vertices $i$ and~$j$. Recently the authors of \cite{AdaBroLug09}
showed in the same set-up that
$\max_{i,j\le n} H_n^{(ij)}/\break\log n \convp\alpha^\star$ where
$\alpha^\star\approx3.5911$ is the unique solution of the equation
$x \log{x}-x=1$.
It would be of interest to investigate such questions in the CM with
exponential weights.

The fundamental difference of first passage percolation on the integer lattice,
or even on the complete graph, is that in our case the underlying graph is
random as well, and we are lead to the delicate relation between the randomness
of the graph together with that of the stochastic process, in this
case first passage percolation, living on it.
Finally, for a slightly different perspective to shortest weight
problems, see \cite{wastlund} where
relations between the random assignment problem and the shortest-weight problem
with exponential edge weights on the complete graph are explored.

%%%%%%%%%%%%%%%%%%%%%%%%%%%%%%%%%%%%%%%%%%%%%%%%%%%%%%%%%%%%%%%%%%%%%%
%%%%%%%%%%%%%%%%%%%%%%%%%%%%%%%%%%%%%%%%%%%%%%%%%%%%%%%%%%%%%%%%%%%%%%
%s4 ###
\section{Overview of the proof and organization of the paper}
\label{sec-overview}
The key idea of the proof is to first grow the shortest-weight
graph (SWG) from vertex 1, until it reaches an appropriate size.
After this, we grow the SWG from vertex 2 until it connects up with the
SWG from vertex 1. The size to which we let the SWG from~1 grow shall be
the same as the \textit{typical size} at which the connection
between the SWG from vertices 1 and 2 shall be made. However, the
connection time
at which the SWG from vertex 2 connects to the SWG from vertex 1
is \textit{random}.

More precisely, we define the SWG from vertex 1, denoted by $\SWG
^{(1)}$, recursively.
The growth of the SWG from vertex 2, which is denoted by $\SWG^{
(2)}$, is similar.
We start with vertex 1 by defining $\SWG_0^{(1)}=\{1\}$. Then we
add the edge and vertex with minimal
edge weight connecting vertex $1$ to one of its neighbors (or itself
when the minimal edge is a self-loop).
This defines $\SWG_{1}^{(1)}$. We obtain $\SWG_m^{(1)}$
from $\SWG_{m-1}^{(1)}$
by adding the edge and end vertex connected to the
$\SWG_{m-1}^{(1)}$
with minimal edge weight. We informally let $\SWG
_m^{(i)}$
denote the SWG from vertex $i\in\{1,2\}$ when $m$ edges (and vertices)
have been added
to it. This definition is \textit{informal}, as we shall need to deal with
self-loops and cycles in a proper way. How we do this is explained in
more detail
in Section~\ref{sec-coupling}. As mentioned before, we first grow
$\SWG_{m}^{(1)}$ to a size $a_n$, which is to be chosen appropriately.
After this, we grow $\SWG_{m}^{(2)}$, and we stop as soon as a
vertex of
$\SWG_{a_n}^{(1)}$ appears in $\{\SWG_m^{(2)}\}
_{m=0}^{\infty}$,
as then the shortest-weight path between vertices 1 and
2 has been found. Indeed, if on the contrary, the shortest weight path
between vertex 1 and vertex 2
contains an edge not contained in the union of the two SWGs when they
meet, then necessarily
this edge would have been chosen in one of the two SWGs at an earlier
stage, since at some
earlier stage this edge must have been incident to one of the SWGs and
had the minimal
weight of all edges incident to that SWG.
In Sections~\ref{sec-coupling} and~\ref{sec-two_flows},
we shall make these definitions precise.

Denote this first common vertex by $A$, and
let $G_i$ be the distance between vertex~$i$ and $A$, that is, the
number of edges
on the minimum weight path from $i$ to $A$.
Then we have that
%
%e4.1 ###
%
\begin{equation} \label{Hn-sum1} \Hn=G_1+G_2,
\end{equation}
while, denoting by $T_i$ the weight of the shortest-weight paths
from $i$ to $A$,
we have
%
%e4.2 ###
%
\begin{equation} \Wn=T_1+T_2.
\end{equation}
Thus, to understand the random variables $\Hn$ and $\Wn$,
it is paramount to understand the random
variables $T_i$ and $G_i$, for $i=1,2$.

Since, for $n\to\infty$, the
topologies of the neighborhoods of vertices 1 and 2
are close to being independent, it seems likely that
$G_1$ and $G_2$, as well as $T_1$ and $T_2$
are close to independent. Since, further,
the CM is locally tree-like, we are lead to
the study of the problem on a tree.

With the above in mind, the paper is organized as follows:

In Section~\ref{sec-flow_tree} we study the flow on a tree. More precisely,
in Proposition~\ref{lemma:CLT-sum1}, we describe
the asymptotic distribution of the length
and weight of the shortest-weight path between
the root and the $m$th added vertex in a
branching process with i.i.d. degrees with offspring
distribution $g$ in~(\ref{eqn:size-bias}). Clearly, the CM
\textit{has} cycles and self-loops, and  thus sometimes deviates from
the tree description.

 In Section~\ref{sec-coupling},
we reformulate the problem of the growth of the SWG from a fixed vertex
as a problem of the SWG on a tree,
where we find a way to deal with cycles by a coupling
argument, so that the arguments in Section~\ref{sec-flow_tree}
apply quite literally. In Proposition~\ref{lemma:CLT-sum2},
we describe the asymptotic distribution of the length
and weight of the shortest-weight path between
a fixed vertex and the $m$th added vertex in the SWG from the CM.
However, observe that the random variables $G_i$ described above are the
generation of a vertex at the time at which the two SWGs collide,
and this time is a \textit{random} variable.

 In Section~\ref{sec-two_flows}, we extend the
discussion to this setting and, in Section~\ref{sec-connedge}, we formulate the necessary
ingredients for the collision time, that is, the time at which
the connecting edge appears, in Proposition~\ref{sec-connedge}.
In Section~\ref{sec-complpf}, we complete the outline.

The proofs of the key propositions
are deferred to Sections~\ref{sec-lemma:CLT-sum1}--\ref{sec-lemma:conn_edge}.

 Technical results needed in the proofs in Sections
\ref{sec-lemma:CLT-sum1}--\ref{sec-lemma:conn_edge}, for example on the
topology of the CM, are deferred to the Appendix~\hyperref[sec-app-A]{A}.

%s4.1 ###
\subsection{Description of the flow clusters in trees}
\label{sec-flow_tree}
We shall now describe the construction of the SWG in
the context of trees. In particular, below, we shall deal with a
flow on a branching process tree, where the offspring is deterministic.

\textit{Deterministic construction}: Suppose we have positive (nonrandom)
integers $d_1, d_2,\ldots.$
Consider the following construction of a branching process
in discrete time:
\begin{Construction}[(Flow from root of tree)]
\label{constr:det}
The shortest-weight graph on a tree with degrees $\{d_i\}_{i=1}^{\infty}$
is obtained as follows:
\begin{enumerate}
\item At time $0$, start with one alive vertex (the initial ancestor).
\item At each time step $i$, pick one of the alive vertices at random,
this vertex dies giving birth to $d_i$ children.
\end{enumerate}
\end{Construction}

In the above construction, the number of offspring $d_i$
is fixed once and for all. For a branching process tree,
the variables $d_i$ are i.i.d. \textit{random} variables. This
case shall be investigated later on, but the case
of deterministic degrees is more general and shall
be important for us to be able to deal with the CM.

Consider a continuous-time branching process defined as follows:
\begin{enumerate}
\item Start with the root which dies immediately giving rise to
$d_1$ alive offspring.
\item Each alive offspring lives for $\operatorname{Exp}(1)$ amount of
time, independent
of all other
randomness involved.
\item When the $m$th vertex dies it leaves behind $d_m$ alive offspring.
\end{enumerate}

\setcounter{footnote}{2}
The split-times (or death-times) of this branching process are denoted by
$T_i, i\ge1$. Note that the Construction~\ref{constr:det} is
equivalent to this continuous branching process, observed
at the discrete times $T_i, i\ge1$. The fact that the chosen alive
vertex is chosen at random follows
from the memoryless property of the exponential random variables that
compete to become the minimal one.
We quote a fundamental result from
\cite{buhler}. In its statement, we let
%
%e4.3 ###
%
\begin{equation} \label{si-def} s_i = d_1+\cdots+d_i -
(i-1).\footnotemark
\end{equation}
\footnotetext{A new \textit{probabilistic} proof is added, since
there is
some confusion between the definition $s_i$ given here, and the
definition of $s_i$ given in \cite{buhler}, below equation (3.1). More
precisely, in \cite{buhler}, $s_i$ is defined as $s_i=d_1+\cdots
+d_i-i$, which is our $s_i-1$.}%

\vspace*{-18pt}
\begin{Proposition}[(Shortest-weight paths on a tree)]
\label{prop:gen}
Pick an alive vertex at time $m\geq1$ uniformly at random among all
vertices alive at this
time. Then
\begin{longlist}[(b)]
\item[(a)] the generation of the $m$th chosen vertex is equal in
distribution to
%
%e4.4 ###
%
\begin{equation} \label{Gm-def} G_m \stackrel{d}{=} \sum_{i=1}^m I_i,
\end{equation}
where $\{I_i\}_{i=1}^{\infty}$ are independent Bernoulli random
variables with
%
%e4.5 ###
%
\begin{equation} \label{bi-def} \prob(I_i = 1) = d_i/s_i;
\end{equation}
\item[(b)] the weight of the shortest-weight path between the root
of the tree and the vertex chosen in the $m$th step is equal
in distribution to
%
%e4.6 ###
%
\begin{equation} \label{Tm-def} T_m \stackrel{d}{=} \sum_{i=1}^m
E_i/s_i,
\end{equation}
where $\{E_i\}_{i=1}^{\infty}$ are i.i.d. exponential random
variables with mean 1.
\end{longlist}
\end{Proposition}

\begin{pf} We shall prove part (a) by induction. The statement is trivial
for \mbox{$m=1$}. We next
assume that~\eqref{Gm-def} holds for $m$ where
$\{I_i\}_{i=1}^{m}$ are independent Bernoulli random
variables satisfying~\eqref{bi-def}.
%Let $(U_1,U_2)$ be a sample of size $2$, taken without
%replacement from all alive vertices at time $m$. Observe that, by
%exchangeability, $U_1$ and $U_2$ have the same marginal distribution.
%Moreover, by the induction
%hypothesis the generation of vertex $U_1$, which we denote by
%$G_{U_1}$ is the sum of $m$ independent Bernoulli random variables
%$I_1,I_2,\ldots,I_m$,
%satisfying \eqref{bi-def}. By the distributional equality,
%the same statement holds for the generation $G_{U_2}$ of $U_2$.
Let $G_{m+1}$ denote the generation of the randomly chosen vertex at
time $m+1$, and consider the
event $\{G_{m+1}=k\}, 1\le k \le m$.
If randomly choosing one of the alive vertices at time $m+1$ results in
one of the
$d_{m+1}$ newly added vertices, then, in order to obtain generation
$k$, the previous uniform
choice, that is, the choice of the vertex which was the last one to
die, must have been a vertex from generation $k-1$.
On the other hand, if a uniform pick is conditioned on not taking
one of the $d_{m+1}$ newly added vertices, then this choice must have been
a uniform vertex from generation $k$. Hence,
we obtain, for $1\le k \le m,$
%
%e4.7 ###
%
\begin{equation} \label{Gm-induc} \prob(G_{m+1}=k)=\frac
{d_{m+1}}{s_{m+1}} \prob(G_{m}=k-1)+ \biggl( 1-\frac{d_{m+1}}{s_{m+1}}
\biggr)\prob(G_{m}=k).
\end{equation}
The proof of part (a) is now immediate from the induction hypothesis.
%Next, recall that $G_{U_1},G_{U_2}$ have the same marginal
%distributions.
%Then, \refeq{Gm-induc} proves that for $1\le k \le m,$
% \begin{equation}
% \label{indstap}
% \prob(G_{m+1}=k)=\prob(G_m+I_{m+1}=k)=\prob(\sum_{i=1}^{m+1}I_i=k).
% \end{equation}
%The statement for $\prob(G_{m+1}=m+1)=\prob(\sum_{i=1}^{m+1}I_i=m+1)$
%follows by taking
%the complement of $\prob(G_{m+1}=m+1)$ and using \eqref{indstap}.
%
The proof of part (b) is as follows.
The minimum of $s_i$ independent $\exp(1)$ random variables has an exponential
distribution with parameter $s_i$, and is hence equal in distribution
to $E_i/s_i$.
We further use the memoryless property of the exponential distribution
which guarantees
that at each of the discrete time steps the remaining lifetimes (or
weights) of the
alive vertices are independent exponential variables with mean 1, independent
of what happened previously.
\end{pf}

We note that, while Proposition~\ref{prop:gen} was applied in
\cite{buhler}, Theorem 3.1, only in the case where the degrees are i.i.d.,
in fact, the results hold more generally for \textit{every} tree
(see, e.g., \cite{buhler}, equation (3.1), and the above proof). This
extension shall
prove to be vital in our analysis.

We next intuitively relate the above result to our setting. Start from
vertex 1,
and iteratively choose the edge with minimal additional weight
attached to the SWG so far.
As mentioned before, because of the properties of the exponential
distribution, the edge with minimal additional weight
can be considered to be picked uniformly at random from all edges
attached to the SWG at that moment.
With high probability,
this edge is connected to a vertex which is not
in the SWG. Let $B_i$ denote the \textit{forward degree} (i.e., the
degree minus 1) of the
vertex to which the $i$th edge is connected. By the results in \cite
{hofs3,hofs1},
$\{B_i\}_{i\geq2}$ are close to being i.i.d. and have
distribution given by~(\ref{eqn:size-bias}). Therefore, we are lead to studying
random variables of the form~(\ref{Gm-def})--(\ref{bi-def}) where $\{
B_i\}_{i=1}^{\infty}$
are i.i.d. random variables. Thus, this means that
we study the \textit{unconditional} law of $G_m$ in~(\ref{Gm-def}),
in the setting where the vector $\{d_i\}_{i=1}^\infty$
is replaced by an i.i.d. sequence of random variables
$\{B_i\}_{i=1}^\infty$. We shall first state a CLT for
$G_m$ and a limit result for $T_m$ in this setting. %In\vspace*{1pt} its statement,
%we shall also
%make use of the random variable $\widetilde T_m$, which is the weight
%of the shortest weight path
%between the root and the \textit{parent} of the $m$th individual in the
%branching process.
%Thus, in particular, $\widetilde T_m\leq T_m,$ and $T_m-\widetilde T_m$
%is the time between
%the addition of the parent of the $m$th individual and the $m$th
%individual itself.

\begin{Proposition}[(Asymptotics for shortest-weight paths on trees)]
\label{lemma:CLT-sum1}
Let $\{B_i\}_{i=1}^{\infty}$ be an i.i.d. sequence of nondegenerate,
positive integer valued,
random variables satisfying
\[
\prob(B_i>k)=k^{2-\tau}L(k),\qquad\tau>2,
\]
for some slowly varying function $k\mapsto L(k)$.
Denote by $\nu=\expec[B_1]$, for $\tau>3$, whereas $\nu=\infty$,
for $\tau\in(2,3)$.
Then,
\begin{longlist}[(b)]
\item[(a)] for $G_m$ given in
\textup{(\ref{Gm-def})}--\textup{(\ref{bi-def})}, with $d_i$
replaced by $B_i$,
there exists a $\beta\geq1$ such that, as $m\rightarrow\infty$,
%
%e4.8 ###
%
\begin{equation} \frac{G_m-\beta\log{m}}{\sqrt{\beta\log m}}\convd
Z, \qquad\mbox{where } Z\sim{\cal N}(0,1)
\end{equation}
a standard normal variable, and where $\beta=\nu/(\nu-1)$ for $\tau>3$,
while $\beta=1$ for $\tau\in(2,3)$;
\item[(b)] for $T_m$ given in \textup{(\ref{Tm-def})}, there exists a random
variable  $X$ such that
%
%e4.9 ###
%
\begin{equation} \label{Tm-conv} T_m-\gamma\log{m}\convd X,
\end{equation}
where $\gamma=1/(\nu-1)$ when $\tau>3$, while $\gamma=0$ when
$\tau\in(2,3)$.
\end{longlist}
\end{Proposition}

Proposition~\ref{lemma:CLT-sum1} is proved in \cite{buhler}, Theorem
3.1,
when $\Var(B_i)<\infty$, which holds when $\tau>4$, but not when
$\tau\in(2,4)$.
We shall prove Proposition~\ref{lemma:CLT-sum1} in Section~\ref
{sec-lemma:CLT-sum1} below.
There, we shall also see that the result persists under weaker
assumptions than
$\{B_i\}_{i=1}^\infty$ being i.i.d., for example, when $\{B_i\}
_{i=1}^\infty$
are \textit{exchangeable} nonnegative integer valued random variables
satisfying certain
conditions. Such extensions shall prove to be useful when dealing with
the actual (forward) degrees in the CM.

%s4.2 ###
\subsection{A comparison of the flow on the CM and the flow on the tree}
\label{sec-coupling}
Proposition~\ref{lemma:CLT-sum1} gives a CLT for
the generation when considering a flow on a tree. In this section, we shall
relate the problem of the flow on the CM to the flow on a tree.
The key feature of this construction is that
\textit{we shall simultaneously grow the graph topology neighborhood
of a vertex,
as well as the shortest-weight graph from it.}
This will be achieved by combining the construction of the CM
as described in Section~\ref{not} with the fact that,
from a given set of vertices and edges, if we grow the shortest-weight
graph, each potential edge is equally likely to be the minimal one.

In the problem of finding the shortest weight path between two
vertices~1 and~2,
we shall grow two SWGs simultaneously from the two vertices
1 and 2, until they meet. This is the problem that we actually need to resolve
in order to prove our main results in Theorems~\ref{main>3}--\ref{main(2,3)}.
The extension to the growth of two SWGs is treated in Section~\ref
{sec-two_flows}
below.

The main difference between the flow on a graph and on a tree is that
on the tree there are no cycles, while on a graph there are.
Thus we shall adapt the growth of the SWG for the CM in such a way
that we obtain a tree (so that the results from Section~\ref{sec-flow_tree}
apply) while we can still retrieve all information about shortest-weight
paths from the constructed graph. This will be achieved by
introducing the notion of \textit{artificial} vertices and stubs.
We start by introducing some notation.

We denote by $\{\SWG_m\}_{m\geq0}$ the SWG process from vertex 1.
We construct this process recursively. We let $\SWG_0$ consist only of
the alive vertex 1, and we let $S_0=1$. We next let $\SWG_1$
consist of the $D_1$ allowed stubs and of the explored vertex 1,
and we let $S_1=S_0+D_1-1=D_1$ denote the number of
allowed stubs. In the sequel of the construction, the allowed
stubs correspond to vertices in the shortest-weight
problem on the tree in Section~\ref{sec-flow_tree}.
This constructs $\SWG_1$. Next, we describe how to construct
$\SWG_m$ from $\SWG_{m-1}$. For this construction, we shall
have to deal with several types of stubs:

\begin{longlist}[(b)]
\item[(a)] The allowed stubs at time $m$, denoted by $\AS_{m}$,
are the stubs that are incident to vertices of the $\SWG_{m}$ and that
have not yet been paired to form an edge; $S_m=|\AS_m|$
denotes their number;
\item[(b)] the free stubs at time $m$, denoted by $\FS_{m}$, are
those stubs
of the
$L_n$ total stubs which have not yet been
paired in the construction of the CM up to and including time $m$;
\item[(c)] the artificial stubs at time $m$, denoted by $\Art_m$, are the
\textit{artificial}
stubs created by breaking ties, as described in more detail below.
\end{longlist}

We note that $\Art_m\subset\AS_m$, indeed, $\AS_m\setminus\FS
_{m}=\Art_m.$
Then, we can construct $\SWG_m$ from $\SWG_{m-1}$ as follows.
We choose one of the $S_{m-1}$ allowed stubs uniformly at random, and then,
if the stub is not artificial, pair it uniformly at random to a free
stub unequal to itself.
Below, we shall consistently call these two stubs the
\textit{chosen} stub and the \textit{paired} stub, respectively.
There are 3 possibilities, depending on what kind of stub we choose
and what kind of stub it is paired to:

\begin{Construction}[(The evolution of SWG for CM as SWG on a tree)]
\label{const-SWGCM}
\textup{(1)}~The chosen stub is real, that is, not artificial, and the
paired stub is not one of the allowed stubs.
In this case, which shall be most likely at the
start of the growth procedure of the SWG, the paired stub
is incident to a vertex outside $\SWG_{m-1}$, we
denote by $B_m$ the forward degree of the vertex incident to the paired stub
(i.e., its degree minus 1) and we define $S_m=S_{m-1}+B_m-1$. Then we remove
the paired and the chosen stub from $\AS_{m-1}$ and add the
$B_m$ stubs incident
to the vertex incident to the paired stub to $\AS_{m-1}$ to
obtain $\AS_m$, we remove the chosen and the paired stubs from $\FS_{m-1}$
to obtain $\FS_{m}$, and $\Art_{m}=\Art_{m-1}$.

\textup{(2)} The chosen stub is real and the paired stub is an
allowed stub. In this case, the paired stub is incident to a
vertex in $\SWG_{m-1}$, and we have created a cycle.
In this case, we create an artificial stub
replacing the paired stub and denote $B_m=0$. Then
we let $S_m=S_{m-1}-1$, remove both the chosen and paired stubs from
$\AS_{m-1}$ and add the artificial stub to obtain $\AS_m$, and
remove the chosen and paired stub from
$\FS_{m-1}$ to obtain $\FS_{m}$, while
$\Art_{m}$ is $\Art_{m-1}$ together with the newly created artificial stub.
In $\SWG_m$, we also add an artificial edge to an artificial vertex in
the place where the chosen stub was, the forward degree of the
artificial vertex being~0.
This is done because a vertex is added each time in the construction on
a tree.%\looseness=1

\textup{(3)} The chosen stub is artificial. In this case, we let $B_m=0$,
$S_m=S_{m-1}-1$ and
remove the chosen stub from $\AS_{m-1}$ and $\Art_{m-1}$ to obtain
$\AS_{m}$ and $\Art_m$, while $\FS_{m}=\FS_{m-1}$.
\end{Construction}

In   Construction~\ref{const-SWGCM}, we always work
on a
tree since we replace an edge which creates a cycle, by
one artificial stub, to replace the paired stub, and an artificial edge
plus an artificial
vertex in the $\SWG_m$ with degree 0, to replace the chosen stub. Note
that the number of allowed
edges at time $m$ satisfies $S_m=S_{m-1}+B_m-1$, where $B_1=D_1$ and,
for $m\geq2$, in cases (2) and (3),
$B_m=0$, while in case (1) (which we expect to occur in most cases),
the distribution of $B_m$ is equal to the forward degree of a vertex
incident to a uniformly chosen stub. Here, the choice of stubs is
without replacement.

The reason for replacing cycles as described above is that we wish to represent
the SWG problem as a problem on a tree, as we now will explain
informally. On a tree with degrees $\{d_i\}_{i=1}^{\infty}$, as in
Section~\ref{sec-flow_tree}, we have that the remaining degree of
vertex $i$
at time $m$ is precisely equal to
$d_i$ minus the number of neighbors that are among the
$m$ vertices with minimal shortest-weight paths
from the root. For first passage percolation on a graph with cycles,
a cycle does not only remove one of the edges of
the vertex incident to it (as on the tree), but also
one edge of the vertex at the other end of the cycle. Thus
this is a \textit{different} problem, and the results from
Section~\ref{sec-flow_tree} do not apply literally.
By adding the artificial stub, edge and vertex, we artificially keep the
degree of the receiving vertex the same, so that we \textit{do}
have the same situation as on a tree, and we can use the results
in Section~\ref{sec-flow_tree}. However, we do need to investigate
the relation between the problem with the artificial stubs
and the original SWG problem on the CM. That is the content of the next
proposition.

In its statement, we shall define the $m$th closest vertex
to vertex 1 in the CM, with i.i.d. exponential weights, as the
unique vertex of which the minimal weight path is the $m$th smallest
among all $n-1$ vertices. Further, at each time $m$, we denote by
\textit{artificial vertices} those vertices which are artificially
created, and we call the other vertices \textit{real vertices}.
Then we let the random time $R_m$ be the first time $j$ that $\SWG_{j}$
consists of $m+1$ real vertices, that is,
%
%e4.10 ###
%
\begin{equation} \label{Rm-def} R_m=\min\{j\geq0 \dvtx\SWG_{j}
\mbox{ contains $m+1$ \mbox{real} vertices} \}.
\end{equation}
The $+1$ originates from the
fact that at time $m=0$, $\SWG_{0}$ consists if 1 real vertex,
namely, the vertex from which we construct the SWG.
Thus, in the above set up, we have that $R_m=m$ precisely
when no cycle has been created in the construction up to time $m$.
Then our main coupling result is as follows:

\begin{Proposition}[(Coupling shortest-weight graphs on a tree and CM)]
\label{prop-coupling}
Jointly for all $m\geq1$, the set of real vertices in $\SWG_{R_m}$ is equal
in distribution to the set of $i$th closest vertices
to vertex 1, for $i=1, \ldots, m$. Consequently:

\begin{longlist}[(b)]
\item[(a)] the generation of the $m$th closest vertex to vertex 1
has distribution $G_{R_m}$ where $G_m$ is defined
in~\textup{(\ref{Gm-def})}--\textup{(\ref{bi-def})} with $d_1=D_1$
and $d_i=B_i,
i\geq2,$ as described
in Construction~\ref{const-SWGCM};
\item[(b)] the weight of the shortest weight path to the $m${th} closest
vertex to vertex 1
has distribution $T_{R_m}$, where $T_m$ is defined in~\textup{(\ref{Tm-def})}
with $d_1=D_1$ and $d_i=B_i, i\geq2,$ as described in
Construction~\textup{\ref{const-SWGCM}}.
\end{longlist}
\end{Proposition}

We shall make use of the nice property that the sequence
$\{B_{R_m}\}_{m=2}^{n}$, which consists of the forward degrees
of chosen stubs that are paired to stubs which are not in the
SWG, is, for the CM, an exchangeable sequence of
random variables (see Lemma~\ref{lem-exchange}
below). This is due to the fact that a free stub is
chosen uniformly at random, and the order of the choices
does not matter. This exchangeability shall prove to be useful in
order to investigate shortest-weight paths in the CM. We now
prove Proposition~\ref{prop-coupling}.

\begin{pf*}{Proof of Proposition~\ref{prop-coupling}}
In growing the
SWG, we give exponential weights to the set $\{\AS_m\}_{m\geq1}$.
After pairing, we identify the exponential weight of
the chosen stub to the exponential weight of the edge which it is part of.
We note that by the memoryless property of the exponential random variable,
each stub is chosen uniformly at random from all the allowed
stubs incident to the SWG at the given time. Further,
by the construction of the CM in Section~\ref{not}, this stub is
paired uniformly at random to one of the available free stubs.
Thus the growth rules of the SWG in Construction~\ref{const-SWGCM}
equal those
in the above description of $\{\SWG_m\}_{m=0}^{\infty}$, unless a
cycle is closed and an artificial stub, edge and vertex are created.
In this case, the artificial stub, edge and vertex
might influence the law of the SWG. However, we note that
the artificial vertices are not being counted in the set of
real vertices, and since artificial vertices have forward degree
$0$, they will not be a part of any shortest path to a real vertex.
Thus the artificial vertex at the end of the artificial edge
does not affect the law of the SWG.
Artificial stubs that are created to replace paired stubs when
a cycle is formed, and which are not yet removed at time $m$, will
be called \textit{dangling ends}.
Now, if we only consider real vertices, then the distribution of
weights and
lengths of the shortest-weight paths between the starting points and those
real vertices are identical. Indeed, we can decorate any graph with as many
dangling ends as we like without changing the shortest-weight
paths to real vertices in the graph.
\end{pf*}

Now that the flow problem on the CM has been translated into
a flow problem on a related tree of which we have explicitly described
its distribution, we may make use of
Proposition~\ref{prop:gen} which shall allow us to extend
Proposition~\ref{lemma:CLT-sum1} to the setting of the CM.
Note that, among others, due to the fact that when we draw
an artificial stub, the degrees are not i.i.d. (and not even exchangeable
since the probability of drawing an artificial stub is likely to
increase in time), we need to extend Proposition~\ref{lemma:CLT-sum1}
to a setting where the degrees are weakly dependent.
In the statement of the result, we recall that $G_m$ is
the height of the $m$th added vertex in the tree problem
above. In the statement below, we write
%
%e4.11 ###
%
\begin{equation} \label{an-def}
a_n=n^{(\tau\wedge3-2)/(\tau\wedge
3-1)} =
\cases{ n^{(\tau-2)/(\tau-1)} &\quad for $\tau\in(2,3)$,\vspace
*{2pt}\cr
n^{1/2}& \quad for $\tau>3$,
}
\end{equation}
where, for $a,b\in{\mathbb R}$, we write $a\wedge b=\min\{a,b\}$.

Before we formulate the CLT
for the hopcount of the shortest-weight graph in the CM, we repeat once more
the setup of the random variables involved. Let
$S_0=1$, $S_1=D_1$, and for $j\ge2$,
%
%e4.12 ###
%
\begin{equation}
\label{setup-6.4}
S_j=D_1+\sum_{i=2}^j (B_i-1),
\end{equation}
where, in case the chosen stub is real, that is, not artificial, and the
paired stub is not one of the allowed stubs, $B_i$ equals the forward degree
of the vertex incident to the $i$th paired stub, whereas
$B_i=0$ otherwise. Finally, we recall that, conditionally on
$D_1,B_2,B_3,\ldots,B_m$,
\begin{eqnarray}
\label{def-Gm1}
G_m&=&\sum_{i=1}^m I_i \qquad\mbox{where}\nonumber\\[-8pt]\\[-8pt]
\prob(I_1=1)&=&1,\qquad
\prob(I_j=1)=B_j/S_j,\qquad2 \leq j \leq m.\nonumber
\end{eqnarray}

\begin{Proposition}[(Asymptotics for shortest-weight paths in the CM)]
\label{lemma:CLT-sum2}
\textup{(a)} Let the law of $G_m$ be given in
\textup{(\ref{def-Gm1})}. Then,
with $\beta\geq1$ as in Proposition \textup{\ref{lemma:CLT-sum1}},
and as long as $m\leq\olmn,$ for any $\olmn$ such that $\log{(\olmn
/a_n)}=o(\sqrt{\log{n}})$,
%
%e4.13 ###
%
\begin{equation} \frac{G_m-\beta\log{m}}{\sqrt{\beta\log m}}\convd
Z \qquad\mbox{where } Z\sim{\cal N}(0,1).
\end{equation}

\textup{(b)} Let the law of $T_m$ be given in \eqref{Tm-def} with $s_i$
replaced by $S_i$ given by \eqref{setup-6.4}, and let $\gamma$ be
as in Proposition~\ref{lemma:CLT-sum1}. Then
there exists a random variable $X$ such that
%
%e4.14 ###
%
\begin{equation} \label{Tm-conv2} T_m-\gamma\log{m}\convd X.
\end{equation}
The same results apply to $G_{R_m}$ and $T_{R_m}$,
that is, in the statements \textup{(a)} and \textup{(b)} the integer
$m$ can be replaced by
$R_m$, as long as $m\leq\olmn$.
\end{Proposition}

Proposition~\ref{lemma:CLT-sum2} implies that the result of
Proposition~\ref{lemma:CLT-sum1} remains true for the
CM whenever $m$ is not too large.
Important for the proof of Proposition~\ref{lemma:CLT-sum2} is the
coupling to a
tree problem in Proposition~\ref{prop-coupling}.
Proposition~\ref{lemma:CLT-sum2} shall be proved in
Section~\ref{sec-lemma:CLT-sum2}. An important ingredient in the proof
will be the comparison
of the variables $\{B_m\}_{m=2}^{m_n}$, for an appropriately chosen $m_n$,
to an i.i.d. sequence. Results in this direction have been proved in
\cite{hofs3,hofs1}, and we shall combine these to the following statement:

\begin{Proposition}[(Coupling the forward degrees to an independent sequence)]
\label{prop-indep}
In the CM with $\tau>2$, there exists a $\rho>0$ such that
the random vector $\{B_m\}_{m=2}^{n^\rho}$ can be coupled to an
\textit{independent} sequence of random variables\vspace*{1pt}
$\{\Bindep_m\}_{m=2}^{n^\rho}$ with probability mass function $g$ in
(\ref{eqn:size-bias})
in such a way that $\{B_m\}_{m=2}^{n^\rho}=\{\Bindep_m\}
_{m=2}^{n^\rho}$
w.h.p.
\end{Proposition}

In Proposition~\ref{prop-indep}, in fact, we can take $\{B_m\}
_{m=2}^{n^\rho}$
to be the forward degree of the vertex to which any collection of
$n^\rho$
distinct stubs has been connected.

%s4.3 ###
\subsection{Flow clusters started from two vertices}
\label{sec-two_flows}
To compute the hopcount, we first grow the SWG from vertex 1 until
time $a_n$,
followed by the growth of the SWG from vertex 2 until the two SWGs meet,
as we now explain in more detail. Denote by $\{\SWG_m^{(i)}\}
_{m=0}^{\infty}$
the SWG from the vertex $i\in\{1,2\}$, and, for $m\geq0$, let
%
%e4.15 ###
%
\begin{equation} \label{SWG12-union} \SWG^{(1,2)}_m=\SWG^{
(1)}_{a_n} \cup\SWG^{(2)}_{m},
\end{equation}
the union of the SWGs of vertex 1 and 2. We shall only consider
values of $m$ where $\SWG^{(1)}_{a_n}$ and
$\SWG^{(2)}_{m}$ are \textit{disjoint}, that is,
they do not contain any common (real) vertices. We shall discuss
the moment when they connect in Section~\ref{sec-connedge} below.

We recall the notation in Section~\ref{sec-coupling},
and, for $i\in\{1,2\}$, denote by $\AS_{m}^{(i)}$ and
$\Art_{m}^{(i)}$ the number of allowed and artificial stubs in
$\SWG^{(i)}_{m}$. We let the set of free stubs $\FS_m$
consist of those stubs which have not yet been
paired in $\SWG^{(1,2)}_m$ in~(\ref{SWG12-union}).
Apart from that, the evolution of $\SWG^{(2)}_{m}$,
following the evolution of $\SWG^{(1)}_{a_n}$,
is identical as in Construction~\ref{const-SWGCM}.
We denote by $S_m^{(i)}=|\AS_{m}^{(i)}|$
the number of allowed stubs in $\SWG^{(i)}_{m}$ for $i\in\{1,2\}$.
We define $B_m^{(i)}$ accordingly.

The above description shows how
we can grow the SWG from vertex 1 followed by the one of vertex 2.
In order to state an adaptation of Proposition~\ref{prop-coupling} to
the setting
where the SWGs of vertex 1 is first grown to size $a_n$, followed by
the growth
of the SWG from vertex 2 until the connecting edge appears,
we let the random time $R_m^{(i)}$ be the first time $l$ such that
$\SWG_{l}^{(i)}$
consists of $m+1$ real vertices. Then our main coupling result for two
simultaneous
SWGs is as follows:

\begin{Proposition}[(Coupling SWGs on two trees and CM from two vertices)]
\label{prop-coupling12}
Jointly for $m\geq0$, as long as the sets of real vertices in
$(\SWG^{(1)}_{a_n},\SWG^{(2)}_{m})$
are \textit{disjoint}, these sets are equal in distribution to the sets
of $j_1$th, respectively $j_2$th, closest vertices
to vertex 1 and 2, respectively, for $j_1=1, \ldots, R_{a_n}^{(1)}$
and $j_2=1, \ldots, R_{m}^{(2)}$, respectively.
\end{Proposition}

%Call the edges linking the
%two SWGs the \textit{potential connecting edges},
%and recall that the process is completed when a first connecting edge
%is
%chosen. Now think of the edges as consisting of two stubs, and
%give each stub connected to the union of the two SWG but whose brother
%stub is not
%yet part of the union of the two SWGs an independent exponential 1
%weight. The weights of all stubs is independent. Then we grow
%the two SWGs simultaneously by first choosing the stub
%with minimal weight incident to the union of the
%two SWGs, and next connecting this stub to a uniformly chosen
%stub that is still available. The precise claim is that this
%description yields SWGs which have the same distribution as
%the original SWGs. The reason is that, when we study the
%flow from the vertices 1 and 2 simultaneously, the edges linking
%the two SWGs are filled at twice the rate. Thus, when we choose
%the next minimal stub, we choose the potential connecting edges
%at \textit{twice} the rate, which is the same as if these edges would
%have weights which are exponential random variables with parameter 2.
%Since an exponential random variable with parameter 2 is the minimum of
%two exponential random variables with parameter 1, we can equally give
%each
%of the two \textit{stubs} of which the potential connecting edge is
%built
%up
%a weight which is an independent exponential random variable with
%parameter 1,
%thus proving the claim.
%We now introduce some notation.

%s4.4 ###
\subsection{The connecting edge}
\label{sec-connedge}
As described above, we grow the two SWGs until
the first stub with minimal weight incident to $\SWG_m^{(2)}$
is paired to a stub incident to $\SWG_{a_n}^{(1)}$.
We call the created edge linking the
two SWGs the \textit{connecting edge}. More precisely, let
%
%e4.16 ###
%
\begin{equation} \label{conn-edge-time} \CE_n=\min\bigl\{m\geq0
\dvtx\SWG
^{(1)}_{a_n}\cap\SWG^{(2)}_{m}\neq\varnothing\bigr\}
\end{equation}
be the first time that $\SWG^{(1)}_{a_n}$ and $\SWG^{(2)}_{m}$
share a vertex. When $m=0$, this means that $2\in\SWG^{(1)}_{a_n}$
(which we shall show happens with small probability), while when $m\geq1$,
this means that the $m$th-stub of $\SWG^{(2)}$
which is chosen and then paired, is paired to a stub from $\SWG^{
(1)}_{a_n}$.
The path found actually is the shortest-weight path between vertices 1
and 2,
since $\SWG^{(1)}_{a_n}$ and $\SWG^{(2)}_{m}$
precisely consists of the closest real vertices to the root $i$,
for $i=1,2$, respectively.

We now study the probabilistic properties of
the connecting edge. Let the edge $e=st$ be incident to $\SWG^{
(1)}_{a_n}$, and
$s$ and $t$ denote its two stubs. Let the vertex incident to $s$ be
$i_s$ and
the vertex incident to $t$ be $i_t$. Assume that $i_s\in\SWG^{
(1)}_{a_n}$,\vspace*{1pt}
so that, by construction, $i_t\notin\SWG^{(1)}_{a_n}$.
Then, conditionally on $\SWG^{(1)}_{a_n}$
and $\{T_i^{(1)}\}_{i=1}^{a_n}$, the weight of $e$ is at least
$T^{(1)}_{a_n}-W^{(1)}_{i_s}$, where $W^{(1)}_{i_s}$ is
the weight of the\vspace*{1.5pt} shortest path
from 1 to $i_s$. By the memoryless property of the exponential distribution,
therefore, the weight on edge $e$ equals $T^{(1)}_{a_n}-W^{
(1)}_{i_s}+E_e$,
where the collection $(E_e)$, for all $e$ incident to $\SWG^{(1)}_{a_n}$
are i.i.d. $\operatorname{Exp}(1)$ random variables. Alternatively,\vspace*{1pt} we
can redistribute
the weight
by saying that the stub $t$ has weight $E_e$, and the stub $s$ has weight
$T^{(1)}_{a_n}-W^{(1)}_{i_s}$. Further, in the growth of
$(\SWG^{(2)}_{m})_{m\geq0}$, we can also think of the
exponential weights
of the edges incident to $\SWG^{(2)}_{m}$ being positioned on the
\textit{stubs} incident to $\SWG^{(2)}_{m}$. Hence, there is no
distinction between the stubs that are part of edges connecting $\SWG
^{(1)}_{a_n}$ and
$\SWG^{(2)}_{m}$ and the stubs that are part of edges incident to
$\SWG^{(2)}_{m}$, but
not to $\SWG^{(1)}_{a_n}$. Therefore, in the growth of
$(\SWG^{(2)}_{m})_{m\geq0}$, we can think of the minimal
weight \textit{stub} incident to $\SWG^{(2)}_{m}$ being chosen
uniformly at random, and then a uniform free stub is chosen
to pair it with. As a result, the distribution of the stubs chosen
\textit{at the time of connection} is equal to any of the other (real)
stubs chosen along the way.
This is a crucial ingredient to prove the scaling of the
shortest-weight path
between vertices 1 and 2.

For $i\in\{1,2\}$, let $H_n^{(i)}$ denote the length of the
shortest-weight path
between vertex $i$ and the common vertex in $\SWG_{a_n}^{(1)}$
and $\SWG_{\CE_n}^{(2)}$,
so that
%
%e4.17 ###
%
\begin{equation} \label{Hn-sum} H_n=H_n^{(1)}+H_n^{(2)}.
\end{equation}
Because of the fact that at time $\CE_n$ we have found the
shortest-weight path, we have that
%
%e4.18 ###
%
\begin{equation} \label{tildeHn-def} \bigl(H_n^{(1)},H_n^{
(2)}\bigr)\stackrel{d}{=}\bigl(G^{(1)}_{a_n+1}-1,G^{(2)}_{\CE
_n}\bigr),
\end{equation}
where $\{G^{(1)}_m\}_{m=1}^{\infty}$ and $\{G^{(2)}_m\}
_{m=1}^{\infty}$ are copies
of the process in~(\ref{Gm-def}), which are \textit{conditioned on
drawing a real stub}.
Indeed, at the time of the connecting edge,
a uniform (real) stub of $\SWG^{(2)}_{m}$ is drawn, and it is
paired to a uniform (real) stub of
$\SWG^{(1)}_{a_n}$. The number of hops in $\SWG^{
(1)}_{a_n}$ to the end of the attached
edge is therefore equal in distribution to $G^{(1)}_{a_n+1}$
conditioned on
drawing a real stub. The $-1$ in~(\ref{tildeHn-def}) arises since the
connecting edge
is counted twice in $G^{(1)}_{a_n+1}+G^{(2)}_{\CE_n}$.
The\vspace*{-1pt} processes $\{G^{(1)}_m\}_{m=1}^{\infty}$ and
$\{G^{(2)}_m\}_{m=1}^{\infty}$ are conditionally independent
given the realizations of $\{B_m^{(i)}\}_{m=2}^n$.

Further, because of the way the weight of the potential connecting edges
has been distributed over the two stubs out of which the connecting
edge is
comprised, we have that
%
%e4.19 ###
%
\begin{equation} \label{tildeWn-def} \Wn=T^{(1)}_{a_n}+ T^{
(2)}_{\CE_n},
\end{equation}
where $\{T^{(1)}_m\}_{m=1}^{\infty}$ and $\{T^{(2)}_m\}
_{m=1}^{\infty}$ are two
copies of the process $\{T_m\}_{m=1}^{\infty}$ in~(\ref{Tm-def}), again
conditioned on drawing a real stub. Indeed, to see~(\ref
{tildeWn-def}), we note that the weight of the connecting edge
is equal to the sum of weights of its two stubs. Therefore, the weight of
the shortest weight path is equal to the sum of the weight within $\SWG
^{(1)}_{a_n}$,
which is equal to $T_{a_n}^{(1)}$, and the weight within $\SWG
^{(2)}_{\CE_n}$,
which is equal to $T^{(2)}_{\CE_n}$.

In the distributions in~(\ref{tildeHn-def}) and~(\ref
{tildeWn-def}) above, we
always condition on drawing a real stub. Since we shall show that this occurs
w.h.p., this conditioning plays a minor role.

%The reason for this difference in weights
%is that an exponential random variable is attached to each
%rather than to each \textit{edge}. For the connecting edge, only the
%weight
%of the \textit{chosen} stub should be counted, not the weight
%of the \textit{paired} stub. This is achieved by taking the weight to
%the
%vertex \textit{incident}
%to the paired stub. We shall show that,
%for the CM, the dependence
%of $\{(G^{\sss(1)}_m, T^{\sss(1)}_m)\}_{m=1}^{\infty}$ and
%$\{(G^{\sss(2)}_m,T^{\sss(2)}_m)\}_{m=1}^{\infty}$ is rather weak, and
%similar
%results apply for $\{\widetilde T^{\sss(i)}_m\}_{m=1}^{\infty}$.
%In the sequel, we shall denote
% \begin{eqnarray}
% \lbeq{tildeWn-def}
% (\widetilde W_n^{\sss(1)}, \widetilde W_n^{\sss(2)})
% &=&\begin{cases}
% (T^{\sss(1)}_{\sss\lceil\CE_n/2\rceil},\widetilde T^{\sss(2)}_{\sss
% \CE_n \mbox{ is odd,}\\
% (\widetilde T^{\sss(1)}_{\sss\lceil\CE_n/2\rceil},T^{\sss(2)}_{\sss
% \CE_n \mbox{ is even,}
% \end{cases}\\
% \lbeq{tildeHn-def}
% (\widetilde H_n^{\sss(1)}, \widetilde H_n^{\sss(2)})
% &=&\begin{cases}
% (G^{\sss(1)}_{\sss\lceil\CE_n/2\rceil},\widetilde G^{\sss(2)}_{\sss
% \CE_n \mbox{ is odd,}\\
% (\widetilde G^{\sss(1)}_{\sss\lceil\CE_n/2\rceil},G^{\sss(2)}_{\sss
% \CE_n \mbox{ is even.}
% \end{cases}
% \end{eqnarray}
%Then, by \refeq{Wn-descr1}--\refeq{Wn-descr2}, we have that
% \eqn{
% \lbeq{Wn-descr-fin}
% \Wn=\widetilde W_n^{\sss(1)}+\widetilde W_n^{\sss(2)},
%
% \Hn=\widetilde H_n^{\sss(1)}+\widetilde H_n^{\sss(2)}.
% }

We shall now intuitively explain why the leading order asymptotics
of $\CE_n$
is given by $a_n$ where $a_n$ is
defined in~(\ref{an-def}). For this, we must know
how many allowed stubs there are, that is, we
must determine how many stubs there are incident to
the union of the two SWGs at any time.
Recall that $S_m^{(i)}$ denotes the number of allowed
stubs in the SWG from vertex $i$ at time $m$. The total number of allowed
stubs incident to $\SWG^{(1)}_{a_n}$ is $S_{a_n}^{(1)}$,
while the number incident to $\SWG^{(2)}_m$ is equal to
$S_m^{(2)}$, and where
%
%e4.20 ###
%
\begin{equation} S_m^{(i)}=D_i+\sum_{l=2}^{m} \bigl(B^{
(i)}_{l}-1\bigr).
\end{equation}
We also write $\Art_m=\Art_{a_n}^{(1)}\cup\Art_{m}^{(2)}$.

Conditionally on $\SWG_{a_n}^{(1)}$ and
$\{(S_{l}^{(2)}, \Art_{l}^{(2)})\}_{l=1}^{m-1}$
and $L_n$, and assuming that $|\Art_m|$, $m$ and $S_m$ satisfy
appropriate bounds, we obtain
%
%e4.21 ###
%
\begin{equation} \prob(\CE_n=m|\CE_n>m-1)\approx\frac
{S_{a_n}^{(1)}}{L_n}.
\end{equation}
%
%while, conditionally on
%$\{(S_{\lceil l/2\rceil}^{\sss(1)},\Art_{\lceil l/2\rceil}^{\sss(1)}
%S_{\lfloor l/2\rfloor}^{\sss(2)}, \Art_{\lfloor l/2\rfloor}^{\sss(2)})
%and $L_n$,
% \eqn{
% \prob(\CE_n=2m+1|\CE_n>2m)\approx\frac{S_{m}^{\sss(2)}}{L_n}.
% }
%We can summarize these formulas by the fact that
%conditionally on $\{(S_{\lceil l/2\rceil}^{\sss(1)},
%S_{\lfloor l/2\rfloor}^{\sss(2)})\}_{l=1}^m$ and $L_n$,
% \eqn{
% \lbeq{CEn-comp}
% \prob(\CE_n=m|\CE_n>m-1)\approx\frac{S_{\lfloor m/2\rfloor}^{
% }
%where $i_m=1+(m\mod2)$.

When $\tau\in(2,3)$ and~(\ref{Fcond(2,3)}) holds,
then $S_{l}^{(i)}/l^{1/(\tau-2)}$ can be expected to converge
in distribution to a stable random variable with parameter
$\tau-2$, while, for $\tau>3$, $S_{l}^{(i)}/l$ converges in
probability to $\nu-1$, where $\nu$ is defined in~(\ref{nu-def}).
We can combine these two statements by saying that
$S_{l}^{(i)}/l^{1/(\tau\wedge3-2)}$ converges in
distribution.
Note that the typical size $a_n$ of $\CE_n$ is such that, uniformly in $n$,
$\prob(\CE_n\in[a_n,2a_n])$ remains in $(\varepsilon, 1-\varepsilon)$,
for some $\vep\in(0,\frac12)$, which is the case when
%
%e4.22 ###
%
\begin{eqnarray} \prob(\CE_n\in[a_n,2a_n])&=&\sum_{m=a_n}^{2a_n}
\prob(\CE_n=m|\CE_n>m-1)\prob(\CE_n>m-1)\nonumber\\[-8pt]\\
[-8pt]\nonumber
&\in&(\vep, 1-\vep)
\end{eqnarray}
uniformly as $n\rightarrow\infty$. By the above discussion, and for
$a_n\leq m\leq2a_n$,
we have $\prob(\CE_n=m|\CE_n>m-1)=\Theta(m^{1/(\tau\wedge
3-2)}/n)=\Theta(a_n^{1/(\tau\wedge3-2)}/n),$
and $\prob(\CE_n>m-1)=\Theta(1)$.
Then we arrive at
%
%e4.23 ###
%
\begin{equation} \label{CE-comp} \prob(\CE_n\in[a_n,2a_n])=\Theta
\bigl(a_n a_n^{1/(\tau\wedge3-2)}/n\bigr),
\end{equation}
which remains uniformly positive and bounded for $a_n$ defined in
(\ref{an-def}).
In turn, this suggests that
%
%e4.24 ###
%
\begin{equation} \label{Mn-lim-heur} \CE_n/a_n \convd M
\end{equation}
for some limiting random variable $M$.

We now discuss what happens when~(\ref{Fcond}) holds for some
$\tau\in(2,3)$, but~(\ref{Fcond(2,3)}) fails. In this case,
there exists a slowly varying function $n\mapsto\ell(n)$ such that
$S_{l}^{(i)}/(\ell(l)l^{1/(\tau-2)})$ converges in distribution.
Then following the above argument shows that the right-hand side
(r.h.s.) of~(\ref{CE-comp}) is
replaced by $\Theta(a_n a_n^{1/(\tau-2)}\ell(a_n)/ n)$ which remains
uniformly positive and bounded for $a_n$ satisfying $a_n^{(\tau
-1)/(\tau-2)}\times\ell(a_n)=n$.
By Bingham, Goldie and Teugels \cite{BinGolTeu89}, Theorem 1.5.12,
there exists a solution $a_n$ to the above equation
which satisfies that it is regularly varying with exponent $(\tau
-2)/(\tau-1)$, so that
%
%e4.25 ###
%
\begin{equation} \label{an-def-RV} a_n=n^{(\tau-2)/(\tau-1)}\ell
^*(n)
\end{equation}
for some slowly varying function $n\mapsto\ell^*(n)$ which depends only
on the distribution function $F$.

%By the above description, we obtain (recall also \refeq{Hn-sum})
% \eqn{
% \lbeq{Hn-sum-Tn}
% \Hn=G^{\sss(1)}_{\sss\lceil\CE_n/2\rceil}+G^{\sss(2)}_{\sss\lfloor
% }
%where $\{G^{\sss(1)}_{m}\}_{m=1}^{\infty}$ and
%$\{G^{\sss(2)}_{m}\}_{m=1}^{\infty}$ are two
%independent sequences with the same distribution as $G_m$
%in \refeq{Gm-def}--\refeq{bi-def}. The -1 in \refeq{Hn-sum-Tn}
%stems from the fact that the connecting edge is
%counted both in $G^{\sss(1)}_{\sss\lceil\CE_n/2\rceil}$ as well as in
%$G^{\sss(2)}_{\sss\lfloor\CE_n/2\rfloor}$, whereas it should only be
%counted once.
%Equation \refeq{Hn-sum-Tn} gives a description of $G_1, G_2$ in
%G_2=G^{\sss(2)}_{\sss\lfloor\CE_n/2\rfloor}$.

In the following proposition, we shall state the necessary
result on $\CE_n$ that we shall need in
the remainder of the proof. In its statement, we shall use the
symbol $o_{\prob}(b_n)$ to denote a random variable $X_n$ which
satisfies that $X_n/b_n\convp0$.

\begin{Proposition}[(The time to connection)]
\label{lemma:conn_edge}
As $n\rightarrow\infty$, under the conditions of Theorems~\ref{main>3}
and~\ref
{main(2,3)} respectively, and with $a_n$ as in~(\ref{an-def}),
%
%e4.26 ###
%
\begin{equation} \label{time-conn-edge} \log{\CE_n}-\log{a_n} =
o_{\prob}\bigl(\sqrt{\log{n}}\bigr).
\end{equation}
Furthermore, for $i\in\{1,2\}$, and with $\beta\geq1$
as in Proposition~\ref{lemma:CLT-sum1},
%
%e4.27 ###
%
\begin{equation} \label{clt-conn-edge} \biggl(\frac{G^{
(1)}_{a_n+1}-\beta\log{a_n}}{\sqrt{\beta\log{a_n}}}, \frac
{G^{(2)}_{\CE_n}-\beta\log{a_n}}{\sqrt{\beta\log{a_n}}}
\biggr) \convd(Z_1,Z_2),
\end{equation}
where $Z_1, Z_2$ are two independent standard normal random variables.
Moreover, with $\gamma$ as in Proposition~\ref{lemma:CLT-sum1},
there exist random variables $X_1, X_2$ such that
%
%e4.28 ###
%
\begin{equation} \label{weight-conn-edge} \bigl(T^{
(1)}_{a_n}-\gamma\log{a_n},T^{(2)}_{\CE_n}-\gamma\log
{a_n} \bigr) \convd(X_1,X_2).
\end{equation}
\end{Proposition}

We note that the main result in~(\ref{clt-conn-edge}) is not a simple
consequence of~(\ref{time-conn-edge}) and Proposition~\ref{lemma:CLT-sum2}.
The reason is that $\CE_n$ is a \textit{random variable},
which a priori depends on $(G^{(1)}_{a_n+1},G^{(2)}_m)$ for
$m\geq0$.
Indeed, the connecting edge is formed out of two stubs which
are not artificial, and thus the choice of stubs is not
completely uniform. However, since there are only few
artificial stubs, we can extend the proof of
Proposition~\ref{lemma:CLT-sum2}
to this case. Proposition~\ref{lemma:conn_edge} shall be proved in
Section~\ref{sec-lemma:conn_edge}.

%s4.5 ###
\subsection{The completion of the proof}
\label{sec-complpf}
By the analysis in Section~\ref{sec-connedge},
we know the distribution of the sizes of the
SWGs at the time when the connecting edge appears.
By Proposition~\ref{lemma:conn_edge}, we know the number
of edges and their weights used in
the paths leading to the two vertices of the
connecting edge together with its fluctuations.
In the final step, we need to combine these results
by averaging both over the \textit{randomness} of the
time when the connecting edge appears (which is a random
variable), as well as over the number of edges in the
shortest weight path when we know the time the
connecting edge appears. Note that by~(\ref{tildeHn-def})
and Proposition~\ref{lemma:conn_edge},
we have, with $Z_1, Z_2$ denoting independent standard normal random variables,
and with $Z=(Z_1+Z_2)/\sqrt{2}$, which is again standard normal,
%
%e4.29 ###
%
\begin{eqnarray} \Hn&\stackrel{d}{=}&G_{a_n+1}^{(1)}+G_{\CE
_n}^{(2)}-1\nonumber\\
&=&2\beta\log{a_n} +Z_1\sqrt{\beta\log
{a_n}}+Z_2\sqrt{\beta\log{a_n}}+o_{\prob}\bigl(\sqrt{\log
{n}}\bigr)\\
&=&2\beta\log{a_n} +Z\sqrt{2\beta\log
{a_n}}+o_{\prob}\bigl(\sqrt{\log{n}}\bigr).\nonumber
\end{eqnarray}
Finally, by~(\ref{an-def}), this gives~(\ref{CLT-hopcount>3}) and
(\ref{CLT-hopcount(2,3)}) with
%
%e4.30 ###
%
\begin{equation} \label{alpha-descr} \alpha=\lim_{n\rightarrow
\infty} \frac{2\beta\log{a_n}}{\log{n}},
\end{equation}
which equals $\alpha=\nu/(\nu-1),$ when $\tau>3,$ since $\beta=\nu
/(\nu-1)$
and\vspace*{-2pt} $\frac{\log{a_n}}{\log{n}}=1/2$, and
$\alpha=2(\tau-2)/(\tau-1),$ when $\tau\in(2,3),$ since $\beta=1$ and
$\frac{\log{a_n}}{\log{n}}=(\tau-2)/(\tau-1).$ This completes the proof
for the hopcount.

In the description of $\alpha$ in~(\ref{alpha-descr}), we note that
when $a_n$ contains a slowly varying function for $\tau\in(2,3)$ as
in~(\ref{an-def-RV}), then
the result in Theorem~\ref{main(2,3)} remains valid with $\alpha\log{n}$
replaced by
%
%e4.31 ###
%
\begin{equation} 2\log{a_n}=\frac{2(\tau-2)}{\tau-1}\log{n}+2\log
{\ell^*(n)}.
\end{equation}

For the weight of the minimal path, we make use
of~(\ref{tildeWn-def}) and~(\ref{weight-conn-edge}) to obtain in a
similar way that
%
%e4.32 ###
%
\begin{equation} W_n-2\gamma\log{a_n}\convd X_1+X_2.
\end{equation}
This completes the proof for the weight of the shortest path.

%s5 ###
\section[Proof of Proposition 4.3]{Proof of Proposition \protect\ref{lemma:CLT-sum1}}
\label{sec-lemma:CLT-sum1}

%s5.1 ###
\subsection[Proof of Proposition 4.3(a)]{Proof of Proposition \protect\ref{lemma:CLT-sum1}(\textup{a})}
\label{sec-pfprop(a)}
%%%%%%%%%%%%%%%%%%%%%%%%%%%%%%%%%%%%%%%%%%%%%%%%%%%%%%%%%%
We start by proving the statement for $\tau\in(2,3)$. Observe that,
in this context,
$d_i=B_i$, and, by~(\ref{si-def}), $B_1+\cdots+B_i=S_i+ i-1$, so that
the sequence $B_j/(S_i+i-1),$ for $j$ satisfying $1\le j\le i,$
is exchangeable for each $i\ge1$. Therefore, we define
%
%e5.1 ###
%
\begin{equation} \label{hatGm-def} \hatG_m=\sum_{i=1}^m \hatI
_i, \qquad\mbox{where\ } \prob(\hatI_i=1|\{B_i\}_{i=1}^{\infty
})=\frac{B_i}{S_i+i-1}.
\end{equation}
Thus, $\hatI_i$ is, conditionally on $\{B_i\}_{i=1}^{\infty}$,
stochastically dominated by $I_i$, for each $i$,
which, since the sequences $\{\hatI_i\}_{i=1}^{\infty}$ and
$\{I_i\}_{i=1}^{\infty}$, conditionally on $\{B_i\}_{i=1}^{\infty}$,
each have independent components, implies that
$\hatG_m$ is stochastically dominated by $G_m$. We take $\hatG_m$ and $G_m$
in such a way that $\hatG_m\leq G_m$ a.s. Then, by the Markov
inequality, for $\kappa_m>0$,
%
%e5.2 ###
%
\begin{eqnarray}
\label{bd-GmhatGm}
\prob(|G_m-\hatG_m|\geq\kappa_m)
&\leq& \kappa_m^{-1} \expec[|G_m-\hatG_m|]=\kappa_m^{-1} \expec
[G_m-\hatG_m]\nonumber\\
&=& \kappa_m^{-1} \sum_{i=1}^{m} \expec\biggl[\frac
{B_i(i-1)}{S_i(S_i+i-1)}\biggr]\\
&=&
\kappa_m^{-1} \sum_{i=1}^{m} \frac{i-1}{i}\expec[1/S_i],\nonumber
\end{eqnarray}
where, in the second equality, we used the exchangeability of
$B_j/(S_i+i-1), 1\le j \le i$.
We will now show that
%
%e5.3 ###
%
\begin{equation}
\label{finitesum}
\sum_{i=1}^{\infty} \expec[1/S_i]<\infty,
\end{equation}
so that for any $\kappa_m\rightarrow\infty$,
we have that $\prob(|G_m-\hatG_m|\leq\kappa_m)\to1$. We can then
conclude that
the CLT for $G_m$ follows from the one for $\hatG_m$.
By Deijfen et al. \cite{DeiEskHofHoo09}, (3.12) for $s=1$, for $\tau
\in(2,3)$ and using that
$S_i=B_1+\cdots+ B_i-(i-1)$, where $\prob(B_1>k)=k^{2-\tau}L(k)$,
there exists a
slowly varying function $i\mapsto l(i)$ such that $\expec[1/S_i]\leq
cl(i)i^{-1/(\tau-2)}$. When $\tau\in(2,3)$, we have that $1/(\tau
-2)>1$, so that~\eqref{finitesum} follows.

We now turn to the CLT for $\hatG_m$.
Observe from the exchangeability of $B_j/(S_i+i-1),$ for $1\le j \le i,$
%we note that
%the sequence $\{\hatI_i\}_{i=1}^{\infty}$ is \textit{independent},
%since,
that for $i_1<i_2<\cdots<i_k$,
%
%e5.4 ###
%
\begin{eqnarray}
\label{ind-hatI}
\prob(\hatI_{i_1}=\cdots=\hatI_{i_k}=1)
&=&\expec\Biggl[\prod_{l=1}^k \frac{B_{i_l}}{S_{i_l}+i_l-1}
\Biggr]\nonumber\\
&=&\expec\Biggl[\frac{B_{i_1}}{S_{i_1}+i_1-1}\prod_{l=2}^k \frac
{B_{i_l}}{S_{i_l}+i_l-1} \Biggr]\\
&=&\frac{1}{i_1}\expec\Biggl[\prod_{l=2}^k \frac
{B_{i_l}}{S_{i_l}+i_l-1} \Biggr]
=\cdots=\prod_{l=1}^k \frac{1}{i_l},\nonumber
\end{eqnarray}
where we used that since $B_1+\cdots+B_j=S_j+j-1$,
\begin{eqnarray*}
\expec\Biggl[\frac{B_{i_1}}{S_{i_1}+i_1-1}\prod_{l=2}^k \frac
{B_{i_l}}{S_{i_l}+i_l-1} \Biggr]
&=&
\frac1{i_1}\sum_{i=1}^{i_1}
\expec\Biggl[\frac{B_{i}}{S_{i_1}+i_1-1}\prod_{l=2}^k \frac
{B_{i_l}}{S_{i_l}+i_l-1} \Biggr]\\
&=&
\frac{1}{i_1}\expec\Biggl[\prod_{l=2}^k \frac
{B_{i_l}}{S_{i_l}+i_l-1} \Biggr].
\end{eqnarray*}
Since $\hatI_{i_1},\ldots,\hatI_{i_k}$ are indicators this implies
that $\hatI_{i_1},\ldots,\hatI_{i_k}$ are independent.
Thus $\hatG_m$ has the same distribution as $\sum_{i=1}^m J_i$
where $\{J_i\}_{i=1}^{\infty}$ are \textit{independent}
Bernoulli random variables with $\prob(J_i=1)=1/i$. It is
a standard consequence of the Lindeberg--L\'evy--Feller CLT
that $(\sum_{i=1}^m J_i-\log{m})/\sqrt{\log{m}}$ is asymptotically
standard normally distributed.
%%%%%%%%%%%%%%%%%%%%%%%%%%%%%%%%%%%%%%%%%%%%%%%%%%%%%%%%%%%%%%%%%%%%
%
\begin{Remark}[(Extension to exchangeable setting)]
\label{rem-CLTexchang}
Note that the CLT for $G_m$ remains valid when
(i) the random variables $\{B_i\}_{i=1}^m$ are \textit{exchangeable},
with the same marginal distribution as in the i.i.d. case, and
(ii) $\sum_{i=1}^{m} \expec[1/S_i]=o(\sqrt{\log{m}})$.
\end{Remark}

The approach for $\tau>3$ is different from that of
$\tau\in(2,3)$. For $\tau\in(2,3)$, we coupled $G_m$ to $\hatG_m$
and proved that $\hatG_m$ satisfies the CLT with the correct
norming constants. For $\tau>3,$ the case we consider now, we
first apply a \textit{conditional} CLT, using the
Lindeberg--L\'evy--Feller condition, stating that, conditionally on
$B_1,B_2,\ldots$ satisfying
%
%e5.5 ###
%
\begin{equation} \label{LLF-cond-CLT} \lim_{m\to\infty} \sum
_{j=1}^m \frac{B_j}{S_j} \biggl(1-\frac{B_j}{S_j} \biggr)=\infty,
\end{equation}
we have that
%
%e5.6 ###
%
\begin{equation} \label{clt-LLF} \frac{G_m-\sum_{j=1}^m B_j/S_j}
{ ( \sum_{j=1}^m {B_j}/{S_j}(1-{B_j}/{S_j})
)^{1/2}}\convd Z,
\end{equation}
where $Z$ is standard normal. The result~(\ref{clt-LLF})
is also contained in \cite{buhler}.

Since $\nu=\expec[B_j]>1$ and $\expec[B_j^a]<\infty,$
for any $a<\tau-2,$ it is not hard to see that the random variable
$\sum_{j=1}^\infty B_j^2/S_j^2$ is positive and has finite first moment,
so that for $m\to\infty$,
%
%e5.7 ###
%
\begin{equation} \label{sec-mom-sum} \sum_{j=1}^m B_j^2/S_j^2=\Op
(1),
\end{equation}
where $\Op(b_m)$ denotes a sequence of random variables $X_m$ for
which $|X_m|/b_m$
is tight.

We claim that
%
%e5.8 ###
%
\begin{equation} \label{suff-cond} \sum_{j=1}^m B_j/S_j-\frac{\nu
}{\nu-1}\log{m}=o_{\prob}\bigl(\sqrt{\log{m}}\bigr).
\end{equation}
Obviously,~(\ref{clt-LLF}),~(\ref{sec-mom-sum}) and~(\ref{suff-cond})
imply Proposition~\ref{lemma:CLT-sum1}(a) when $\tau>3$.
%and first show that \refeq{suff-cond} is a sufficient
%condition for Proposition~\ref{lemma:CLT-sum1}(a) when $\tau>3$.
%For this, we note that \refeq{clt-LLF}, \refeq{sec-mom-sum} and
%imply that
% \eqn{
% \lbeq{clt2}
% \frac{G_m-\sum_{j=1}^m \frac{B_j}{S_j}}
% {\sqrt{\frac{\nu}{\nu-1}\log{m}}}\convd Z,
% }
%and we then split
% \eqn{
% \lbeq{clt3}
% \frac{G_m-\frac{\nu}{\nu-1}\log{m}}
% {\sqrt{\frac{\nu}{\nu-1}\log{m}}}
% =\frac{G_m-\sum_{j=1}^m \frac{B_j}{S_j}}{\sqrt{\frac{\nu}{\nu-1}
% +\frac{\sum_{j=1}^m \frac{B_j}{S_j}-\frac{\nu}{\nu-1}\log{m}}
% {\sqrt{\frac{\nu}{\nu-1}\log{m}}}.
% }
%By \refeq{clt2}, the first term converges in distribution to a standard
%normal random variable, while by \refeq{suff-cond}, the second term
%converges in probability to 0. This completes the proof subject to

In order to prove~(\ref{suff-cond}), we split
\begin{eqnarray} \label{split} &&\sum_{j=1}^m B_j/S_j-\frac{\nu
}{\nu
-1}\log{m}\nonumber\\[-8pt]\\[-8pt]
&&\qquad= \Biggl(\sum_{j=1}^m (B_j-1)/S_j-\log{m} \Biggr)+ \Biggl
(\sum
_{j=1}^m 1/S_j-\frac{1}{\nu-1}\log{m} \Biggr),\nonumber
\end{eqnarray}
and shall prove that each of these two terms on the r.h.s. of~(\ref
{split}) is
$o_{\prob}(\sqrt{\log{m}})$. For the first term, we note from
the strong law of large
numbers that
% \eqn{
% \sum_{j=1}^m \log{ (\frac{S_{m}}{S_{m-1}} )}
% =\log{S_m}-\log{S_0}=\log{S_m},
% }
%since $S_0=1$. Now, by the law of large numbers,
% \eqn{
% \log{S_m}-\log{m}
% =\log{(S_m/m)}\convas\log{(\nu-1)}.
% }
%Thus,
%
%e5.9 ###
%
\begin{equation} \sum_{j=1}^m \log\biggl(\frac{S_{j}}{S_{j-1}}
\biggr)=\log{S_m}-\log{S_0}=\log{m}+\Op(1).
\end{equation}
Also, since $-\log{(1-x)}=x+O(x^2)$, we have that
\begin{eqnarray} \sum_{j=1}^m \log{ (S_{j}/S_{j-1} )} &=&-\sum
_{j=1}^m \log\bigl(1-(B_j-1)/S_j \bigr)\nonumber\\[-8pt]\\[-8pt]
&=&\sum_{j=1}^m (B_j-1)/S_j
+O \Biggl(\sum_{j=1}^m (B_j-1)^2/S_j^2 \Biggr).\nonumber
\end{eqnarray}
Again, as in~(\ref{sec-mom-sum}), for $m\rightarrow\infty$,
%
%e5.10 ###
%
\begin{equation} \sum_{j=1}^m (B_j-1)^2/S_j^2=\Op(1),
\end{equation}
so that
%
%e5.11 ###
%
\begin{equation} \label{lln(a)} \sum_{j=1}^m (B_j-1)/S_j-\log{m}=\Op
(1).
\end{equation}
In order to study the second term on the right-hand side of~(\ref{split}),
we shall prove a slightly
stronger result than necessary, since we shall also use this later on.
Indeed, we shall show that there exists a random variable $Y$ such that
%
%e5.12 ###
%
\begin{equation} \label{as-conv-sum} \sum_{j=1}^m 1/S_j-\frac{1}{\nu
-1}\log{m}\convas Y.
\end{equation}
In fact, the proof of~(\ref{as-conv-sum}) is a consequence of \cite
{AthKar67}, Theorem~1,
since $\expec[(B_i-1)\log(B_i-1)]<\infty$ for $\tau>3$. We decided
to give a separate
proof of~(\ref{as-conv-sum}) which can be easily adapted to the
exchangeable case.

To prove~(\ref{as-conv-sum}), we write
%
%e5.13 ###
%
\begin{equation} \sum_{j=1}^m 1/S_j-\frac{1}{\nu-1}\log{m} =\sum
_{j=1}^m \frac{(\nu-1) j-S_j}{S_j(\nu-1)j}+\Op(1),
\end{equation}
so that in order to prove~(\ref{as-conv-sum}),
it suffices to prove that, uniformly in $m\geq1$,
%
%e5.14 ###
%
\begin{equation} \label{suff-to-bound} \sum_{j=1}^m \frac{|S_j-(\nu
-1) j|}{S_j(\nu-1)j}<\infty\qquad \mathrm{a.s.}
\end{equation}
Thus, if we further make use of the fact that $S_j\geq\eta j$
except for at most finitely many $j$ (see also Lemma~\ref
{lem-cond-LD-est} below),
then we obtain that
%
%e5.15 ###
%
\begin{equation} \qquad\Bigg|\sum_{j=1}^m \frac{1}{S_j}-\frac
{1}{\nu
-1}\log{m}\Bigg| \leq\sum_{j=1}^m \frac{|S_j-(\nu-1) j|}{S_j(\nu
-1)j}+\Op(1) \leq C\sum_{j=1}^m \frac{|S^*_j|}{j^2},
\end{equation}
where $S^*_j=S_j-\expec[S_j]$, since $\expec[S_j]=(\nu-1)j+1$.
We now take the expectation, and conclude that for any $a>1$,
Jensen's inequality for the convex function $x\mapsto x^a$, yields
%use the H\"older inequality,
%to obtain, that for any $a>1$,
%
%e5.16 ###
%
\begin{equation} \expec[|S_j^*|] \leq\expec[|S_j^*|^a]^{1/a}.
\end{equation}

To bound the last expectation, we will use a consequence of the
Marcinkiewicz--Zygmund inequality (see, e.g.,
\cite{Gut05}, Corollary~8.2, page~152).
Taking $1<a<\tau-2$, we have that
$\expec[|B_1|^a]<\infty$, since $\tau>3$, so that
%
%e5.17 ###
%
\begin{equation} \label{Marcinkiewicz-Zygmund} \expec\Biggl[\sum
_{j=1}^m \frac{|S_j^*|}{j^{2}} \Biggr] \leq\sum_{j=1}^m \frac
{\expec
[|S_j^*|^a]^{1/a}}{j^{2}}\leq\sum_{j=1}^m \frac{c_a^{1/a} \expec
[|B_1|^a]^{1/a}}{j^{2-1/a}} <\infty.
\end{equation}
This completes the proof of~(\ref{as-conv-sum}).

\begin{Remark}[(Discussion of exchangeable setting)]
\label{rem-CLT-tau>3}
When the random variables $\{B_i\}_{i=1}^m$ are \textit{exchangeable},
with the same marginal distribution as in the i.i.d. case, and with
$\tau>3$,
we note that to prove a CLT for $G_m$, it suffices to prove~(\ref{sec-mom-sum})
and~(\ref{suff-cond}). The proof of~(\ref{suff-cond}) contains two steps,
namely,~(\ref{lln(a)}) and~(\ref{suff-to-bound}). For the CLT to
hold, we
in fact only need that the involved quantities are $o_{\prob
}(\sqrt{\log{m}})$,
rather than $\Op(1)$. For this, we note that:

\begin{longlist}[(b)]
\item[(a)] the argument to prove~(\ref{lln(a)}) is rather flexible, and
shows that if (i)~$\log S_m/ m =o_{\prob}(\sqrt{\log{m}})$ and
if (ii) the condition in~(\ref{sec-mom-sum}) is satisfied with
$\Op(1)$ replaced by $o_{\prob}(\sqrt{\log{m}})$, then~(\ref
{lln(a)}) follows
with $\Op(1)$ replaced by $o_{\prob}(\sqrt{\log{m}})$;
\item[(b)] for the proof of~(\ref{suff-to-bound}) we will make use of
stochastic domination
and show that each of the stochastic bounds will satisfy~(\ref
{suff-to-bound}) with
$\Op(1)$ replaced by $o_{\prob}(\sqrt{\log{m}})$ (compare
Lemma~\ref{Mar-Zyg}).
%the proof of \refeq{suff-to-bound} makes
%use of the i.i.d. nature of $B_1, \ldots, B_m$
%in a more fundamental way, due to the use of the
%Marcinkiewicz-Zygmund inequality in \refeq{Marcinkiewicz-Zygmund},
%which requires the i.i.d. assumption.
%However, the fact that $\sum_{j=1}^m \frac{1}{S_j}-\frac{1}{\nu-1}
%follows when
% \eqn{
% \lbeq{cond-two-CLT}
% \sum_{j=1}^m \frac{S_j-(\nu-1) j}{S_j(\nu-1)j}=\op(\sqrt{\log{m}}).
% }
\end{longlist}
\end{Remark}

%s5.2 ###
\subsection[Proof of Proposition 4.3(b)]{Proof of Proposition \protect\ref{lemma:CLT-sum1}(\textup{b})}
\label{sec-pfprop(b)}
We again start by proving the result for $\tau\in(2,3)$.
It follows from~(\ref{Tm-def}) and the independence of $\{E_i\}_{i\geq
1}$ and
$\{S_i\}_{i\geq1}$ that, for the proof of~(\ref{Tm-conv}), it is
sufficient to show
that
%
%e5.18 ###
%
\begin{equation} \label{bound-mean-ratioS-indep} \sum_{i=1}^\infty
\expec[1/S_i]<\infty,
\end{equation}
which holds due to \eqref{finitesum}. %The argument for $\widetilde
%T_m$ is similar, with the same
%limit. Indeed, on the one hand, since $\widetilde T_m\leq T_m$, the
%limit of $\widetilde T_m$ cannot be larger
%than that of $T_m$. On the other hand, since $G_m\rightarrow\infty$
%and $m\mapsto T_m$ is
%increasing,\vspace*{1pt} we have that the $\widetilde T_m\geq T_k,$ w.h.p. for
%any $k$. Therefore,
%$\widetilde T_m$ must have the same limit as $T_m$.\vspace*{1pt}
%
%The extension of this result to $\tau>3$, where the weak limits of
%$T_m$ and $\widetilde T_m$ are \textit{different},
%is deferred to Appendix~\hyperref[sec-weak-conv-weight-tau3]{C}.

%s6 ###
\section[Proof of Proposition 4.6]{Proof of Proposition \protect\ref{lemma:CLT-sum2}}
\label{sec-lemma:CLT-sum2}

In this section, we extend the proof of Proposition~\ref{lemma:CLT-sum1}
to the setting where the random vector $\{B_i\}_{i=2}^m$ is \textit{not}
i.i.d., but rather corresponds to the vector of forward degrees in the CM.

\label{conditioning}
In the proofs for the CM, we shall always \textit{condition} on the fact
that the vertices under consideration are part of the giant component.
As discussed below~(\ref{nu-def}), in this case, the giant component
has size $n-o(n)$, so that each vertex is in the giant component
w.h.p. Further, this conditioning ensures that $S_j>0$
for every $j=o(n)$.

%The proofs in Section~\ref{sec-lemma:CLT-sum1} have been
%deliberately made as flexible as possible, in that we make use of
%the independence only when absolutely necessary.
We recall that the set up of the random variables involved in Proposition
\ref{lemma:CLT-sum2} is given in~(\ref{setup-6.4})
and~(\ref{def-Gm1}). The random variable $R_m$, defined in~(\ref
{Rm-def}), is the first time $t$ the SWG$_t$ consists of
$m+1$ real vertices.

\begin{Lemma}[(Exchangeability of $\{B_{R_m}\}_{m=1}^{n-1}$)]
\label{lem-exchange}
Conditionally on $\{D_i\}_{i=1}^{n}$, the sequence of
random variables $\{B_{R_m}\}_{m=1}^{n-1}$ is exchangeable, with
marginal probability distribution
%
%e6.1 ###
%
\begin{equation} \label{BT-law} \prob_n(B_{R_1}=j)=\sum_{i=2}^n
\frac{(j+1)\indic{D_i=j+1}}{L_n-D_1},
\end{equation}
where $\prob_n$ denotes the conditional probability given $\{D_i\}_{i=1}^n$.
\end{Lemma}

\begin{pf} We note that, by definition, the random variables
$\{B_{R_m}\}_{m=1}^{n-1}$ are equal to the forward degrees (where
we recall that the forward degree is equal to the degree minus 1) of
a vertex chosen from all vertices unequal to 1, where a vertex $i$
is chosen with probability proportional to its degree, that is, vertex
$i\in\{2, \ldots, n\}$ is chosen with probability
$P_i=D_i/(L_n-D_1)$. Let $K_2, \ldots, K_n$ be the vertices chosen;
then the sequence $K_2, \ldots, K_n$ has the same distribution as
draws with probabilities $\{P_i\}_{i=2}^n$ \textit{without
replacement}. Obviously, the sequence $(K_2, \ldots,
K_n)$ is exchangeable, so that the sequence $\{B_{R_m}\}_{m=1}^{n-1}$,
which can be identified as $B_{R_m}=D_{K_{m+1}}-1$, inherits this
property.
\end{pf}

We continue with the proof of Proposition~\ref{lemma:CLT-sum2}.
By Lemma~\ref{lem-exchange}, the sequence $\{B_j\}_{j=2}^m$
is exchangeable, when we condition on $|\Art_j|=0$ for all $j\leq m$.
Also, $|\Art_j|=0$ for all $j\leq m$ holds precisely when $R_m=m$.
In Lemma~\ref{lem-first-art} in Appendix~\hyperref[sec-app-A]{A}, the
probability that $R_{m_n}=m_n$,
for an appropriately chosen $m_n$, is investigated.
We shall make crucial use of this lemma to study $G_{m_n}$.

\begin{pf*}{Proof of Proposition~\ref{lemma:CLT-sum2}} Recall that by
definition
$\log(\olmn/a_n)=\break o(\sqrt{\log{n}})$. Then, we split, for some
$\ulmn$ such that
$\log(a_n/\ulmn)=o(\sqrt{\log{n}})$,
%
%e6.2 ###
%
\begin{equation} \label{Gm-split} G_{\olmn}=\widetilde G_{\ulmn}+
[G_{\olmn}-\widetilde G_{\ulmn}],
\end{equation}
where $\widetilde G_{\ulmn}$ has the same marginal distribution
as $G_{\ulmn}$, but also satisfies that $\widetilde G_{\ulmn}\leq
G_{\olmn},$
a.s. By construction, the sequence of random variables $m\mapsto G_m$ is
stochastically increasing, so that this is possible by the fact
that random variable~$A$ is stochastically smaller than $B$ if and only
if we can couple\vspace*{1.5pt}
$A$ and $B$ to $(\hat{A}, \hat{B})$ such that $\hat{A}\leq\hat
{B},$ a.s.

Denote by ${\cal A}_m=\{R_m=m\}$ the event that the first artificial
stub is chosen after time $m$. Then, by Lemma~\ref{lem-first-art}, we
have that
$\prob({\cal A}_{\ulmn}^c)=o(1)$. Thus, by\vspace*{1pt} intersecting with
${\cal A}_{\ulmn}$ and its complement, and then using the Markov
inequality, we find
for any $c_n=o(\sqrt{\log{n}})$,
%
%e6.3 ###
%
\begin{eqnarray} \prob(|G_{\olmn}-\widetilde G_{\ulmn}|\geq c_n)
&\leq&
\frac{1}{c_n} \expec[|G_{\olmn}-\widetilde G_{\ulmn}|\indicwo
{{\cal A}_{\ulmn}} ]+o(1)\nonumber\\
&=&\frac{1}{c_n} \expec
\bigl[[G_{\olmn}-\widetilde G_{\ulmn}]\indicwo{{\cal A}_{\ulmn}}
\bigr]+o(1)\\
&=& \frac{1}{c_n}\sum_{i=\ulmn+1}^{\olmn} \expec
\biggl[\frac{B_i}{S_i} \indicwo{{\cal A}_{\ulmn}} \biggr]+o(1).
\nonumber
\end{eqnarray}
We claim that
%
%e6.4 ###
%
\begin{equation} \label{mean-BiSi} \sum_{i=\ulmn+1}^{\olmn}\expec
\biggl[\frac{B_i}{S_i}\indicwo{{\cal A}_{\ulmn}} \biggr] =o\bigl
(\sqrt
{\log{n}}\bigr).
\end{equation}
Indeed, to see~(\ref{mean-BiSi}), we note that $B_i=0,$ when
$i\neq R_j$ for some $j$. Also, when ${\cal A}_{\ulmn}$
occurs, then $R_{\ulmn}=\ulmn$. Thus, using also that $R_m\geq m$, so that
$R_i\le\olmn$ implies that $i\leq\olmn$,
\begin{eqnarray}
\sum_{i=\ulmn+1}^{\olmn}\expec\biggl[\frac{B_i}{S_i}\indicwo
{{\cal A}_{\ulmn}} \biggr]
&\leq&\sum_{i=\ulmn+1}^{\olmn} \expec\biggl[\frac{B_{R_i}}{S_{R_i}}
\indicwo{\{\ulmn+1\le R_i\le\olmn\}} \biggr]\nonumber\\[-8pt]\\[-8pt]
&\leq&\sum_{i=\ulmn+1}^{\olmn}\frac{1}{i-1}
\expec\biggl[\frac{S_{R_i}+R_i}{S_{R_i}}\indicwo{\{\ulmn+1\le
R_i\le
\olmn\}} \biggr],\nonumber
\end{eqnarray}
the latter following from the exchangeability of $\{B_{R_i}\}
_{i=2}^{n-1}$, because
\[
S_{R_i}=D_1+\sum_{j=2}^{R_i}(B_j-1)=D_1+ \sum_{j=2}^i B_{R_j}-(R_i-1),
\]
so that
%
%e6.5 ###
%
\begin{equation} \sum_{j=2}^i B_{R_j}= S_{R_i}-D_1+R_i-1\leq
S_{R_i}+R_i.
\end{equation}
In Lemma~\ref{boundexpS1} of the Appendix~\hyperref[sec-app-A]{A} we
show that there exists
a constant~$C$ such that for $i\le\olmn$,
%
%e6.6 ###
%
\begin{equation} \label{Sratio-bd} \expec\biggl[\frac
{S_{R_i}+R_i}{S_{R_i}}\indicwo{\{\ulmn+1\le R_i\le\olmn\}}
\biggr]\le C,
\end{equation}
so that, for an appropriate chosen $c_n$ with $\log{(\olmn/\ulmn
)}/c_n\to0$,
%
%e6.7 ###
%
\begin{equation} \prob(|G_{\olmn}-\widetilde G_{\ulmn}|\geq
c_n ) \leq\frac{C}{c_n}\sum_{i=\ulmn+1}^{\olmn}\frac{1}{i-1}
\leq\frac{C\log{(\olmn/\ulmn)}}{c_n}=o(1),
\end{equation}
since $\log{(\olmn/\ulmn)}=o(\sqrt{\log{n}})$.
Thus, the CLT for $G_{\olmn}$ follows from the one  for
$\widetilde G_{\ulmn}$ which, since the marginal of $\widetilde
G_{\ulmn}$
is the same as the one of $G_{\ulmn}$, follows from the one for
$G_{\ulmn}$. By Lemma~\ref{lem-first-art}, we further have that
with high probability, there has not been any artificial stub up
to time $\ulmn$, so that, again with high probability,
$\{B_{m}\}_{m=2}^{\ulmn}=\{B_{R_m}\}_{m=2}^{\ulmn}$,
the latter, by Lemma~\ref{lem-exchange}, being an exchangeable
sequence.\looseness=1

We next adapt the proof of Proposition~\ref{lemma:CLT-sum1}
to exchangeable sequences under certain conditions. We start with $\tau
\in(2,3)$,
which is relatively the more simple case.
Recall the definition of $G_m$ in \eqref{def-Gm1}.
We define, for $i\geq2$,
%
%e6.8 ###
%
\begin{equation} \label{hatS-def} \hatS_i=\sum_{j=2}^i
B_j=S_i+i-1-D_1.
\end{equation}
Similarly to the proof of Proposition~\ref{lemma:CLT-sum1} we now
introduce
%
%e6.9 ###
%
\begin{equation} \label{hatG-def}\qquad\hspace*{10pt} \hatG_m=1+\sum_{i=2}^m \hatI
_i, \qquad\mbox{where\ } \prob(\hatI_i=1|\{B_i\}
_{i=2}^m)=B_i/\hatS_i,\qquad2 \leq i \leq m.
\end{equation}
Let $\hat{Q}_i=B_i/\hatS_i, Q_i=B_i/S_i$. Then, by a standard coupling
argument, we can couple $\hatI_i$ and $I_i$ in such a way that $\prob
(\hatI_i\neq I_i|\{B_i\}_{i=2}^m)=|\hat{Q}_i-Q_i|.$

The CLT for $\hatG_m$ follows because, also in the exchangeable setting,
$\hatI_2,\ldots,\hatI_m$ are independent and, similar to~\eqref{bd-GmhatGm},
%
%e6.10 ###
%
\begin{eqnarray} &&\prob(|G_m-\hatG_m|\geq\kappa_n )\nonumber\\
&&\qquad\leq
\kappa_n^{-1} \expec[|G_m-\hatG_m|] \leq\kappa_n^{-1} \expec
\Biggl[\sum
_{i=1}^m |I_i-\hatI_i|\Biggr]\nonumber\\
&&\qquad=\kappa_n^{-1} \sum_{i=2}^{m} \expec[|\hat
{Q}_i-Q_i|]\nonumber\\
&&\qquad=\kappa_n^{-1} \sum_{i=2}^{m} \expec
\biggl[B_i\frac{|S_i-\hatS_i|}{S_i\hatS_i} \biggr]\nonumber\\ &&\qquad\leq\kappa
_n^{-1} \sum
_{i=2}^{m} \expec\biggl[B_i\frac{D_1+(i-1)}{S_i\hatS_i} \biggr]\nonumber\\
&&\qquad=\kappa_n^{-1} \sum_{i=2}^{m} \frac{1}{i-1}\expec\biggl[\frac
{D_1+(i-1)}{S_i} \biggr]\\
&&\qquad=\kappa_n^{-1} \sum_{i=2}^{m} \biggl(\expec
[1/S_i]+\frac{1}{i-1}\expec[D_1/S_i] \biggr)\nonumber\\
&&\qquad\leq\kappa
_n^{-1} \sum_{i=2}^{m} \biggl(\expec[1/(S_i-D_1+2)]+\frac
{1}{i-1}\expec[D_1/(S_i-D_1+2)] \biggr),\nonumber
\end{eqnarray}
where we used that $D_1\geq2$ a.s. We take $m=\ulmn$, as discussed
above. Since $D_1$
is independent of $S_i-D_1+2$ for $i\geq2$ and $\expec[D_1]<\infty$, we
obtain the CLT for $G_{\ulmn}$ from the one
for $\hatG_{\ulmn}$ when, for $\tau\in(2,3)$,
%
%e6.11 ###
%
\begin{equation} \label{cond-CLT-(2,3)}\quad\sum_{i=1}^{\ulmn}
\expec
[1/\Sigma_i]=O(1),\qquad \mbox{where } \Sigma_i=1+\sum
_{j=2}^i (B_j-1),\qquad i\ge1.
\end{equation}
In Lemma~\ref{boundexpS1} of the Appendix~\hyperref[sec-app-A]{A} we
will prove that for $\tau
\in(2,3)$, the
statement~\eqref{cond-CLT-(2,3)} holds. The CLT for $G_{R_{\olmn}}$
follows in an
identical way.

We continue by studying the distribution of $T_m$ and $\widetilde T_m$,
for $\tau\in(2,3)$.
We recall that $T_m=\sum_{i=1}^m E_i/S_i$ [see \eqref{Tm-def}]. In
the proof of
Proposition~\ref{lemma:CLT-sum1}(b) for $\tau\in(2,3)$, we have made
crucial use of
(\ref{bound-mean-ratioS-indep}), which is now replaced by~(\ref
{cond-CLT-(2,3)}).
We split
%
%e6.12 ###
%
\begin{equation} \label{Tm-split-indep} T_m=\sum_{i=1}^m E_i/S_i=\sum
_{i=1}^{n^\rho} E_i/S_i+\sum_{i>n^\rho}^m E_i/S_i.
\end{equation}
The mean of the second term converges to 0 for each $\rho>0$ by Lemma
\ref{boundexpS1}, while the first
term is by Proposition~\ref{prop-indep} w.h.p. equal to $\sum
_{i=1}^{n^\rho} E_i/\Sindep_i$,
where $\Sindep_i=\sum_{j=1}^i \Bindep_j$, and where $\Bindep
_1=D_1$, while $\{\Bindep_i\}_{i=2}^{n^\rho}$
is an i.i.d. sequence of random variables with probability mass
function $g$ given in
(\ref{eqn:size-bias}), which is independent from $D_1$.
Thus, noting that also $\sum_{i>n^\rho}^m E_i/\Sindep_i\convp0$,
and with
%
%e6.13 ###
%
\begin{equation} \label{X-def} X=\sum_{i=1}^{\infty} E_i/\Sindep_i,
\end{equation}
we obtain that $T_m\convd X$. The random variable $X$ has the
interpretation of
the explosion time of the continuous-time branching process, where the
degree of
the root has distribution function $F$, while
the degrees of the other vertices is an i.i.d. sequence of random
variables with probability mass\vspace*{-1pt} function $g$ given in
(\ref{eqn:size-bias}). %A~similar argument holds for $\widetilde T_m$,
%with $\widetilde X\stackrel{d}{=}X$.
This completes the proof of Proposition~\ref{lemma:CLT-sum2} for $\tau
\in(2,3)$, and we turn
to the case $\tau>3$.

For $\tau>3$, we follow the steps in the proof of Proposition~\ref
{lemma:CLT-sum1}(a)
for $\tau>3$ as closely as possible.
Again, we apply a conditional CLT as in~(\ref{clt-LLF}), to obtain the CLT
when~(\ref{LLF-cond-CLT}) holds. From Lemma~\ref{lem-reciproke} we
conclude that
(\ref{Sratio-bd}) also holds when $\tau>3$. Hence, as before, we may assume
by Lemma~\ref{lem-first-art}, that w.h.p., there has not been any
artificial stub up
to time $\ulmn$, so that, again w.h.p.,
$\{B_{m}\}_{m=2}^{\ulmn}=\{B_{R_m}\}_{m=2}^{\ulmn}$,
the latter, by Lemma~\ref{lem-exchange}, being an exchangeable sequence.
For the exchangeable sequence $\{B_{m}\}_{m=2}^{\ulmn}$ we will then
show that
%
%e6.14 ###
%
\begin{equation} \label{sec-mom-sum-exch} \sum_{j=2}^{\ulmn}
B_j^2/S_j^2=\Op(1).
\end{equation}
The statement~(\ref{sec-mom-sum-exch}) is proven in Lemma~\ref
{lem-bd-squareBSratio}.

As in the proof of Proposition~\ref{lemma:CLT-sum1}(a), the claim that
%
%e6.15 ###
%
\begin{equation}
\label{suff-cond2}
\sum_{j=2}^{\ulmn} B_j/S_j-\frac{\nu}{\nu-1}\log{\ulmn}=o_{
\prob}\bigl(\sqrt{\log{\ulmn}}\bigr)
\end{equation}
is sufficient for the CLT when $\tau>3$. Moreover, we have shown in Remark
\ref{rem-CLT-tau>3} that
\eqref{suff-cond2} is satisfied when
%
%e6.16 ###
%
\begin{equation}
\label{SL-excha}
\log{(S_{\ulmn}/\ulmn)}=
o_{\prob}\bigl(\sqrt{\log{\ulmn}}\bigr)
\end{equation}
and
%
%e6.17 ###
%
\begin{equation}
\label{lln(b)-excha}
\sum_{j=1}^{\ulmn} \frac{S_j-(\nu-1) j}{S_j(\nu-1)j}=o_{\prob
}\bigl(\sqrt{\log{\ulmn}}\bigr).
\end{equation}
The proofs of \eqref{SL-excha} and \eqref{lln(b)-excha} are given in
Lemmas~\ref{lemma-lgS}
and~\ref{Mar-Zyg} of Appendix~\hyperref[sec-app-A]{A}, respectively.
Again, the proof
for $G_{R_{\olmn}}$ is identical.
\end{pf*}

For the results for $T_m$ and $\widetilde T_m$ for $\tau>3$, we refer to
Appendix~\hyperref[sec-weak-conv-weight-tau3]{C}.

%s7 ###
\section[Proof of Proposition 4.9]{Proof of Proposition \protect\ref{lemma:conn_edge}}
\label{sec-lemma:conn_edge}
In this section, we prove Proposition~\ref{lemma:conn_edge}.
We start by proving that $\log{\CE_n/a_n}=o_{\prob}(\sqrt{\log{n}})$,
where $\CE_n$ is the time at which the connecting edge appears between
the SWGs of vertices 1 and 2 [recall~(\ref{conn-edge-time})], as
stated in~(\ref{time-conn-edge}).
As described in Section~\ref{sec-connedge},
we shall condition vertices 1 and 2 to be in the giant component,
which occurs w.h.p. and guarantees that $S^{(i)}_{m}>0$
for any $m=o(n)$ and $i\in\{1,2\}$. After this, we complete
the proof of~(\ref{clt-conn-edge})--(\ref{weight-conn-edge})
in the case where $\tau\in(2,3)$, which turns out to be relatively
simplest, followed by a proof of~(\ref{clt-conn-edge}) for $\tau>3$.
The proof of~(\ref{weight-conn-edge}) for $\tau>3$, which is more
delicate, is deferred to Appendix~\hyperref[sec-weak-conv-weight-tau3]{C}.

We start by identifying the distribution of $\CE_n$.
In order for $\CE_n=m$ to occur,
apart from further requirements, the minimal stub from $\SWG^{(2)}_{m}$
must be real, that is, it may not be artificial.
This occurs with probability equal to\vspace*{1pt}
$1-|\Art^{(2)}_{m}|/S^{(2)}_{m}$.

By Construction~\ref{const-SWGCM}, the number of allowed stubs
incident to
$\SWG^{(2)}_{m}$ equals $S^{(2)}_m$, so the number
of real stubs equals $S^{(2)}_m-|\Art_m^{(2)}|$.
Similarly, the number of allowed stubs incident to
$\SWG^{(1)}_{a_n}$ equals $S^{(2)}_m$, so the number
of real\vspace*{-1pt} stubs equals $S^{(1)}_{a_n}-|\Art_{a_n}^{(1)}|$.
Further, the number of free stubs equals $|\FS_m|=L_n-a_n-m-S_m+|\Art_m|$,
where we recall that $S_m=S^{(1)}_{a_n}+S^{(2)}_m$ and
$\Art_m=\Art_{a_n}^{(1)}\cup\Art_{m}^{(2)}$,
and is hence bounded above by $L_n$ and below by
$L_n-a_n-m-S_m$. When the minimal-weight stub is indeed
real, then it must be attached to one of the real
allowed stubs incident to $\SWG^{(1)}_{a_n}$, which
occurs with conditional probability given $\SWG^{(1,2)}_m$ and~$L_n$
equal to
%
%e7.1 ###
%
\begin{equation} \label{CEn-comp-rep} \frac{S_{a_n}^{(1)}-|\Art
_{a_n}^{(1)}|} {L_n-a_n-m-S_m+|\Art_m|}.
\end{equation}
Thus, in order to prove
Proposition~\ref{lemma:conn_edge}, it suffices to investigate
the\vspace*{1pt} limiting behavior of $L_n$, $S_m^{(i)}$ and $|\Art_m|$.
By the law of large numbers, we known that
$L_n-\mu n=o_{\prob}(n)$ as $n\rightarrow\infty$.
To study $S_m^{(i)}$ and $|\Art_m|$, we shall make use of results from
\cite{hofs3,hofs1}. Note that we can write
$S_m^{(i)}=D_i+B_2^{(i)}+\cdots+B_m^{(i)}-(m-1),$
where $\{B_m^{(i)}\}_{m=2}^{\infty}$ are close to
being independent. See \cite{hofs3}, Lemma~A.2.8,
for stochastic domination results on $\{B_m^{(i)}\}_{m=2}^{\infty}$
and their sums in terms of i.i.d. random variables,
which can be applied in the case of $\tau>3$.
See \cite{hofs1}, Lemma~A.1.4, for bounds on tail
probabilities for sums and maxima of random variables
with certain tail properties.

The next step to be performed is to give criteria in terms of the
processes $S_m^{(i)}$ which guarantee that the estimates
in Proposition~\ref{lemma:conn_edge} follow. We shall start
by proving that with high probability $\CE_n\geq\ulmn$, where
$\ulmn=\vep_na_n$, where $\vep_n\downarrow0$. This proof makes use
of, and
is quite similar to, the proof of Lemma~\ref{lem-first-art} given in
Appendix~\hyperref[sec-app-A]{A}.

\begin{Lemma}[(Lower bound on time to connection)]
\label{lem-CEn-lb}
Let $\ulmn/a_n=o(1)$. Then
%
%e7.2 ###
%
\begin{equation} \prob(\CE_n\leq\ulmn)=o(1).
\end{equation}
\end{Lemma}

\begin{pf}
%We further write
% \eqn{
% {\cal A}_{m}={\cal A}_{m}^{\sss(1)}\cap{\cal A}_{m}^{\sss(2)}.
% }
Denote
%
%e7.3 ###
%
\begin{equation} {\cal E}_n=\bigl\{S_{a_n}^{(1)} \leq a_n M_n\bigr\},
\end{equation}
where $M_n=\eta_n^{-1}$ for $\tau>3$,
while $M_n=\eta_n^{-1}n^{(3-\tau)/(\tau-1)}$ for $\tau\in(2,3)$,
and where $\eta_n\downarrow0$
sufficiently slowly. Then, by~(\ref{bd-mean-Sm-tau>3}) for $\tau>3$ and
(\ref{aim-lb-Sm(2,3)-ub}) for $\tau\in(2,3)$,
%
%e7.4 ###
%
\begin{equation} \label{Encompl} \prob({\cal E}_n^c)=o(1),
\end{equation}
since $\ulmn=o(a_n)$.
%Thus, in particular,
%when ${\cal A}_{\ulmn}$ occurs, then, for all $m\leq\ulmn$,
% \eqn{
% \lbeq{CEn-comp2}
% \Qprob^{\sss(m)}_n(\CE_n=m|\CE_n>m-1)=\frac{S_{\lfloor m/2\rfloor}^{
% {L_n-m-S_m}, i_m=1+(m \mbox{mod} 2),
% }
%where we write $\Qprob^{\sss(m)}_n$ for the conditional distribution
%given $\SWG^{\sss(1,2)}_m$ and $\{D_i\}_{i=1}^n$.
%We start with
% \eqan{
% \prob(\CE_n\leq m_n)&=\sum_{m=2}^{\ulmn} \prob(\CE_n=m|\CE_n>m-1)
% &\leq\sum_{m=2}^{m_n} \prob(\CE_n=m|\CE_n>m-1)+o(1).
% }
By the law of total probability,
%
%e7.5 ###
%
\begin{equation}\qquad\hspace*{6pt}
\prob(\CE_n\leq\ulmn)\leq\prob({\cal E}_n^c)+\sum
_{m=1}^{\ulmn} \prob(\{\CE_n=m\}\cap{\cal E}_n|\CE_n>m-1)\prob
(\CE_n>m-1).
\end{equation}
Then, we make use of~(\ref{CEn-comp-rep}) and~(\ref{Encompl}), to
arrive at
%
%e7.6 ###
%
\begin{equation} \label{CEn-lb} \prob(\CE_n\leq\ulmn) \leq\sum
_{m=1}^{\ulmn} \expec\biggl[\indicwo{{\cal E}_n}\frac{S_{a_n}^{
(1)}} {L_n-m-S_m} \biggr]+o(1).
\end{equation}
As in the proof of Lemma~\ref{lem-first-art}, we have that $m\le\ulmn
=o(n)$ and $S_m=o(n)$,
while $L_n\geq n$. Thus,~(\ref{CEn-lb}) can be simplified to
%
%e7.7 ###
%
\begin{equation} \label{CEn-lb1} \prob(\CE_n\leq\ulmn)\leq\frac
{1+o(1)}{n}\ulmn\expec \bigl[\indicwo{{\cal E}_n}S_{a_n}^{
(1)} \bigr]+o(1) \leq\frac{1+o(1)}{n\eta_n}a_n\ulmn.
\end{equation}
When choosing $\eta_n\downarrow0$ sufficiently slowly,
for example as $\eta_n=\sqrt{\ulmn/a_n}$, we obtain that $\prob(\CE
_n\leq\ulmn)=o(1)$
whenever $\ulmn=o(a_n).$
\end{pf}

We next state an upper bound on $\CE_n$.

\begin{Lemma}[(Upper bound on time to connection)]
\label{lem-CEn-ub}
Let $\olmn/ a_n\to\infty$, then,
%
%e7.8 ###
%
\begin{equation} \prob(\CE_n>\olmn)=o(1).
\end{equation}
\end{Lemma}

\begin{pf}
We start by giving an explicit formula for $\prob(\CE_n>m)$. As before,
$\Qprob^{(m)}_n$ is the conditional distribution given $\SWG
^{(1,2)}_m$ and
$\{D_i\}_{i=1}^n$. Then, by
Lemma~\ref{lem-prod-CEn},
%
%e7.9 ###
%
\begin{equation} \label{prod-CEn1} \prob(\CE_n>m)=\expec\Biggl
[\prod
_{j=1}^m \Qprob^{(j)}_n(\CE_n>j|\CE_n>j-1) \Biggr].
\end{equation}
Equation~(\ref{prod-CEn1}) is identical in spirit to \cite{hofs3},
Lemma 4.1,
where a similar identity was used for the graph distance in the CM.
Now, for any sequence $\vep_n\rightarrow0$, let
%
%e7.10 ###
%
\begin{equation} \label{Bcal-n-def} {\cal B}_n= \Biggl\{\frac{c\olmn
}{n} \big|\Art_{a_n}^{(1)}\big| +\frac{c}{n}S_{a_n}^{(1)}\sum
_{m=1}^{\olmn} \frac{|\Art_{m}^{(2)}|}{S_{m}^{(2)}}\leq
\vep_n \Biggr\}.
\end{equation}

By Lemma~\ref{cor-conv-zero-error}, the two terms appearing in the
definition of
${\cal B}_n$ in~(\ref{Bcal-n-def}) converge to zero in probability, so that
$\prob({\cal B}_n)=1-o(1)$ for some $\vep_n\rightarrow0$. Then we bound
%
%e7.11 ###
%
\begin{equation} \label{prod-CEn1b} \prob(\CE_n>m)\leq\expec
\Biggl[\indicwo{{\cal B}_n}\prod_{j=1}^m \Qprob^{(j)}_n(\CE
_n>j|\CE
_n>j-1) \Biggr] +\prob({\cal B}_n^c).
\end{equation}
We continue by noticing that according to~(\ref{CEn-comp-rep}),
%
%e7.12 ###
%
\begin{equation} \label{Qprob-m} \Qprob^{(m)}_n(\CE_n=m|\CE
_n>m-1)=\frac{S_{a_n}^{(1)}-|\Art_{a_n}^{(1)}|} {|\FS
_m|} \biggl(1-\frac{|\Art_{m}^{(2)}|}{S_{m}^{(2)}} \biggr),
\end{equation}
where $|\FS_m|$ is the number of real free stubs which is available at
time $m$.
Combining~(\ref{prod-CEn1b}) and~(\ref{Qprob-m}) we arrive at
%
%e7.13 ###
%
\begin{equation}\qquad\quad\prob(\CE_n>\olmn)= \expec\Biggl
[\indicwo{{\cal
B}_n}\prod_{m=1}^{\olmn} \biggl(1-\frac{S_{a_n}^{(1)}- |\Art
_{a_n}^{(1)}|}{|\FS_m|} \biggl(1-\frac{|\Art_{m}^{(2)}|}
{S_{m}^{(2)}} \biggr) \biggr) \Biggr]+o(1).
\end{equation}
Since $|\FS_m|\leq L_n\leq n/c$, w.h.p., for some $c>0,$ and using
that $1-x\leq{\mathrm e}^{-x}$, we can further bound
\begin{eqnarray} \prob(\CE_n>\olmn) &\leq& \expec\Biggl[\indicwo
{{\cal
B}_n}\exp\Biggl\{-\frac{c}{n} \bigl(S_{a_n}^{(1)}- \big|\Art
_{a_n}^{(1)}\big| \bigr)\sum_{m=1}^{\olmn} \biggl(1-\frac{|\Art
_{m}^{(2)}|} {S_{m}^{(2)}} \biggr) \Biggr\} \Biggr]+o(1)
\nonumber\\[-8pt]\\[-8pt]
&\leq&\expec\biggl[\indicwo{{\cal B}_n}\exp\biggl\{-\frac{c\olmn
}{n}S_{a_n}^{(1)} \biggr\} \biggr]+e_n+o(1),\nonumber
\end{eqnarray}
where
%
%e7.14 ###
%
\begin{equation} \label{twee-termen}
\qquad
e_n =O \Biggl(\expec\Biggl[\indicwo
{{\cal B}_n} \Biggl(\frac{c\olmn}{n}\big|\Art_{a_n}^{(1)}\big|
+\frac
{c}{n}S_{a_n}^{(1)}\sum_{m=1}^{\olmn} \frac{|\Art_{m}^{
(2)}|}{S_{m}^{(2)}} \Biggr) \Biggr] \Biggr)=O(\vep_n).
\end{equation}
Hence,
%
%e7.15 ###
%
\begin{equation} \prob(\CE_n>\olmn)\leq\expec\biggl[\exp\biggl\{
-\frac{c\olmn}{n} S_{a_n}^{(1)} \biggr\} \biggr]+o(1).
\end{equation}
When $\tau>3$, by Lemma~\ref{lem-cond-LD-est} in the
Appendix~\hyperref[sec-app-A]{A},
we have that, w.h.p., and for some
$\eta>0$,
%
%e7.16 ###
%
\begin{equation} \label{aim-lb-Sm>3} S_{a_n}\geq\eta a_n,
\end{equation}
so that
%
%e7.17 ###
%
\begin{equation} \prob(\CE_n>\olmn) \leq\exp\biggl\{-\frac{c\eta
a_n \olmn}{n} \biggr\}+o(1)=o(1),
\end{equation}
as long as $\olmn a_n/n=\olmn/\sqrt{n}\to\infty$.
For $\tau\in(2,3)$, by~(\ref{aim-lb-Sm(2,3)-lb}) in Lemma~\ref
{lem-bdsS-(2,3)}, and using
that $n^{1/(\tau-1)}/n=1/a_n$, we have
for every $\vep_n\to0$,
%
%e7.18 ###
%
\begin{equation} \prob(\CE_n>\olmn) \leq\exp\biggl\{-\frac
{c\olmn\vep
_n}{a_n}\biggr\}+o(1)=o(1),
\end{equation}
whenever $\vep_n \olmn/a_n\rightarrow\infty$. By adjusting $\vep
_n$, it is hence sufficient
to assume that $\olmn/a_n\to\infty$.
\end{pf}

Lemmas~\ref{lem-CEn-lb} and~\ref{lem-CEn-ub} complete the proof of
(\ref{time-conn-edge}) in Proposition~\ref{lemma:conn_edge}. We
next continue with the proof of~(\ref{clt-conn-edge}) in
Proposition~\ref{lemma:conn_edge}. We start by showing that
$\prob(\CE_n=0)=o(1)$. Indeed, $\CE_n=0$ happens precisely when
$2\in\SWG_{a_n}^{(1)}$, which, by exchangeability,
occurs with probability at most $a_n/n=o(1)$.

For $\CE_n\geq1$, we note that at time $\CE_n$, we draw a real stub.
%Thus, the random time $\CE_n=R_{m}^{\sss(2)}$ for a certain $m$.
Consider the pair $(G_{a_n+1}^{(1)},G_{\CE_n}^{(2)})$
conditionally on $\{\CE_n=m\}$ for a certain $m$. The event $\{\CE
_n=m\}$
is \textit{equal} to the event that the last chosen stub in
$\SWG^{(2)}_{m}$ is paired to a stub incident to
$\SWG^{(1)}_{a_n}$, while this is not the case for
all previously chosen\vspace*{-2pt} stubs. For $j=1, \ldots, m$, and $i\in\{1,2\}$,
denote by $I_{j}^{(i)}$ the $j$th real vertex\vspace*{-2pt} added to $\SWG
^{(i)}$,
and denote by $V_m^{(i)}$ the number
of real vertices in $\SWG_m^{(i)}$. Then, for $m\geq1$, the
event $\{\CE_n=m\}$
is \textit{equal} to the event that the last chosen stub in
$\SWG^{(2)}_{m}$ is paired to a stub incident to
$\SWG^{(1)}_{a_n}$, and
%
%e7.19 ###
%
\begin{equation} \bigl\{I_j^{(1)}\bigr\}_{j=1}^{V_{a_n}^{
(1)}}\cap\bigl\{I_j^{(2)}\bigr\}_{j=1}^{V_{m}^{(2)}}=\varnothing
.
\end{equation}
As a result, conditionally on $\{\CE_n=m\}$ and $V_{a_n}^{
(1)}=k_1, V_{m}^{(2)}=k_2$,
the vector consisting of both $\{B_{R_{j_1}}^{(1)}\}
_{j_1=1}^{k_1}$ and
$\{B_{R_{j_2}}^{(2)}\}_{j_2=1}^{k_2}$ is an exchangeable vector,
with law
is equal to that of $k_1+k_2$ draws from $\{D_i-1\}_{i=3}^n$ without
replacement, where, for $i\in[n]\setminus\{1,2\}$, $D_i-1$ is drawn
with probability
equal to $D_i/(L_n-D_1-D_2)$. The above explains the role of the random
stopping time
$\CE_n$.

We continue by discussing the limiting distributions of
$(H_n^{(1)},H_n^{(2)})$ in order
to prove~(\ref{clt-conn-edge}). For this, we note that if we condition on
$\{\CE_n=m\}$ for some $m$ and on $\SWG_m^{(1,2)}$, then, by
(\ref{tildeHn-def})
$(H_n^{(1)},H_n^{(2)})\stackrel{d}{=}(G_{a_n+1}^{
(1)},G_{m}^{(2)})$,
where\vspace*{-1pt} the \textit{conditional} distribution of
$(G_{a_n+1}^{(1)},G_{m}^{(2)})$
is as two \textit{independent}\vspace*{-1pt} copies of $G$ as
described in~(\ref
{Gm-def}), where
$\{d_j\}_{j=1}^{a_n}$ in~(\ref{Gm-def})
is given by $d_1=D_1$ and $d_j=B^{(1)}_j, j\ge2$,\vspace*{-1pt} while,
$H_n^{(2)}=G^{(2)}_{m}$,
where $d_1=D_2$ and $d_j=B^{(2)}_j, j\ge2$. Here, we make\vspace
*{-1pt} use
of the fact
that $H_n^{(1)}$ is the distance from vertex 1 to the vertex to which
the paired stub is connected to, which has the same distribution as
the distance from vertex 1 to the vertex which has been added at time
$a_n+1$, minus 1, since the paired stub is again a
\textit{uniform} stub (conditioned to be real).

Thus, any possible dependence of $(H_n^{(1)},H_n^{(2)})$
arises through the dependence of the vectors
$\{B_{j}^{(1)}\}_{j=2}^{\infty}$ and $\{B_{j}^{(2)}\}
_{j=2}^{\infty}$.
However, the proof of Proposition~\ref{lemma:CLT-sum2} shows that
certain weak dependency of
$\{B_{j}^{(1)}\}_{j=2}^{\infty}$ and $\{B_{j}^{(2)}\}
_{j=2}^{\infty}$
is allowed.\vspace*{1pt}

We start by completing the proof for $\tau\in(2,3)$ which is the
more simple one.
Recall the split in~(\ref{hatGm-def}), which was fundamental in
showing the CLT
for $\tau\in(2,3)$. Indeed, let $\{\hatI_{j}^{(1)}\}
_{j=1}^{\infty}$
and $\{\hatI_{j}^{(2)}\}_{j=1}^{\infty}$ be two sequences of
indicators, with $\hatI_{1}^{(1)}=\hatI_{1}^{(2)}=1$,
which are, conditionally on $\{B_{j}^{(1)}\}_{j=2}^{\infty}$ and
$\{B_{j}^{(2)}\}_{j=2}^{\infty}$,
independent with, for $i\in\{1,2\}$,
%
%e7.20 ###
%
\begin{equation} \label{hatI-joint-def} \prob\bigl(\hatI_{j}^{
(i)}=1|\bigl\{B_{j}^{(i)}\bigr\}_{j=2}^{\infty} \bigr)=B_{j}^{(i)}
/\bigl(S_{j}^{(i)}+j-1-D_i\bigr).
\end{equation}
Then, the argument in~(\ref{ind-hatI}) can be straightforwardly
adapted to show that
the unconditional distributions of $\{\hatI_{j}^{(1)}\}
_{j=2}^{\infty}$
and $\{\hatI_{j}^{(2)}\}_{j=2}^{\infty}$ are that of\vspace*{1pt} two
independent sequences
$\{J_{j}^{(1)}\}_{j=2}^{\infty}$
and $\{J_{j}^{(2)}\}_{j=2}^{\infty}$ with $\prob(J_{j}^{
(i)}=1)=1/(j-1)$.
Thus,\vspace*{1pt} by the independence, we immediately obtain that since $\CE
_n\rightarrow\infty$ with $\log(\CE_n/a_n)=o_{\prob}(\sqrt
{\log{n}})$,
%
%e7.21 ###
%
\begin{equation} \biggl(\frac{H^{(1)}_n-\beta\log{a_n}}{\sqrt
{\beta\log{a_n}}}, \frac{H^{(2)}_{n}-\beta\log{a_n}}{\sqrt
{\beta\log{a_n}}} \biggr) \convd(Z_1,Z_2).
\end{equation}
The argument to show that, since $\CE_n\leq\olmn$,
$(H_n^{(1)},H_n^{(2)})$
can be well\vspace*{1pt} approximated by
$(G_{\ulmn}^{(1)},G_{\ulmn}^{(2)})$
[recall~(\ref{Gm-split})]
only depends on the marginals of $(H_n^{(1)},H_n^{(2)})$,
and thus remains valid verbatim. We conclude that~(\ref
{clt-conn-edge}) holds.

We next prove~(\ref{weight-conn-edge}) for $\tau\in(2,3)$. For
this, we again use
Proposition~\ref{prop-indep} to note that the forward degrees $\{B_j\}
_{j=3}^{n^\rho}$
can be coupled to i.i.d.\vspace*{-2pt} random variables $\{\Bindep_j\}
_{j=3}^{n^\rho}$, which are independent
from $B_1=D_1, B_2=D_2$. Then we can follow the proof of Proposition
\ref{lemma:CLT-sum2}(b) for $\tau\in(2,3)$
verbatim,\vspace*{-1pt} to obtain that $(T^{(1)}_{a_n},T^{(2)}_{\CE
_n})\convd(X_1,X_2)$, where
$X_1, X_2$ are two\vspace*{1pt} independent copies of $X$ in~(\ref
{X-def}). This
completes the proof of
Proposition~\ref{lemma:conn_edge} when $\tau\in(2,3)$.

We proceed with the proof of Proposition~\ref{lemma:conn_edge}
when $\tau>3$ by
studying $(H_n^{(1)},H_n^{(2)})$.
We follow the proof of Proposition~\ref{lemma:CLT-sum2}(a), paying
particular attention to
the claimed independence of the limits $(Z_1,Z_2)$ in~(\ref
{clt-conn-edge}). The proof of
Proposition~\ref{lemma:CLT-sum2}(a) is based on a \textit{conditional}
CLT, applying
the Lindeberg--L\'evy--Feller condition. Thus, the conditional limits
$(Z_1,Z_2)$ of
\begin{eqnarray} &&\biggl(\frac{H_n^{(1)}-\sum_{j=2}^{a_n}
B_j^{(1)}/S_j^{(1)}} { (\sum_{j=2}^{a_n}
({B_j^{(1)}}/{S_j^{(1)}})(1-{B_j^{(1)}}/{S_j^{
(1)}}) )^{1/2}},\nonumber\\[-8pt]\\[-8pt]
&&\hspace*{9pt}\frac{H_n^{(2)}-\sum_{j=2}^{\CE_n}
B_j^{(2)}/S_j^{(2)}} { (\sum_{j=2}^{\CE_n}
({B_j^{(2)}}/{S_j^{(2)}})(1-{B_j^{(2)}}/{S_j^{
(2)}}) )^{1/2}} \biggr)\nonumber
\end{eqnarray}
are clearly independent. The proof then continues by showing that the
asymptotic mean and variance
can be replaced by $\beta\log{n}$, which is a computation based on
the marginals $\{B_j^{(1)}\}_{j=2}^{\infty}$
and $\{B_j^{(2)}\}_{j=2}^{\infty}$ only, and, thus, these results
carry over verbatim, when we further
make use of the fact that, w.h.p., $\CE_n\in[\ulmn,\olmn]$ for
any $\ulmn,\olmn$ such that
$\log{(\olmn/\ulmn)}=o(\sqrt{\log{n}})$. This completes the
proof of~(\ref{clt-conn-edge}) for $\tau>3$.
The proof of~(\ref{weight-conn-edge}) for $\tau>3$ is a bit more
involved, and is deferred to
Appendix~\hyperref[sec-weak-conv-weight-tau3]{C}.

\setcounter{Theorem}{0}
\setcounter{equation}{0}
\begin{appendix}
\section*{Appendix A: Auxiliary lemmas for CLTs in CM}\label{sec-app-A}

In this appendix, we denote by $B_1=D_1$, the degree of vertex $1$
and $B_2,\ldots,B_m$, $m<n,$ the forward degrees of the shortest
weight graph $\SWG_m$. The forward degree $B_k$ is chosen recursively from
the set $\FS_k$, the set of free stubs at time $k$. Further we
denote by
\[
S_k=D_1+\sum_{j=2}^k (B_j-1),
\]
the number of allowed stubs at time $k$. As before the random variable
$R_m$ denotes
the first time that the shortest path graph from vertex $1$ contains
$m+1$ real vertices.
Consequently
\[
B_{R_2},\ldots,B_{R_m},
\]
$m<n,$ can be seen as a sample without replacement from the degrees
\[
D_2-1,D_3-1,\ldots,D_n-1.
\]

%s7.1 ###
\subsection{The first artificial stub}
\label{sec-first-art-stub}
We often can and will replace the sample
$B_2,\ldots,B_{\ulmn},$ by the sample
$B_{R_2},\ldots,B_{R_{\ulmn}}$.
The two samples have, w.h.p., the same distribution if the first
artificial stub appears after time $\ulmn$. This will be the content
of our first lemma.

%%%%%%%%%%%%%%%%%%%%%%%%%%%%%%%%%%%%%%%%%%%%%%
%
\begin{Lemma}[(The first artificial stub)]
\label{lem-first-art}
Let $\ulmn/a_n\to0$.
Then,
%
%e7.1 ###
%
\begin{equation} \prob(R_{\ulmn}>\ulmn)=o(1).
\end{equation}
\end{Lemma}

\begin{pf} For the event $\{R_{\ulmn}>\ulmn\}$ to happen it is
mandatory that
for some $m\le\ulmn$, we have $R_m>m$, while $R_{m-1}=m-1$. Hence
%
%e7.2 ###
%
\begin{equation} \prob_n(R_{\ulmn}>\ulmn)= \sum_{m=2}^{\ulmn}
\prob_n(R_m>m, R_{m-1}=m-1).
\end{equation}
Now, when $R_m>m, R_{m-1}=m-1$, one of the $S_{m-1}$ stubs incident
to $\SWG_{m-1}$ has been drawn, so that
%
%e7.3 ###
%
\begin{equation}\qquad\prob_n(R_m>m,R_{m-1}=m-1) =\expec_n \biggl
[\frac
{S_{m-1}}{L_n-S_{m-1}-2m}\indic{R_{m-1}=m-1} \biggr].
\end{equation}
Since $\ulmn=o(n)$, we claim that, with high probability,
$S_{m-1}=o(n)$. Indeed, the maximal degree is $\Op(n^{1/(\tau-1)})$,
so that, for $m\leq\ulmn$,
%
%e7.4 ###
%
\begin{equation} \label{Sm-bd} S_{m}\leq\Op\bigl(m n^{1/(\tau
-1)}\bigr)\leq
\Op\bigl(\ulmn n^{1/(\tau-1)}\bigr)=o_{\prob}(n),
\end{equation}
since, for $\tau>3$, $a_n=n^{1/2}$ and $n^{1/(\tau-1)}=o(n^{1/2})$,
while, for $\tau\in(2,3)$, $a_n=n^{(\tau-2)/(\tau-1)}$, so that
$\ulmn n^{1/(\tau-1)}=o(n)$. Moreover, $L_n\geq n$, so that
%
%e7.5 ###
%
\begin{equation} \label{Artbd1} \prob_n(R_m>m, R_{m-1}=m-1) \leq
\frac{C}{n}\expec_n\bigl[S_{m-1}\indic{R_{m-1}=m-1}\bigr].
\end{equation}
By the remark preceding this lemma, since $R_{m-1}=m-1$, we have that
$S_{m-1}=D_1+\sum_{j=2}^{m-1} (B_{R_{j}}-1 )$, so that, by
Lemma~\ref{lem-exchange},
%
%e7.6 ###
%
\begin{equation} \label{Artbd2} \prob_n(R_m>m, R_{m-1}=m-1) \leq
\frac{C}{n}D_1+\frac{C(m-2)}{n}\expec_n[B_{R_2}].
\end{equation}
The first term converges to 0 in probability, while the expectation in the
second term, by~(\ref{BT-law}), equals
%
%e7.7 ###
%
\begin{equation} \label{BT-expec} \expec_n[B_{R_2}]=\sum_{i=2}^n
\frac{D_i(D_i-1)}{L_n-D_1}.
\end{equation}
When $\tau>3$, this has a bounded expectation, so that,
for $a_n=\sqrt{n}$,
%
%e7.8 ###
%
\begin{equation}
\label{Artbd3}
\quad\prob(R_{\ulmn}>\ulmn)\leq\sum_{m=2}^{\ulmn}
\frac{C}{n}\expec[D_1]+\sum_{m=2}^{\ulmn}\frac{C(m-2)}{n}\expec[B_{R_2}]
\leq C\frac{\ulmn^2}{n}\to0.
\end{equation}
When $\tau\in(2,3)$, however, then $\expec[D_i^2]=\infty$, and we
need to be a bit more careful.
In this case, we obtain from
(\ref{Artbd2}) that
%
%e7.9 ###
%
\begin{equation} \label{Artbd4} \prob_n(R_{\ulmn}>\ulmn) \leq
C\frac{\ulmn^2}{n}\expec_n[B_{R_2}].
\end{equation}
From \eqref{BT-expec}, and since $L_n-D_1\ge n-1$,
%
%e7.10 ###
%
\begin{equation} \label{Artbd5} \expec_n[B_{R_2}]\leq\frac{C}{n-1}
\sum_{i=2}^n D_i(D_i-1) \leq\frac{C}{n-1} \sum_{i=2}^n D_i^2.
\end{equation}
%
%Since the tail of $D_i$ is regularly varying with exponent $1-\tau$,
%we find
% $$
% \prob(D_i^2>x)=\prob(D_i>\sqrt{x})=x^{(1-\tau)/2}L(\sqrt{x}).
% $$
From~(\ref{Fcond(2,3)}), we obtain that
$x^{(\tau-1)/2}\prob(D_i^2>x)\in[c_1, c_2]$ uniformly in $x\geq0$,
and since $D_1,D_2,\ldots,D_n$ is i.i.d., we can conclude that
$n^{-2/(\tau-1)}\sum_{i=2}^n D_i^2$ converges
to a proper random variable. Hence, since
$a_n/\ulmn\rightarrow\infty$ we obtain, w.h.p.,
%
%e7.11 ###
%
\begin{equation} \label{Artbds} \expec_n[B_{R_2}]\leq\frac
{a_n}{\ulmn}n^{2/(\tau-1)-1}=\frac{a_n}{\ulmn}n^{(3-\tau)/(\tau
-1)}.
\end{equation}
Combining \eqref{Artbd4} and \eqref{Artbds}, and using that
$a_n=n^{(\tau-2)/(\tau-1)}$ we obtain that, w.h.p.,
%
%e7.12 ###
%
\begin{eqnarray} \label{Artbd6}\qquad \prob_n(R_{\ulmn}>\ulmn)
&\leq& Ca_n
\ulmn n^{{2}/({\tau-1})-1}\nonumber\\[-8pt]\\[-8pt]
 \qquad&=&C\frac{\ulmn}{a_n} n^{{2(\tau
-2)}/({\tau-1})+({3-\tau})/({\tau-1})-1}
= C\frac{\ulmn
}{a_n}=o_{\prob}(1).\nonumber
\end{eqnarray}
This proves the claim.
\end{pf}
%
%%%%%%%%%%%%%%%%%%%%%%%%%%%%%%%%%%%%%%%%%%%%%%%%%

%s7.2 ###
\subsection[Coupling the forward degrees to an i.i.d. sequence:
Proposition 4.7]{Coupling the forward degrees to an i.i.d. sequence:
Proposition \protect\ref{prop-indep}}
We will now prove Proposition \protect\ref{prop-indep}.
To this end, we denote the order statistics of the degrees by
%
%e7.13 ###
%
\begin{equation} \label{stoch-ord-1} D_{(1)}\leq D_{(2)}\leq
\cdots\leq D_{(n)}.
\end{equation}
Let $m_n\rightarrow\infty$ and consider the i.i.d. random variables
${\underline X_1},{\underline X_2},\ldots,
{\underline X_{m_n}}$, where
${\underline X_i}$ is taken \textit{with replacement} from the stubs
%
%e7.14 ###
%
\begin{equation} \label{stoch-ord-2} D_{(1)}-1,D_{
(2)}-1,\ldots, D_{(n-m_n)}-1,
\end{equation}
that is, we sample \textit{with replacement} from the original forward degrees
$D_1-1,D_2-1,\ldots,D_n-1$, where the $m_n$ largest degrees are discarded.
Similarly, we consider the i.i.d. random variables ${\overline
X_1},{\overline X_2},\ldots,
{\overline X_{m_n}}$, where\vadjust{\goodbreak}
${\overline X_i}$ is taken \textit{with replacement} from the stubs
%
%e7.15 ###
%
\begin{equation} \label{stoch-ord-3} D_{(m_n+1)}-1,D_{
(m_n+2)}-1,\ldots, D_{(n)}-1,
\end{equation}
that is, we sample \textit{with replacement} from the original forward
degrees $D_1-1,D_2-1
,\ldots,D_n-1$,
where the $m_n$ smallest degrees are discarded. Then, obviously, we obtain
a stochastic ordering $ \underline{X}_i\leq_{st} B_i\leq_{st}
\overline{X}_i$, compare
\cite{hofs3}, Lemma~A.2.8. As a consequence, we can couple
$\{B_i\}_{i=2}^{m_n}$ to $m_n$ i.i.d. random variables $\{\underline
{X}_i\}_{i=1}^{m_n-1},
\{\overline{X}_i\}_{i=1}^{m_n-1}$ such that, a.s.,
%
%e7.16 ###
%
\begin{equation} \label{stochastische-ordening} \underline
{X}_{i-1}\leq B_i\leq\overline{X}_{i-1}.
\end{equation}
The random variables $\{\underline{X}_i\}_{i=1}^{m_n-1}$, as well as
$\{\overline{X}_i\}_{i=1}^{m_n-1}$
are i.i.d., but their distribution depends on $m_n$, since they are
draws \textit{with replacement}
from $D_1-1, \ldots, D_n-1$ where the largest $m_n$, respectively
smallest $m_n$, degrees have been
removed [recall~(\ref{stoch-ord-2})]. Let the \textit{total variation
distance}
between two probability mass functions $p$ and $q$ on ${\mathbb N}$ be
given by
%
%e7.17 ###
%
\begin{equation} d_{\mathrm{TV}}(p,q)=\frac12 \sum_{k=0}^\infty
|p_k-q_k|.
\end{equation}
We shall show that, with $\underline{g}$ and $\overline{g}$,
respectively, denoting the
probability mass functions of $\underline{X}_i$ and $\overline{X}_i$,
respectively,
there exists $\rho'>0$ such that w.h.p.
%
%e7.18 ###
%
\begin{equation} \label{dTV-aim1} d_{\mathrm{TV}}\bigl(\underline{g}^{
(n)},g\bigr)\leq n^{-\rho'},\qquad d_{\mathrm{TV}}\bigl(\overline{g}^{
(n)},g\bigr)\leq n^{-\rho'}.
\end{equation}
This proves the claim for any $\rho<\rho'$, since~(\ref{dTV-aim1})
implies that
$d_{\mathrm{TV}}(\underline{g}^{(n)},\overline{g}^{
(n)})\leq2n^{-\rho'}$, so that
we can couple $\{\underline{X}_i\}_{i=1}^{m_n-1}$ and $\{\overline
{X}_i\}_{i=1}^{m_n-1}$
in such a way that $\prob(\{\underline{X}_i\}_{i=1}^{m_n}=\{\overline
{X}_i\}_{i=1}^{m_n})\leq2m_n n^{-\rho'}
=o(1),$ when $m_n=n^\rho$ with $\rho'<\rho$. In particular, this
yields that we can couple
$\{B_i\}_{i=2}^{m_n}$ to $\{\underline{X}_i\}_{i=1}^{m_n-1}$ in such a
way that
$\{B_i\}_{i=2}^{m_n}=\{\underline{X}_i\}_{i=1}^{m_n-1}$ w.h.p.
Then, again from
(\ref{dTV-aim1}), we can couple $\{\underline{X}_i\}_{i=1}^{m_n-1}$
to a sequence
of i.i.d. random variables $\{\Bindep_i\}_{i=1}^{m_n-1}$ such that
$\{\underline{X}_i\}_{i=1}^{m_n-1}
=\{\Bindep_i\}_{i=1}^{m_n-1}$ w.h.p. Thus,~(\ref{dTV-aim1})
completes the proof of Proposition~\ref{prop-indep}.

To prove~(\ref{dTV-aim1}), we bound
%
%e7.19 ###
%
\begin{equation} d_{\mathrm{TV}}\bigl(\underline{g}^{(n)},g\bigr)
\leq
d_{\mathrm{TV}}\bigl(\underline{g}^{(n)},g^{(n)}\bigr) +d_{\mathrm
{TV}}\bigl(g^{(n)},g\bigr),
\end{equation}
and a similar identity holds for $d_{\mathrm{TV}}(\overline{g}^{
(n)},g)$, where
%
%e7.20 ###
%
\begin{equation} \label{gn-def} g^{(n)}_k=\frac{1}{L_n}\sum
_{j=1}^n (k+1)\indic{D_j=k+1}.
\end{equation}
In \cite{hofs3}, (A.1.11), it is shown that there exists $\alpha_2,
\beta_2>0$ such that
%
%e7.21 ###
%
\begin{equation} \prob\bigl(d_{\mathrm{TV}}\bigl(g^{(n)},g\bigr
)\geq n^{-\alpha
_2}\bigr)\leq n^{-\beta_2}.
\end{equation}
Thus, we are left to investigate $d_{\mathrm{TV}}(\underline{g}^{
(n)},g^{(n)})$ and
$d_{\mathrm{TV}}(\overline{g}^{(n)},g^{(n)})$. We bound
%
%e7.24 ###
%e7.23 ###
%e7.22 ###
%
\begin{eqnarray} d_{\mathrm{TV}}\bigl(\underline
{g}^{(n)},g^{(n)}\bigr)
&=& \frac12 \sum_{k=0}^{\infty}\big|\underline{g}^{(n)}_k-g^{
(n)}_k\big|\nonumber\\
&\leq&\sum_{k=0}^{\infty}(k+1) \biggl(\frac
{1}{{\underline L}_n}- \frac{1}{L_n} \biggr)\sum_{j=1}^{n-m_n}
\indic
{D_j=k+1}\nonumber\\
&&{} + \sum_{k=0}^{\infty}(k+1)\frac{1}{{\underline L}_n}\sum
_{j=n-m_n+1}^{n} \indic{D_{(j)}=k+1}\\
&\leq& \biggl(\frac
{L_n-{\underline L}_n}{L_n{\underline L}_n} \biggr)\sum
_{j=1}^{n-m_n}D_{(j)} +\frac{1}{{\underline L}_n}\sum
_{j=n-m_n+1}^{n} D_{(j)}\nonumber\\
 &\leq&2 \biggl(\frac{L_n-{\underline
L}_n}{{\underline L}_n} \biggr)
=\frac{2}{{\underline L}_n}\sum
_{j=n-m_n+1}^{n} D_{(j)},\nonumber
\end{eqnarray}
where ${\underline L}_n=\sum_{j=1}^{n-m_n}D_{(j)}$.
Define $b_n=\Theta(n/m_n)^{1/(\tau-1)}$. Then, from $1-F(x)=x^{-(\tau
-1)}L(x),$
and concentration results for the binomial distribution, we have, w.h.p.,
$D_{(n-m_n+1)}\geq b_n$, so that, w.h.p.,
%
%e7.25 ###
%
\begin{equation} \label{bound-order-stats} \frac{L_n-{\underline
L}_n}{{\underline L}_n} =\frac{1}{{\underline L}_n} \sum
_{j=n-m_n+1}^{n} D_{(j)} \leq\frac{1}{{\underline L}_n} \sum
_{j=1}^{n} D_{j} \indic{D_j\geq b_n}.
\end{equation}
Now, in turn, by the Markov inequality,
%
%e7.26 ###
%
\begin{eqnarray} \label{MI-order-stats} \prob\Biggl(\frac
{1}{{\underline L}_n} \sum_{j=1}^{n} D_{j} \indic{D_j\geq b_n}\geq
n^{\vep}b_n^{2-\tau} \Biggr) &\leq& n^{-\vep}b_n^{\tau-2}\expec
\Biggl[\frac{1}{{\underline L}_n} \sum_{j=1}^{n} D_{j} \indic
{D_j\geq
b_n} \Biggr]\nonumber\\[-8pt]\\[-8pt]
&\leq& Cn^{-\vep},\nonumber
\end{eqnarray}
so that
%
%e7.27 ###
%
\begin{equation} \prob\bigl(d_{\mathrm{TV}}\bigl(\underline{g}^{
(n)},g^{(n)}\bigr)\geq n^{\vep}b_n^{-(\tau-2)} \bigr) =o(1).
\end{equation}
Thus, w.h.p., $d_{\mathrm{TV}}(\underline{g}^{(n)},g^{
(n)})\leq n^{\vep}(m_n/n)^{(\tau-2)/(\tau-1)}$,
which proves~(\ref{dTV-aim1}) when we take $m_n=n^\rho$ and
$\rho'=(1-\rho)(\tau-2)/(\tau-1)-\vep>0$. The upper bound for
$d_{\mathrm{TV}}(\overline{g}^{(n)},g^{(n)})$ can
be treated similarly.

%We let ${\underline G}^{\sss(n)}$ denote the
%distribution function of ${\underline X_1}$. The distribution function
%${\underline G}^{\sss(n)}$ clearly depends on $\olmn$, but this
%dependence is omitted from
%the notation. We finally denote the distribution function of $B$
%having the
%size-biased distribution in \refeq{eqn:size-bias} by $G(x)=\prob(B\leq
%x)=\sum_{j\leq x} g_j$.

%s7.3 ###
\subsection{Auxiliary lemmas for $2<\tau<3$}
In this section we treat some lemmas that complete the proof of
Proposition~\ref{lemma:CLT-sum2}(a) for $\tau\in(2,3)$. In particular,
we shall verify condition (ii) in Remark~\ref{rem-CLTexchang}.

\begin{Lemma}[(A bound on the expected value of $1/S_i$)]
\label{boundexpS1}
Fix $\tau\in(2,3)$. For $\ulmn, \olmn$ such that
$\log{(a_n/\ulmn)}$, $\log{(\olmn/a_n)}=o(\sqrt{\log n})$
and for $b_n$ such that $b_n\to\infty$,
%
%e7.28 ###
%
\begin{eqnarray} \mbox{(\textup{i})}&&\hspace*{8pt} \sum_{i=1}^{\ulmn}
\expec
[1/\Sigma_i]=O(1),\nonumber\\
\mbox{(\textup{ii})}&&\hspace*{8pt}\sum_{i=b_n}^{\ulmn} \expec[1/\Sigma
_i]=o(1)\quad \mbox{and}\\
\mbox{(\textup{iii})}&&\hspace*{8pt} \sup_{i\leq\olmn}\expec
\bigl[(R_i/S_{R_i})\indicwo{\{\ulmn+1\le R_i\le\olmn\}}\bigr
]<\infty.\nonumber
\end{eqnarray}
\end{Lemma}

\begin{pf} Let $\ulmn=o(a_n)$. Let
%
%e7.29 ###
%
\begin{equation} M_i=\max_{2\le j\le i} (B_j-1).
\end{equation}
Then, we use that, for $1\leq i\le\ulmn$,
%
%e7.30 ###
%
\begin{equation} \quad\Sigma_i\equiv1+\sum_{j=2}^i (B_j-1)\geq\max
_{2\le j\le i} (B_j-1)-(i-2) =M_i-(i-2).
\end{equation}
Fix $\delta>0$ small, and split
%
%e7.31 ###
%
\begin{eqnarray} \expec[1/\Sigma_i ]&\leq&\expec
\bigl[1/\Sigma_i\indic{\Sigma_i\leq i^{1+\delta}} \bigr] +\expec
\bigl[1/\Sigma_i\indic{\Sigma_i>i^{1+\delta}} \bigr]\nonumber\\
[-8pt]\\[-8pt]
&\leq&\prob(\Sigma
_i\leq i^{1+\delta}) + i^{-(1+\delta)}.\nonumber
\end{eqnarray}
Now, if $\Sigma_i\leq i^{1+\delta}$, then $M_i\leq i^{1+\delta
}+i\leq2i^{1+\delta}$, and
$\Sigma_j\leq i^{1+\delta}+i\leq2i^{1+\delta}$ for all $j\leq i$.
As a result,
for each $j\leq i$, the conditional probability that
$B_j-1>2i^{1+\delta}$,
given $\Sigma_{j-1}\leq2i^{1+\delta}$ and $\{D_s\}_{s=1}^n$ is at least
%
%e7.32 ###
%
\begin{eqnarray} \frac{1}{L_n}\sum_{s=1}^n D_s\indic
{D_s>2i^{1+\delta}} &\geq&2i^{1+\delta}\sum_{s=1}^n \indic
{D_s>2i^{1+\delta}}/L_n\nonumber\\[-8pt]\\[-8pt]
&=&2i^{1+\delta} \operatorname{BIN}
\bigl(n,1-F(2i^{1+\delta}) \bigr)/L_n.\nonumber
\end{eqnarray}
Further, by~(\ref{Fcond(2,3)}), for some $c>0$, $n[1-F(2i^{1+\delta
})]\geq2cn i^{-(1+\delta)(\tau-1)}$,
so that, for $i\leq\ulmn=o(n^{(\tau-2)/(\tau-1)})$, $n
i^{-(1+\delta)(\tau-1)}\geq n^{\vep}$
for some $\vep>0$. We shall use Azuma's inequality that states that
for a binomial random variable
$\operatorname{BIN}(N,p)$ with parameters $N$ and $p$, and all $t>0$,
%
%e7.33 ###
%
\begin{equation} \label{binbd}
\prob\bigl(\operatorname{BIN}(N,p)\leq Np-t\bigr)\leq
\exp\biggl\{-\frac{2t^2}{N} \biggr\}.
\end{equation}
As a result,
%
%e7.34 ###
%
\begin{eqnarray}&&\prob\Bigl(\operatorname{BIN}\bigl
(n,1-F(2i^{1+\delta})
\bigr)\leq\expec\bigl[\operatorname{BIN} \bigl(n,1-F(2i^{1+\delta
}) \bigr)\bigr]/2 \Bigr)\nonumber\\[-8pt]\\[-8pt]
 &&\qquad\leq
{\mathrm e}^{-n[1-F(2i^{1+\delta})]/2}
\leq{\mathrm e}^{-n^{\vep}},\nonumber
\end{eqnarray}
so that, with probability at least $1-{\mathrm e}^{-n^\vep}$,
%
%e7.35 ###
%
\begin{equation} \label{sum-Di-bd} \frac{1}{L_n}\sum_{s=1}^n
D_s\indic{D_s>2i^{1+\delta}}\geq ci^{-(1+\delta)(\tau-2)}.
\end{equation}
Thus, the probability that in the first $i$ trials, no vertex with
degree at least $2i^{1+\delta}$ is chosen is bounded above by
%
%e7.36 ###
%
\begin{equation} \bigl(1-ci^{-(1+\delta)(\tau-2)} \bigr)^i+{\mathrm
e}^{-n^\vep} \leq{\mathrm e}^{-ci^{1-(1+\delta)(\tau-2)}}+{\mathrm
e}^{-n^\vep},
\end{equation}
where we used the inequality $1-x\leq{\mathrm e}^{-x}, x\ge0$.
Finally, take $\delta>0$ so small that $1-(1+\delta)(\tau-2)>0$;
then we arrive at
%
%e7.37 ###
%
\begin{equation} \label{Sigma-i-bd} \expec[1/\Sigma_i ]\leq
i^{-(1+\delta)}+{\mathrm e}^{-ci^{1-(1+\delta)(\tau-2)}}+{\mathrm
e}^{-n^\vep},
\end{equation}
which, when summed over $i\leq\ulmn$, is $O(1)$. This proves (i). For
(ii), we note that,
for any $b_n\rightarrow\infty$, the sum of the r.h.s. of~(\ref
{Sigma-i-bd}) is $o(1)$.
This proves (ii).

To prove (iii), we take $\log{(a_n/\ulmn)}, \log{(\olmn
/a_n)}=o(\sqrt{\log n})$. We bound the expected value by
\[
\olmn\expec\bigl[(1/S_{R_i})\indicwo{\{\ulmn+1\le R_i\le\olmn\}
}\bigr].
\]
For $\ulmn+1\leq i \leq\olmn$,
%
%e7.38 ###
%
\begin{equation} S_i=D_1+\sum_{j=2}^i (B_j-1)\geq1+\sum_{j=2}^i
(B_j-1)=\Sigma_i,
\end{equation}
and the above derived bound for the expectation $\expec[1/\Sigma_i]$
remains valid
for $\ulmn+1\leq i \leq\olmn$, since also for $i\leq\olmn$, we
have $ni^{-(1+\delta)(\tau+1)}\geq n^{\vep}$; moreover since
the r.h.s. of~(\ref{Sigma-i-bd}) is decreasing in $i$, we obtain
%
%e7.39 ###
%
\begin{equation} \label{Sigma-overal-bnd} \expec[1/\Sigma_i
] \leq\ulmn^{-(1+\delta)}+{\mathrm e}^{-c\ulmn^{1-(1+\delta
)(\tau-2)}}+{\mathrm e}^{-n^\vep}.
\end{equation}
Consequently,
%
%e7.40 ###
%
\begin{eqnarray} &&\olmn\expec\bigl[(1/S_{R_i})\indicwo{\{\ulmn
+1\le
R_i\le\olmn\}}\bigr]\nonumber\\[-8pt]\\[-8pt]
&&\qquad\le\olmn \bigl( \ulmn^{-(1+\delta)}+{\mathrm
e}^{-c\ulmn^{1-(1+\delta)(\tau-2)}}+{\mathrm e}^{-n^\vep}
\bigr)=o(1),\nonumber
\end{eqnarray}
using that $\log{(a_n/\ulmn)}, \log{(\olmn/a_n)}=o(\sqrt{\log n})$.
This proves (iii).
\end{pf}

\begin{Lemma}[(Bounds on $S_{\olmn}$)]
\label{lem-bdsS-(2,3)}
Fix $\tau\in(2,3)$. Then, w.h.p., for every $\vep
_n\rightarrow0$,
%
%e7.41 ###
%
\begin{equation} \label{aim-lb-Sm(2,3)-lb} S_{a_n}\geq\vep_n
n^{1/(\tau-1)},
\end{equation}
while, w.h.p., uniformly for all $m\leq\olmn$,
%
%e7.42 ###
%
\begin{equation} \label{aim-lb-Sm(2,3)-ub} \expec_n[S_m]\leq\vep
_n^{-1} m n^{(3-\tau)/(\tau-1)}.
\end{equation}
\end{Lemma}

\begin{pf}
We prove~(\ref{aim-lb-Sm(2,3)-lb}) by noting that, by
(\ref{sum-Di-bd}) and
the fact that $\vep_n\downarrow0$,
%
%e7.43 ###
%
\begin{equation}\quad\sum_{i=1}^n D_{i}\indic{D_i\geq\vep_n
n^{1/(\tau
-1)}} \geq cn \bigl(\vep_n n^{1/(\tau-1)} \bigr)^{-(\tau-2)} =
(cn/a_n ) \vep_n^{-(\tau-2)}.
\end{equation}
Therefore, the probability to choose none of these vertices with degree
at least
$\vep_n n^{1/(\tau-1)}$ before time $a_n$ is bounded by
%
%e7.44 ###
%
\begin{equation} \bigl(1-c \vep_n^{-(\tau-2)} n^{-(2-\tau)/(\tau
-1)}\bigr)^{a_n}\leq{\mathrm e}^{-c \vep_n^{-(\tau-2)}}=o(1)
\end{equation}
for any $\vep_n \downarrow0$. In turn, this implies that, w.h.p.,
$S_{a_n}\geq\vep_n n^{1/(\tau-1)}-a_n \geq\vep_n n^{1/(\tau-1)}/2$,
whenever $\vep_n$ is such that $\vep_n n^{1/(\tau-1)}\geq2a_n$.

To prove~(\ref{aim-lb-Sm(2,3)-ub}), we use that, w.h.p.,
$D_{(n)}\leq\vep_n^{-1} n^{1/(\tau-1)}$ for any $\vep
_n\rightarrow0$.
Thus, w.h.p., using the inequality $L_n>n,$
%
%e7.45 ###
%
\begin{equation} \label{expec-n-Sm} \expec_n[S_m]\leq m \expec
_n[B_{2}] \leq\frac{m}{n}\sum_{j=1}^{n} D_j(D_j-1) \indic{D_j\leq
\vep_n^{-1} n^{1/(\tau-1)}}.
\end{equation}
Thus, in order to prove the claimed uniform bound, it suffices to give
a bound on
the above sum that holds w.h.p. For this, the expected value of the
sum on the r.h.s. of~(\ref{expec-n-Sm}) equals
%
%e7.47 ###
%e7.46 ###
%
\begin{eqnarray}&& \expec\Biggl[\sum_{j=1}^{n} D_j(D_j-1) \indic
{D_j\leq
\vep_n^{-1} n^{1/(\tau-1)}} \Biggr]\nonumber\\
&&\qquad\leq n\sum_{j=1}^{\vep_n^{-1}
n^{1/(\tau-1)}} j\prob(D_1>j)\\
&&\qquad\leq c_2n\sum_{j=1}^{\vep
_n^{-1} n^{1/(\tau-1)}} j^{2-\tau} \leq\frac{c_2}{3-\tau}n \vep
_n^{-(3-\tau)} n^{(3-\tau)/(\tau-1)}.\nonumber
\end{eqnarray}
Since $\tau\in(2,3)$, $\vep_n^{\tau-2}\rightarrow\infty$, so that
uniformly for all $m\leq\olmn$, by the Markov inequality,
%
%e7.49 ###
%e7.48 ###
%
\begin{eqnarray} &&\prob\bigl(\expec_n[S_m]\geq\vep_n^{-1} m
n^{(3-\tau
)/(\tau-1)} \bigr)\nonumber\\
&&\qquad\leq
\vep_n m^{-1}n^{-(3-\tau)/(\tau-1)}\expec
\bigl[\expec_n[S_m]\indic{\max_{j=1}^n D_j\leq\vep_n^{-1}
n^{1/(\tau-1)}} \bigr]\\
&&\qquad\leq c_2\vep_n^{-(2-\tau)}=o(1).\nonumber
\end{eqnarray}
This completes the proof of~(\ref{aim-lb-Sm(2,3)-ub}).
\end{pf}

%%%%%%%%%%%%%%%%%%%%%%%%%%%%%%%%%%%%%%%%%%%%%%%%%%%%%%%%

%s7.4 ###
\subsection{Auxiliary lemmas for $\tau>3$}
\label{lem-cond-LD-est}
In the lemmas below we use the coupling~\eqref
{stochastische-ordening}. We define the partial sums
${\underline S}_i$ and ${\overline S}_i$ by
%
%e7.50 ###
%
\begin{equation} \label{def-uSenos} {\underline S}_i=\sum_{j=1}^{i-1}
({\underline X}_i-1 ), \qquad {\overline S}_i=\sum
_{j=1}^{i-1} ({\overline X}_i-1 ), \qquad i\geq2.
\end{equation}
As a consequence of \eqref{stochastische-ordening}, we obtain for
$i\ge2$,
%
%e7.51 ###
%
\begin{equation} {\underline S}_i\le\sum_{j=2}^i (B_j-1 )\le
{\overline S}_i\qquad \mbox{a.s.}
\end{equation}

\begin{Lemma}[(A conditional large deviation estimate)]
Fix $\tau>2$. Then w.h.p., there exist a $c>0$ and $\eta>0$
sufficiently small,
such that for all $i\geq0$, and w.h.p.,
%
%e7.52 ###
%
\begin{equation} \label{cond-LD-estim} \prob_n(\underline{S}_i\leq
\eta i)\leq{\mathrm e}^{-c i}.
\end{equation}
The same bound applies to $\overline{S}_i$.
\end{Lemma}

\begin{pf}
We shall prove~(\ref{cond-LD-estim}) using a conditional large
deviation estimate, and
an analysis of the moment generating function of ${\underline X}_1$, by
adapting the proof of
the upper bound in Cram\'er's theorem. Indeed, we rewrite and bound,
for any $t\geq0$,
%
%e7.53 ###
%
\begin{equation} \prob_n(\underline{S}_i\leq\eta i) =\prob
_n({\mathrm e}^{-t\underline{S}_i}\geq{\mathrm e}^{-t\eta i})\leq
({\mathrm e}^{t\eta}\phi_n(t) )^i,
\end{equation}
where $\phi_n(t)=\expec_n[{\mathrm e}^{-t({\underline X}_1-1)}]$ is
the (conditional) moment generating function of
${\underline X}_1-1$. Since ${\underline X}_1-1\geq0$, we have that
${\mathrm e}^{-t({\underline X}_1-1)}\leq1$,
and ${\underline X}_1\convd B$, where $B$ has the size-biased
distribution in~(\ref{eqn:size-bias}).
Therefore, for every $t\geq0$, $\phi_n(t)\convd\phi(t)$, where
$\phi(t)=\expec[{\mathrm e}^{-t(B-1)}]$ is
the Laplace transform of
$B$. Since this limit is a.s. constant, we even obtain that $\phi
_n(t)\convp\phi(t)$. Now, since $\expec[B]=\nu>1$,
for each $0<\eta<\expec[B]-1$, there exists a $t^*>0$ and $\vep>0$
such that ${\mathrm e}^{-t^*\eta}\phi(t^*)\leq1-2\vep$.
Then, since ${\mathrm e}^{t^*\eta}\phi_n(t^*)\convp{\mathrm
e}^{t^*\eta}\phi(t^*)$, w.h.p. and for all $n$\vspace*{1pt} sufficiently large,
$|{\mathrm e}^{t^*\eta}\phi_n(t^*)-{\mathrm e}^{t^*\eta}\phi
(t^*)|\leq\vep$, so that ${\mathrm e}^{-t^*\eta}\phi_n(t^*)\leq
1-\vep<1.$
The proof for $\overline{S}_i$ follows since $\overline{S}_i$ is
stochastically larger than
$\underline{S}_i$. This completes the proof.
\end{pf}

%Fix $\tau>3$. Let $\olmn$ be such that $\log{(\olmn/a_n)}=o(\sqrt{
%Then, w.h.p.,
%
%a.s.

\begin{Lemma}
\label{lem-reciproke}
Fix $\tau>3$. For $\ulmn, \olmn$ such that $\log{(\olmn/a_n)}$,
$\log{(a_n/\ulmn)}=o(\sqrt{\log n})$,
%
%e7.54 ###
%
\begin{eqnarray}
\sup_{i\leq
\olmn}\expec\bigl[R_i/S_{R_i}\indicwo{\{\ulmn+1\le R_i\le
\olmn\}}\bigr]<\infty.
\end{eqnarray}
\end{Lemma}

\begin{pf} Take $\ulmn+1\leq k \leq\olmn$ and recall the definition
of $\Sigma_k<S_k$ in~(\ref{cond-CLT-(2,3)}). For $\eta>0$,
\begin{eqnarray*}
\expec[k/\Sigma_k]&=&
\expec[k/\Sigma_k]\indicwo{\{\Sigma_k<\eta k\}}
+\expec[k/\Sigma_k]\indicwo{\{\Sigma\geq\eta k\}}\\
&\leq&
\expec[k/\Sigma_k]\indicwo{\{\Sigma_k<\eta k\}}+\eta^{-1}\\
&\le&
k\prob(\Sigma_k<\eta k)+\eta^{-1}\le k\prob({\underline S}_k<\eta
k)+\eta^{-1},
\end{eqnarray*}
since $\Sigma_k=1+\sum_{j=2}^k (B_j-1)>{\underline S}_k$, a.s.
Applying the large deviation estimate from the previous lemma, we obtain
\[
\expec[k/\Sigma_k]\leq\eta^{-1}+k{\mathrm e}^{-c_2 k}
\]
for each $\ulmn+1\le k \leq\olmn$. Hence,
%
%e7.55 ###
%
\begin{equation} \sup_{i\leq\olmn}\expec\bigl[R_i/S_{R_i}\indicwo
{\{
\ulmn+1\le R_i\le\olmn\}}\bigr] \leq\eta^{-1}+\olmn{\mathrm e}^{-c_2
\ulmn}.
\end{equation}
\upqed
\end{pf}

\begin{Lemma}
\label{lem-bd-squareBSratio}
Fix $\tau>3$, and let $\ulmn$ be such that $\log{(a_n/\ulmn
)}=o(\sqrt{\log{n}})$.
Then, for each sequence $C_n\to\infty$,
%
%e7.56 ###
%
\begin{equation}
\label{sec-mom-sum-exch2}
\prob_n \Biggl(\sum_{j=2}^{\ulmn} B_j^2/S_j^2>C_n \Biggr)\convp0.
\end{equation}
Consequently,
%
%e7.57 ###
%
\begin{equation} \sum_{j=2}^{\ulmn} B_j^2/S_j^2=\Op(1).
\end{equation}
\end{Lemma}

\begin{pf} If we show that the conditional expectation of
$\sum_{j=2}^{\ulmn} B_j^2/S_j^2$, given $\{D_i\}_{i=1}^n$,
is finite, then \eqref{sec-mom-sum-exch2} holds. Take $a\in(1,\min
(2,\tau-2))$; this is possible since
$\tau>3$. We bound
%
%e7.58 ###
%
\begin{eqnarray} \label{amoment}
\expec_n \biggl[ \biggl(\frac
{B_j}{S_j} \biggr)^2 \biggr] &\leq& 2 \biggl(\expec_n \biggl[
\biggl(\frac
{B_j-1}{S_j} \biggr)^2 \biggr] \biggr) +2\expec_n \biggl[\frac1{(S_j)^2}
\biggr]\nonumber\\[-8pt]\\[-8pt]
&\leq& 2 \biggl(\expec_n \biggl[ \biggl(\frac{B_j-1}{S_j}
\biggr)^a \biggr] \biggr) +2\expec_n \biggl[\frac1{(S_j)^a}
\biggr] .\nonumber
\end{eqnarray}
By stochastic domination and Lemma~\ref{lem-cond-LD-est}, we find
that, w.h.p., using
$a>1$,
\[
\sum_{j=2}^{\ulmn}
\expec_n \biggl[
\frac1{(S_j)^{a}}
\biggr]<\infty.
\]

We will now bound \eqref{sec-mom-sum-exch2}.
Although, by definition
\[
S_j=D_1+\sum_{i=2}^j (B_i-1)
\]
for the asymptotic statements that we discuss here, we may as well
replace this definition
by
%
%e7.59 ###
%
\begin{equation}
\label{alternatieve-def}
S_j=\sum_{i=2}^j (B_i-1),
\end{equation}
and use exchangeability, so that
\[
\expec_n \biggl[ \biggl(\frac{B_j-1}{S_j} \biggr)^a \biggr]=
\expec_n \biggl[ \biggl(\frac{B_2-1}{S_j} \biggr)^a \biggr],
\]
since for each $j$, we have $\frac{B_j-1}{S_j}\stackrel{d}{=}\frac
{B_1}{S_j}$.
Furthermore, for $j\geq2$,
\[
\expec_n \biggl[ \biggl(\frac{B_2-1}{S_j} \biggr)^a \biggr]
\leq
\expec_n \biggl[ \biggl(\frac{B_2-1}{S_{3,j}} \biggr)^a \biggr],
\]
where $S_{3,j}=(B_3-1)+\cdots+(B_j-1)$. Furthermore,
we can replace $S_{3,j}$ by ${\underline S}_{3,j}=({\underline
X_3}-1)+\cdots+
({\underline X}_j-1)$,
which are mutually independent and sampled from $D_{(1)}-1,\ldots
,D_{(\ulmn)}-1$, as above and which are also independent of
$B_2$. Consequently,
%
%e7.60 ###
%
\begin{eqnarray}
\label{proofamoment}
\sum_{j=2}^{\ulmn}\expec_n \biggl[ \biggl(\frac{B_j-1}{S_j}
\biggr)^2 \biggr] &\leq&
\sum_{j=2}^{\ulmn}\expec_n \biggl[ \biggl(\frac{B_j-1}{S_j}
\biggr)^a
\biggr]\nonumber\\
&=&
\sum_{j=2}^{\ulmn}\expec_n \biggl[ \biggl(\frac{B_2-1}{S_j}
\biggr)^a
\biggr]\nonumber\\
& \leq&\expec_n \biggl[ \biggl(\frac{B_2-1}{S_2} \biggr)^a \biggr]+
\sum_{j=3}^{\ulmn}\expec_n \biggl[ \biggl(\frac{B_2-1}{S_{3,j}}
\biggr)^a \biggr]\\
&\leq&1+
\sum_{j=3}^{\ulmn}\expec_n \biggl[ \biggl(\frac
{B_2-1}{{\underline
S}_{3,j}} \biggr)^a \biggr]\nonumber\\
& =&1+
\expec_n[(B_2-1)^a]\sum_{j=3}^{\ulmn}\expec_n \biggl[ \biggl
(\frac
{1}{{\underline S}_{3,j}} \biggr)^a \biggr].\nonumber
\end{eqnarray}
Finally, the expression $\sum_{j=3}^{\ulmn} \expec_n
[1/{\underline S}^a_{2,j} ]$
can be shown to be finite as above.
\end{pf}

%Fix $\tau>3$, then whp,
%
%X}_1]|$ is
%bounded by
%$$
%+
% (\frac1{{\underline L}_n}-\frac1{L_n} )
%$$
%Since w.h.p. $D_{\sss(n)}\le n^{\chi}$, where $\chi=1/(\tau-1)+
%term is bounded by $n^{-1} m_n n^{2\chi}$ and hence for $\tau>5$. {\bf
%Gerard:
%This is too crude, what is wrong?}. The second term we rewrite as
%$$
%$$
%Obviously $\frac{{\underline L}_n}{{\underline L}_n}\to1$, w.h.p..
%

\begin{Lemma}[(Logarithmic asymptotics of $S_{\olmn}$)]
\label{lemma-lgS}
Fix $\tau>3$, and let $\ulmn$ be such that $\log{(a_n/\ulmn
)}=o(\sqrt{\log{n}})$. Then,
%
%e7.61 ###
%
\begin{equation}
\label{SL-excha2}
\log{S_{\ulmn}}-\log{\ulmn}=
o_{\prob}\bigl(\sqrt{\log{\ulmn}}\bigr).
\end{equation}
\end{Lemma}

\begin{pf} As in the previous lemma we define w.l.o.g. $S_j$ by \eqref
{alternatieve-def}.
Then,
\[
S_j\le_{st}{\overline S}_j,\vadjust{\goodbreak}
\]
where ${\overline S}_j$ is a sum of i.i.d. random variables
${\overline X}_i-1$, where the
${\overline X}_i$ are sampled from $D_1,\ldots,D_n$ with replacement,
where $\ulmn$ of the
vertices with the smallest degree(s) have been removed. Using the
Markov inequality,
%
%e7.62 ###
%
\begin{eqnarray}
\label{upbound-for-log}
\prob_n \bigl(\log(S_{\ulmn}/\ulmn)>c_n \bigr)&=&
\prob_n (S_{\ulmn}/\ulmn>{\mathrm e}^{c_n} )\nonumber\\[-8pt]\\[-8pt]
&\le&{\mathrm e}^{-c_n}\expec_n[{\overline S}_{\ulmn}/\ulmn
]={\mathrm
e}^{-c_n}\expec_n[{\overline X}_{i}-1].\nonumber
\end{eqnarray}
%
%By definition,
% $$
% \expec_n[{\overline X}_{i}]=\frac1{{\overline L}_n}
% \sum_{j=\olmn+1}^n D_{\sss(j)}(D_{\sss(j)}-1)
% \le
% \frac{L_n}{{\overline L}_n}
% \frac1{L_n}
% \sum_{j=1}^n D_{\sss(j)}(D_{\sss(j)}-1),
% $$
%where
% $$
% {\overline L}_n=D_{\sss(\olmn+1)}+\cdots+D_{\sss(n)}=L_n- (D_{
% $$
%since
% $$
% \frac{D_{\sss(1)}+\cdots+D_{\sss(\olmn)}}{n}\leq\frac{D_{1}+
% =\Op(\olmn/n).
% $$
%The upper bound is completed by noting that, for $\tau>3$,
We shall prove below that, for $\tau>3$,
$\expec_n[{\overline X}_{1}]\convp\nu$ so that
%
%e7.63 ###
%
\begin{equation} \label{bd-mean-Sm-tau>3} \expec_n[S_m]\leq\nu m
\bigl(1+o_{\prob}(1)\bigr).
\end{equation}
Indeed,
from \cite{hofs3}, Proposition A.1.1, we know that there are $\alpha
,\beta>0$,
such that
%
%e7.64 ###
%
\begin{equation} \label{nun-n-conv} \prob(|\nu_n-\nu|>n^{-\alpha
})\leq n^{-\beta},
\end{equation}
where
%
%e7.65 ###
%
\begin{equation}
\label{def-nu}\qquad
\nu_n= \sum_{j=1}^\infty j g_j^{(n)}=\sum_{j=1}^\infty j(j+1)
\frac1{L_n}
\sum_{i=1}^n
\indic{D_i=j+1}
=\frac1{L_n}\sum_{i=1}^n D_i(D_i-1).
\end{equation}
Define $\overline{\nu}_n=\expec_n[{\overline X}_1]$. Then we claim that
there exists $\alpha,\beta>0$ such that
%
%e7.66 ###
%
\begin{equation} \label{olnun-n-conv} \prob(|\overline{\nu}_n-\nu
_n|>n^{-\alpha})\leq n^{-\beta}.
\end{equation}
To see~(\ref{olnun-n-conv}), by definition of
$\overline{\nu}_n=\expec_n[{\overline X}_1]$,
%
%e7.67 ###
%
\begin{eqnarray}
|\overline{\nu}_n-\nu_n|&=&
\Bigg|\frac1{{\overline L}_n}\sum_{i=\ulmn+1}^n D_{(i)}\bigl(D_{
(i)}-1\bigr) -\frac1{L_n}
\sum_{i=1}^n D_i(D_i-1)\Bigg|
\nonumber\\
&\leq& \Bigg|\frac1{{\overline L}_n}\sum_{i=\ulmn+1}^n D_{
(i)}\bigl(D_{(i)}-1\bigr)
-\frac1{L_n}\sum_{i=\ulmn+1}^n D_{(i)}\bigl(D_{(i)}-1\bigr)
\Bigg|\\
&&{} + \Bigg|\frac1{L_n}\sum_{i=\ulmn+1}^n D_{(i)}\bigl
(D_{(i)}-1\bigr)
-\frac1{L_n}\sum_{i=1}^n D_{(i)}\bigl(D_{(i)}-1\bigr)\Bigg
|.\nonumber
\label{ineq-degrees}
\end{eqnarray}
The first term on the r.h.s. of \eqref{ineq-degrees} is with
probability at least
$1-n^{-\beta}$ bounded above by $n^{-\alpha}$, w.h.p., since it is
bounded by
\[
\biggl(\frac{L_n-{\overline L}_n}{{\overline L}_n} \biggr)\frac1{L_n}
\sum_{i=1}^n D_i(D_i-1),\vadjust{\goodbreak}
\]
and since, using~(\ref{bound-order-stats}) and~(\ref
{MI-order-stats}), $(L_n-{\overline L}_n)/{\overline L}_n=o_{\prob
}(n^{-\alpha})$
for some $\alpha>0$. The second term on the r.h.s. of \eqref
{ineq-degrees} is bounded by
%
%e7.68 ###
%
\begin{equation} \label{small-order-stats} \frac1{L_n} \sum
_{j=1}^{\ulmn} D_{(j)}^2\leq\frac1{L_n} \sum_{j=1}^{\ulmn}
D_j^2=o_{\prob}(n^{-\alpha}),
\end{equation}
since $\tau>3$. This completes the proof of~(\ref{olnun-n-conv}).
Combining \eqref{upbound-for-log} with $c_n=o(\sqrt{\log{\ulmn}})$
and the fact that $\expec_n[{\overline X}_{1}]\convp\nu$,
we obtain an upper bound for the left-hand side of \eqref{SL-excha2}.

For the lower bound, we simply make use of the fact that, by Lemma~\ref
{lem-cond-LD-est}
and w.h.p., $S_{\ulmn}\geq\eta\ulmn$, so that
$\log{S_{\ulmn}}-\log{\ulmn}\geq\log{\eta}=
o_{\prob}(\sqrt{\log{\ulmn}})$.
\end{pf}

%%%%%%%%%%%%%%%%%%%%%%%%%%%%%%%%%%%%%%%%%%%%%%%%%%%%%%

\begin{Lemma}Fix $\tau>3$, and let $\ulmn$ be such that $\log
{(a_n/\ulmn)}=o(\sqrt{\log{n}})$. Then,
\label{Mar-Zyg}
%
%e7.69 ###
%
\begin{equation}
\label{lln(b)-excha2}
\sum_{j=1}^{\ulmn} \frac{S_j-(\nu-1) j}{S_j(\nu-1)j}
= \sum_{j=1}^{\ulmn} \biggl[\frac{1}{(\nu-1)j}-\frac{1}{S_j}
\biggr]
=o_{\prob}\bigl(\sqrt{\log{\ulmn}}\bigr).
\end{equation}
\end{Lemma}

\begin{pf}
We can stochastically bound the sum \eqref{lln(b)-excha2} by
%
%e7.70 ###
%
\begin{equation} \qquad\sum_{j=1}^{\ulmn} \biggl[\frac{1}{(\nu
-1)j}-\frac
{1}{{\underline S}_j} \biggr]\leq\sum_{j=1}^{\ulmn} \biggl[\frac
{1}{(\nu-1)j}-\frac{1}{S_j} \biggr]\leq\sum_{j=1}^{\ulmn}
\biggl[\frac{1}{(\nu-1)j}-\frac{1}{{\overline S}_j} \biggr].
\end{equation}
We now proceed by proving \eqref{lln(b)-excha2} both with $S_j$
replaced by
${\overline S}_j$,\vspace*{1.5pt} and with $S_j$ replaced by
${\underline S}_j$.
In the proof of Lemma~\ref{lemma-lgS} we have shown that $\expec
_n[{\overline X}_1]$\vspace*{1pt}
converges, w.h.p., to $\nu$.
Consequently, we can copy the proof of Proposition~\ref{lemma:CLT-sum1}(a)
to show that, w.h.p.,
%
%e7.71 ###
%
\begin{equation}
\label{lln(b)-excha3}
\sum_{j=1}^{\ulmn} \frac{{\overline S}_j-(\nu-1) j}{{\overline
S}_j(\nu-1)j}
=o_{\prob}\bigl(\sqrt{\log{\ulmn}}\bigr).
\end{equation}
Indeed, assuming that ${\overline S}_j>\vep j$ for all $j>j_0$,
independent of $n$
(recall Lem\-ma~\ref{lem-cond-LD-est}), we can use the bound
%
%e7.72 ###
%
\begin{eqnarray}
\label{borelbound}
\sum_{j=j_0}^{\ulmn} \frac{|{\overline S}_j-(\nu-1) j|}{{\overline
S}_j(\nu-1)j}
&\leq& C\sum_{j=j_0}^{\ulmn} \frac{|{\overline S}_j-(\nu-1)
j|}{j^2}\nonumber\\[-8pt]\\[-8pt]
&\leq& C\sum_{j=j_0}^{\ulmn} \frac{|{\overline S}^*_j|}{j^2}+\Op
(|\nu
-\overline{\nu}_n|\log{\olmn}),\nonumber
\end{eqnarray}
where ${\overline S}^*_j={\overline S}_j-(\overline{\nu}_n-1)j,$ is
for fixed $n$
the sum of i.i.d. random variables with mean $0$. Combining~(\ref
{nun-n-conv}) and~(\ref{olnun-n-conv}),
we obtain that $\Op(|\nu-\nu_n|\log{\olmn})=o_{\prob}(1)$,
so we are left to bound the first contribution
in \eqref{borelbound}.

According to the Marcinkiewicz--Zygmund inequality [recall~(\ref
{Marcinkiewicz-Zygmund})],
for $a\in(1,2)$,
\begin{eqnarray*}
\expec_n[|{\overline S}^*_j|^a]&\leq& B^*_a \expec\Biggl[\sum_{k=1}^j
\bigl({\overline X}_k-({\overline\nu}_n-1)\bigr)^2 \Biggr
]^{a/2}\\
&\le&
B^*_a\sum_{k=1}^j \expec_n [
|{\overline X}_k-(\overline{\nu}_n-1)]|^a ]
=jC_a \expec_n [
|{\overline X}_1-(\overline{\nu}_n-1)|^a ],
\end{eqnarray*}
%
%where the second inequality follows from the fact that, for $0\le
%r=a/2\le1$,
%%
%(x+y)^r\le(|x|+|y|)^r\le|x|^r+|y|^r.
%
When we take $1<a<\tau-2$, where $\tau-2>1$, then uniformly in $n$,
we have that $\expec_n[|{\overline X}_1-\overline{\nu}_n|^a]<c_a$
because
\begin{eqnarray*}
\expec_n[|{\overline X}_1|^a]&=&\sum_{s=1}^\infty s^a g_s^{(n)}
=\frac1{L_n}\sum_{i=1}^n D_i^a(D_i-1)\\
&\leq&\frac1{L_n}\sum_{i=1}^n
D_i^{a+1}\convas
\frac{\expec[D_1^{a+1}]}{\mu}<\infty,
\end{eqnarray*}
since $a<\tau-2$, so that
%
%e7.73 ###
%
\begin{eqnarray} \expec_n \Biggl[\sum_{j=1}^{\olmn} \frac
{|{\overline S}_j^*|}{j^{2}} \Biggr] &\leq&\sum_{j=1}^{\olmn} \frac
{\expec_n[|{\overline S}_j^*|^a]^{1/a}}{j^{2}}\nonumber\\[-8pt]\\[-8pt]
&=& \sum_{j=1}^{\olmn}
\frac{(c_a)^{1/a} \expec_n[|{\overline X}_1-(\nu
_n-1)|^a]^{1/a}}{j^{2-1/a}} <\infty,\nonumber
\end{eqnarray}
since $a>1$, and the last bound being true a.s. and uniform in $n$.
The proof for ${\underline S}_j$ is identical, where now,
instead of~(\ref{small-order-stats}),
we use that there exists $\alpha>0$ such that, w.h.p.,
%
%e7.74 ###
%
\begin{equation} \label{large-order-stats} \frac1{L_n} \sum
_{j=n-\olmn+1}^{n} D_{(j)}^2=o_{\prob}(n^{-\alpha}),
\end{equation}
using the argument in~(\ref{bound-order-stats})--(\ref{MI-order-stats}).
\end{pf}
%
%%%%%%%%%%%%%%%%%%%%%%%%%%%%%%%%%%%%%%%%%%%%%%%%%%%%%%%%%%%%%%%%%%%

\setcounter{equation}{0}
\setcounter{Theorem}{0}
\renewcommand{\theequation}{B.\arabic{equation}}
\renewcommand{\theTheorem}{B.\arabic{Theorem}}
\section*{Appendix B: On the deviation from a tree}\label{appB}

In this section, we do the necessary preliminaries needed for the proof
of Proposition~\ref{lemma:conn_edge} in Section~\ref{sec-lemma:conn_edge}.
One of the ingredients is writing $\prob(\CE_n>m)$ as the expectation
of the product of conditional probabilities [see \eqref{prod-CEn1} and
Lemma~\ref{lem-prod-CEn}].
A~second issue of Section~\ref{sec-lemma:conn_edge} is
to estimate the two error terms in \eqref{twee-termen}. We will deal
with these two
error terms in Lemma~\ref{cor-conv-zero-error}. Lemma~\ref{lem-Rm-ub}
is a preparation
for Lemma~\ref{cor-conv-zero-error} and gives an upper bound for the
expected number
of artificial stubs, which in turn is bounded by the expected number of
closed cycles.

In the statement of the following lemma, we recall that
$\Qprob^{(j)}_n$ denotes the conditional distribution
given $\SWG^{(1,2)}_j$ and $\{D_i\}_{i=1}^n$.

\begin{Lemma}[(Conditional product form tail probabilities $\CE_n$)]
\label{lem-prod-CEn}
%
%e7.75 ###
%
\begin{equation} \label{prod-CEn} \prob(\CE_n>m)=\expec\Biggl
[\prod
_{j=1}^m \Qprob^{(j)}_n(\CE_n>j|\CE_n>j-1) \Biggr].
\end{equation}
\end{Lemma}

\begin{pf} By the tower
property of conditional expectations, we can write
%
%e7.76 ###
%
\begin{eqnarray}\prob(\CE_n>m)&=&\expec\bigl[\Qprob^{(1)}_n(\CE
_n>m) \bigr]\nonumber\\[-8pt]\\[-8pt]
&=&\expec\bigl[\Qprob^{(1)}_n(\CE_n>1)\Qprob^{
(1)}_n(\CE_n>m|\CE_n>1) \bigr].\nonumber
\end{eqnarray}
Continuing this further, for all $1\leq k\leq m$,
%
%e7.78 ###
%e7.77 ###
%
\begin{eqnarray}
&&\Qprob^{(k)}_n(\CE_n>m|\CE_n>k)\nonumber\\
&&\qquad=\expec^{
(k)}_n \bigl[\Qprob^{(k+1)}_n(\CE_n>m|\CE_n>k) \bigr]\\
&&\qquad=\expec^{(k)}_n \bigl[\Qprob^{(k+1)}_n(\CE_n>k+1|\CE
_n>k)\Qprob^{(k+1)}_n(\CE_n>m|\CE_n>k+1) \bigr],\nonumber
\end{eqnarray}
where $\expec^{(k)}_n$ denotes the expectation w.r.t. $\Qprob
^{(k)}_n$.
In particular,
%
%e7.81 ###
%e7.80 ###
%e7.79 ###
%
\begin{eqnarray} \prob(\CE_n>m)&=&\expec\bigl[\Qprob^{(1)}_n(\CE
_n>m) \bigr]\nonumber\\
&=&\expec\bigl[\Qprob^{(1)}_n(\CE_n>1)\expec
^{(1)}_n \bigl[\Qprob^{(2)}_n(\CE_n>2|\CE_n>1)\nonumber\\
&&{}\hspace*{89pt}\times\Qprob
^{(2)}_n(\CE_n>m|\CE_n>2) \bigr] \bigr]\\
&=&\expec
\bigl[\Qprob^{(1)}_n(\CE_n>1)\Qprob^{(2)}_n(\CE_n>2|\CE
_n>1)\nonumber\\
&&{}\hspace*{8pt}\times
\Qprob^{(2)}_n(\CE_n>m|\CE_n>2) \bigr],\nonumber
\end{eqnarray}
where the last equality follows since $\Qprob^{(1)}_n(\CE_n>1)$
is measurable w.r.t. $\Qprob^{(2)}_n$
and the tower property. Continuing this indefinitely, we arrive at
(\ref{prod-CEn}).
\end{pf}

\begin{Lemma}[(The number of cycles closed)]
\label{lem-Rm-ub}
\textup{(a)} Fix $\tau\in(2,3)$. Then, w.h.p., there exist $\olmn$ with
$\olmn/ a_n\to\infty$ and $C>0$ such that
for all $m\leq\olmn$ and all $\vep_n\downarrow0$,
%
%e7.82 ###
%
\begin{equation} \label{m-Rm-bd(2,3)} \expec_n\bigl[R_{m}^{
(i)}-m\bigr]\leq\vep_n^{-1} \biggl(\frac{m}{a_n} \biggr)^2,
\qquad i=1,2.
\end{equation}

\textup{(b)} Fix $\tau>3$. Then, there exist $\olmn$ with $\olmn
/a_n\to
\infty$ and $C>0$ such that
for all $m\leq\olmn$,
%
%e7.83 ###
%
\begin{equation} \label{m-Rm-bd>3} \expec\bigl[R_{m}^{(i)}-m\bigr
]\leq
Cm^2/n,\qquad i=1,2.
\end{equation}
\end{Lemma}

\begin{pf} Observe that
%
%e7.84 ###
%
\begin{equation} \label{Art-bd} R_{m}^{(i)}-m\leq\sum_{j=1}^{m}
U_j,
\end{equation}
where $U_j$ is the indicator that a cycle is closed at time $j$.
Since closing a cycle means choosing an allowed stub, which occurs with
conditional probability at most $S_{j-1}^{(i)}/(L_n-2j-1)$, we find
that
%
%e7.85 ###
%
\begin{equation} \label{cond-expec} \expec\bigl[U_j| S_{j-1}^{
(i)},L_n\bigr]=S_{j-1}^{(i)}/(L_n-2j-1),
\end{equation}
so that
%
%e7.86 ###
%
\begin{equation} \label{BEnine} \expec_n\bigl[R_{m}^{(i)}-m\bigr
]\leq\sum
_{j=1}^{m} \expec_n[U_j] =\sum_{j=1}^{m} \expec_n\bigl[S_{j-1}^{
(i)}/(L_n-2j-1)\bigr].
\end{equation}
When $\tau>3$, and using that, since $j\leq\olmn=o(n)$,
we have $L_n-2j-1\geq2n-2j-1\geq n$ a.s., we arrive at
%
%e7.87 ###
%
\begin{equation} \expec\bigl[R_{m}^{(i)}-m\bigr] \leq\frac
{1}{n}\sum
_{j=1}^{m} \expec\bigl[S_{j-1}^{(i)}\bigr] \leq\frac{\mu
}{n}+\frac
{1}{n} \sum_{j=2}^{m} C(j-1) \leq Cm^2/n.
\end{equation}
When $\tau\in(2,3)$, we have to be a bit more careful. In this case,
we apply~(\ref{aim-lb-Sm(2,3)-ub})
to the r.h.s. of~(\ref{BEnine}), so that, w.h.p., and uniformly in $m$,
%
%e7.88 ###
%
\begin{equation} \expec_n\bigl[R_{m}^{(i)}-m\bigr]\leq\frac
{m^2}{n} \vep
_n^{-1} n^{(3-\tau)/(\tau-1)} =\vep_n^{-1} \biggl(\frac{m}{a_n}
\biggr)^2.
\end{equation}
This proves~(\ref{m-Rm-bd(2,3)}).
\end{pf}

\begin{Lemma}[(Treatment of error terms)]
\label{cor-conv-zero-error}
As $n\rightarrow\infty$, there exists $\olmn$ with $\olmn/ a_n\to
\infty,$ such that
%
%e7.89 ###
%
\begin{equation} \frac{\olmn}{n}\big|\Art_{a_n}^{(1)}\big|=o_{
\prob}(1), \qquad \frac{S_{a_n}^{(1)}}{n}\sum_{m=1}^{\olmn}
\frac{|\Art_{m}^{(2)}|}{S_{m}^{(2)}} =o_{\prob}(1).
\end{equation}
\end{Lemma}

\begin{pf}
We start with the first term. By Lemma~\ref{lem-Rm-ub}, for $\tau>3$,
%
%e7.90 ###
%
\begin{equation} \expec\bigl[\big|\Art_{m}^{(i)}\big|\bigr] \leq
\expec
\bigl[R_{m}^{(i)}-m\bigr]\leq Cm^2/n,\qquad m\le\olmn.
\end{equation}
As a result, we have that
%
%e7.91 ###
%
\begin{equation} \frac{\olmn}{n}\expec\bigl[\big|\Art
_{a_n}^{(1)}\big|\bigr]
\leq C\olmn^3/n^2=o(1).
\end{equation}
Again by Lemma~\ref{lem-Rm-ub}, but now for $\tau\in(2,3)$,
w.h.p. and uniformly in $m\leq\olmn$, where $\olmn$ is determined
in Lemma~\ref{lem-Rm-ub},
%
%e7.92 ###
%
\begin{equation} \frac{\olmn}{n}\expec_n\bigl[\big|\Art
_{a_n}^{(1)}\big|\bigr]
\leq\frac{\olmn}{n}\vep_n^{-1}=o(1),
\end{equation}
whenever $\vep_n^{-1}\rightarrow\infty$ sufficiently slowly.

Using \eqref{Art-bd} and $|\Art_m|\leq R_m-m$, and
using also that, w.h.p. and for all $j\leq m$,
$S_{m}^{(2)}\geq S_{j-1}^{(2)}$,
we obtain that
%
%e7.94 ###
%e7.93 ###
%
\begin{eqnarray} \label{mean-ratio-m-rep} \expec_n \biggl[\frac
{|\Art
_{m}^{(2)}|} {S_{m}^{(2)}} \biggr]
& \leq&
\expec_n \biggl[\frac
{\sum_{j=1}^{m} U_j} {S_{m}^{(2)}} \biggr]=\sum_{j=1}^{m}\expec
_n \biggl[\frac{U_j} {S_{m}^{(2)}} \biggr]\nonumber\\
&\leq&\sum
_{j=1}^{m}\expec_n \biggl[\frac{U_j} {S_{j-1}^{(2)}} \biggr]\leq
\sum_{j=1}^{m}\expec_n[1/(L_n-2j-1)]\\
&\leq& m/n,\nonumber
\end{eqnarray}
where we used~(\ref{cond-expec}) in the one-but-last inequality.

When $\tau>3$, we thus further obtain
%
%e7.95 ###
%
\begin{equation} \frac{1}{n}\sum_{m=1}^{\olmn}\expec\biggl[\frac
{|\Art_{m}^{(2)}|}{S_{m}^{(2)}} \biggr] \leq\frac{1}{n}\sum
_{m=1}^{\olmn}m/n=O(\olmn^2/n^2),
\end{equation}
so that, also using the bound on $S_{a_n}^{(1)}$ that holds w.h.p. as
proved in~(\ref{Encompl}),
%
%e7.96 ###
%
\begin{equation} \frac{1}{n}\sum_{m=1}^{\olmn} \frac{S_{a_n}^{
(1)}|\Art_{m}^{(2)}|}{S_{m}^{(2)}} =o_{\prob}(1).
\end{equation}

When $\tau\in(2,3)$, by~(\ref{mean-ratio-m-rep}),
%
%e7.97 ###
%
\begin{equation} \sum_{m=1}^{\olmn}\expec_n \biggl[\frac{|\Art
_{m}^{(2)}|}{S_{m}^{(2)}} \biggr] \leq\sum_{m=1}^{\olmn}
m/n\leq\olmn^2/n,
\end{equation}
so that, again using the bound on $S_{a_n}^{(1)}$ that holds w.h.p. as
proved in~(\ref{Encompl}),
%
%e7.98 ###
%
\begin{eqnarray} \frac{S_{a_n}^{(1)}}{n}\sum_{m=1}^{\olmn} \frac
{|\Art_{m}^{(2)}|}{S_{m}^{(2)}} &=&\Op\bigl(\eta_n^{-1}
n^{-2+(3-\tau)/(\tau-1)}a_n\olmn^2\bigr)\nonumber\\[-8pt]\\[-8pt]
&=&\Op\bigl(\eta_n^{-1} (\olmn
/a_n )^2 n^{-1/(\tau-1)}\bigr)=o_{\prob}(1),\nonumber
\end{eqnarray}
since $a_n=n^{(\tau-2)/(\tau-1)}$ and whenever
$\olmn/a_n, \eta_n^{-1}\rightarrow\infty$ sufficiently slowly such that
$n^{-1/(\tau-1)}\eta_n^{-1} (\olmn/a_n )^2=o(1)$.
\end{pf}

\setcounter{equation}{0}
\setcounter{Theorem}{0}
\renewcommand{\theequation}{C.\arabic{equation}}
\renewcommand{\theTheorem}{C.\arabic{Theorem}}
\section*{Appendix C: Weak convergence of the weight for $\tau>3$}\label{sec-weak-conv-weight-tau3}

In this section we prove Propositions~\ref{lemma:CLT-sum1}(b)
and~\ref{lemma:CLT-sum2}(b), for $\tau>3$. Moreover,
we show weak convergence of $C_n/a_n$ and prove~(\ref
{weight-conn-edge}) for $\tau>3$.
We start with Proposition~\ref{lemma:CLT-sum1}(b).

For this, we rewrite $T_m$ [compare \eqref{Tm-def}, with $s_i$
replaced by $S_i$],
%
%e7.99 ###
%
\begin{equation} \label{decom-Tm} T_m-\frac{1}{\nu-1}\log{m} =\sum
_{i=1}^m \frac{E_i-1}{S_i} + \Biggl[\sum_{i=1}^m \frac
{1}{S_i}-\frac
{1}{\nu-1}\log{m}. \Biggr]
\end{equation}
The second term on the r.h.s. of \eqref{decom-Tm}
converges a.s. to some $Y$ by~(\ref{as-conv-sum});
thus, it suffices to prove that $\sum_{i=1}^m (E_i-1)/S_i$ converges a.s.
For this, we use that the second moment equals, due to the independence
of $\{E_i\}_{i=1}^{\infty}$ and $\{S_i\}_{i=1}^{\infty}$
and the fact that $\expec[E_i]=\Var(E_i)=1$,
%
%e7.100 ###
%
\begin{equation} \expec\Biggl[ \Biggl(\sum_{i=1}^m \frac
{E_i-1}{S_i} \Biggr)^2 \Biggr] =\expec\Biggl[\sum_{i=1}^m 1/S_i^2
\Biggr],
\end{equation}
which converges uniformly in $m$. This shows that
%
%e7.101 ###
%
\begin{equation} T_m-\frac{1}{\nu-1}\log{m} \convd\sum
_{i=1}^\infty\frac{E_i-1}{S_i}+Y,
\end{equation}
which completes the proof for $T_m$ for $\tau>3$.

We continue the proof of Proposition~\ref{lemma:conn_edge} by showing
that, for $\tau>3$,
(\ref{Mn-lim-heur}) holds.

\begin{Lemma}[(Weak convergence of connection time)]
\label{lem-weak-conv-CE}
Fix $\tau>3$, then,
%
%e7.102 ###
%
\begin{equation} \label{CE-weak-conv} \CE_n/a_n\convd M,
\end{equation}
where $M$ has an exponential distribution with mean $\mu/(\nu
-1)$, that is,
%
%e7.103 ###
%
\begin{equation} \label{M-distr-def} \prob(M>x)=\exp\biggl\{-\frac
{\nu-1}{\mu} x \biggr\}.
\end{equation}
\end{Lemma}

\begin{pf}
The proof is somewhat sketchy; we leave the details to the reader.
%, and is a more precise version of the proof of Lemma~\ref{lem-CEn-ub}.
We again make use of the product structure in Lemma~\ref{lem-prod-CEn}
[recall~(\ref{prod-CEn})],
and simplify~(\ref{Qprob-m}), by taking complementary probabilities, to
%
%e7.104 ###
%
\begin{equation} \label{Qprob-m-rep} \Qprob^{(m)}_n(\CE
_n>m+1|\CE_n>m)\approx1-S_{a_n}^{(1)}/L_n.
\end{equation}
For $m\leq\olmn$, error terms that are left out can easily be seen to
be small by
Lemma~\ref{cor-conv-zero-error}. We next simplify by substitution of
$L_n=\mu n$, and using that ${\mathrm e}^{-x}\approx1-x$, for $x$
small, to obtain
that
%
%e7.105 ###
%
\begin{equation} \Qprob^{(m)}_n(\CE_n>m+1|\CE_n>m)\approx\exp
\bigl\{ -S_{a_n}^{(1)}/(\mu n) \bigr\}.
\end{equation}
Substituting the above approximation into~(\ref{prod-CEn}) for $m=a_n
x$ yields
%
%e7.106 ###
%
\begin{equation} \prob(\CE_n> a_n x)\approx\expec\biggl[\exp
\biggl\{
-a_n x\frac{S_{a_n}^{(1)}}{\mu n} \biggr\} \biggr]= \exp\biggl\{
-\frac{(\nu-1)}{\mu n}a_n^2 x \biggr\},
\end{equation}
where we approximate $S_{m}^{(1)}\approx(\nu-1)m$. Since $a_n=\sqrt{n}$,
we arrive at~(\ref{CE-weak-conv})--(\ref{M-distr-def}).
%\rightqed
\end{pf}

We now complete the proof of~(\ref{weight-conn-edge}) for $\tau>3$.
It is not hard to prove from~(\ref{decom-Tm}) that
%
%e7.107 ###
%
\begin{equation} \bigl(T^{(1)}_{a_n}- \gamma\log{a_n}, T^{
(2)}_{\CE_n}-\gamma\log{\CE_n} \bigr) \convd(X_1,X_2),
\end{equation}
where $(X_1, X_2)$ are two independent random variables with
distribution given by
%
%e7.108 ###
%
\begin{eqnarray}\label{eqn:dist-x} X_1&=&
\sum_{i=1}^\infty\frac{E_i-1}{\Sindep_i} +\lim
_{m\to\infty} \Biggl[ \Biggl(\sum_{i=1}^m 1/\Sindep_i \Biggr
)-\log
{m} \Biggr]\nonumber\\[-8pt]\\[-8pt]
&=&\sum_{i=1}^\infty\frac{E_i-1}{\Sindep_i} +\sum
_{i=1}^\infty\biggl(\frac{1}{\Sindep_i}-\frac{1}{(\nu-1)i} \biggr
) +
\gammaeuler, \nonumber
\end{eqnarray}
where $\gammaeuler$ is the Euler--Mascheroni constant.
By Lemma~\ref{lem-weak-conv-CE},
%
%e7.109 ###
%
\begin{equation} \bigl(T^{(1)}_{a_n}-\gamma\log{a_n},T^{
(2)}_{\CE_n}-\gamma\log{a_n} \bigr) \convd(X_1, X_2+\gamma
\log{M} ),
\end{equation}
where $M$ is the weak limit of $\CE_n/a_n$ defined in~(\ref{Mn-lim-heur}).
%By Lemma~\ref{lem-tildeTm-conv}, a similar argument applies to $
%$\widetilde T^{\sss(2)}_{\sss\lfloor\CE_n/2\rfloor}$, and shows that
% \eqn{
% (\widetilde T^{\sss(1)}_{\sss\lceil\CE_n/2\rceil}
% -\gamma\log{a_n},\widetilde T^{\sss(2)}_{\sss\lfloor\CE_n/2
% \convd(X_1+\gamma\log{(M/2)}-E_1/\nu, X_2+\gamma\log{(M/2)}-E_2/
% }
%where $E_1, E_2$ are two independent exponential random variables with
%mean 1, independent of all
%other random variables involved.
We conclude that
%
%e7.110 ###
%
\begin{equation} \Wn-\gamma\log{n} \convd V=X_1+X_2+\gamma\log{M}.
\end{equation}
Since $(\nu-1)M/\mu$ is an exponential variable with mean 1,
$\Lambda=\log{((\nu-1)M/\mu)}$ has a Gumbel distribution.

Finally let us derive the distribution of $X_i$. The random variables
$X_i$ are related to
a random variable $W$, which appears as a limit in a supercritical
continuous-time
branching process
as described in Section~\ref{sec-flow_tree}.
Indeed, denoting by $Z(t)$ the number of alive individuals in
a continuous-time branching process where the root has degree $D$
having distribution function $F$,
while all other vertices in the tree have degree
$\{\Bindep_i\}_{i=2}^{\infty}$, which are i.i.d.\
random variables with probability mass function $g$ in~(\ref{eqn:size-bias}).
Then, $W$ arises as
%
%e7.111 ###
%
\begin{equation} \label{Zt-mart} Z(t){\mathrm e}^{-(\nu-1)t} \convas
W.
\end{equation}
We note the following general results about the limiting distributional
asymptotics of continuous-time branching processes.

\begin{Proposition}[(The limiting random variables)]
\textup{(a)} The limiting random variable $W$ has the following
explicit construction:
%
%e7.112 ###
%
\begin{equation}
W = \sum_{j=1}^{D} \widetilde{W}_j {\mathrm e}^{-(\nu-1)\xi_j}.
\label{eqn:w-dist}
\end{equation}
Here $D$ has distribution $F$, $\xi_i$ are i.i.d. exponential random
variables with mean one independent of
$\widetilde{W}_i,$ which are independent and identically distributed with
Laplace transform $\phi(t) = \expec({\mathrm e}^{-t \widetilde{W}})$
given by the formula
%
%e7.113 ###
%
\begin{equation}
\qquad\phi^{-1}(x) = (1-x) \exp\biggl\{\int_1^x \biggl(\frac{\nu-1}{h(s)
-s} + \frac{1}{1-s} \biggr)\,ds \biggr\},\qquad 0< x\leq1,\label
{eq:lap-W}
\end{equation}
and $h(\cdot)$ is the probability generating function of the size-biased
probability mass function $g$ [see~(\ref{eqn:size-bias})].

\textup{(b)} Let $T_m$ be the random variables defined as
%
%e7.114 ###
%
\begin{equation} T_m = \sum_{i=1}^m E_i/\Sindep_i,
\end{equation}
where $E_i$ are i.i.d. exponential random variables with mean one,
and recall that $\Sindep_i$ is a random walk where the first step has
distribution $D$ where $D\sim F$ and the remaining increments have
distribution $B-1$ where $B$ has the size biased distribution. Then
%
%e7.115 ###
%
\begin{equation} T_m - \frac{\log{m}}{\nu-1} \convas-\frac{\log
(W/(\nu-1) )}{\nu-1},
\end{equation}
where $W$ is the martingale limit in~(\ref{Zt-mart}) in part \textup{(a)}.

\textup{(c)} The random variables $X_i$, $i=1,2$, are i.i.d. with $X_i
\stackrel{d}{=}
-\frac{\log(W/(\nu-1) )}{\nu-1}.$

\end{Proposition}
\begin{pf} These results follow from results about continuous-time
branching processes (everything relevant to this result is taken from
\cite{athreya}). Part (b) is proved in \cite{athreya}, Theorem 2, page 120.
To prove part (a) recall the continuous-time version of the
construction described in Section~\ref{sec-flow_tree}, where we shall
let $D\sim F$ denote the number of offspring of the initial root and,
for $i\geq2$, $B_i\sim g$, the size-biased biased probability mass
function~(\ref{eqn:size-bias}). Then note that for any $t$
sufficiently large we can decompose $Z(t)$, the number of alive nodes
at time $t$ as
%
%e7.116 ###
%
\begin{equation}
Z(t) {\mathrm e}^{-(\nu-1)t} = \sum_{i=1}^D
\widetilde{Z}_i({t-\xi_i}) {\mathrm e}^{-(\nu-1)t}.
\label{eq:dist-ft}
\end{equation}
Here $D$, $\xi_i$ and the processes $\widetilde{Z}_i(\cdot)$ are all
independent of each other, $D\sim F$ denotes the number of offspring of
the root, $\xi_i$ are lifetimes of these\vspace*{1.5pt} offspring and are distributed
as i.i.d. exponential random variables with mean 1 and $\widetilde
{Z}_j(\cdot)$, corresponding to the subtrees attached below offspring
$j$ of the root, are independent continuous-time branching processes
where each individual lives for an exponential mean $1$ amount of time
and then dies, giving birth to a random number of offspring where the
number of offspring has distribution $B\sim g$ as in~(\ref{eqn:size-bias}).

Now known results (see \cite{athreya}, Theorem~1, page 111 and
Theorem~3, page 116) imply that
\[
\widetilde{Z}_i(t){\mathrm e}^{-(\nu-1)t} \convas\widetilde{W}_i,
\]
where $ \widetilde{W}_i $ have Laplace transform given by~(\ref{eq:lap-W}).
Part (a) now follows by comparing~(\ref{eqn:w-dist}) with~(\ref{eq:dist-ft}).

Part (c) follows from part (b) and observing that
\[
T_m - \frac{1}{(\nu-1)}\log{m} = \sum_{i=1}^m \frac{E_i-1}{\Sindep
_i} + \sum_{i=1}^m \frac{1}{\Sindep_i} - \frac{1}{(\nu-1)}\log{m},
\]
and a comparison with~(\ref{eqn:dist-x}). This completes the proof.
\end{pf}

Thus, with $\Lambda$ a Gumbel distribution, the explicit
distribution of
the re-centered minimal weight paths is given by
%
%e7.117 ###
%
\begin{equation}\qquad\quad V= -\frac{\log(W_1/(\nu-1) )}{\nu
-1}-\frac{\log (W_2/(\nu-1) )}{\nu-1} + \gamma\Lambda
-\gamma\log{(\nu-1)/\mu},
\end{equation}
since $\log M=\Lambda-\log((\nu-1)/\mu)$. Rearranging terms
establishes the claims on the limit $V$ below Theorem~\ref{main>3}, and
completes the proof of~(\ref{weight-conn-edge})
in Proposition~\ref{lemma:conn_edge}(b) for $\tau>3$.
\end{appendix}

\section*{Acknowledgments}
The authors thank Allan Sly for help with a preliminary version of the paper.
G. H. thanks the Collegium Budapest for the opportunity to work on the
revision of this paper during his visit.

% imsref loaded by smiklovaite, 2010-03-04 16:26:09
%

\printaddresses


\begin{thebibliography}{37}

%b1 ###
\bibitem{AdaBroLug09}
%
\begin{barticle}[auto:SpringerTagBib|2009-01-14|16:51:27]
\bauthor{\bsnm{Addario-Berry},~\bfnm{L.}\binits{L.}},
\bauthor{\bsnm{Broutin},~\bfnm{N.}\binits{N.}} \AND
\bauthor{\bsnm{Lugosi},~\bfnm{G.}\binits{G.}}
(\byear{2010}).
\btitle{The longest minimum-weight path in a complete graph.
Preprint}.
\bjournal{Combin. Probab. Comput.}
\bvolume{19}
\bpages{1--19}.
\end{barticle}
%
\endbibitem

%b2 ###
\bibitem{AthKar67}
%
\begin{barticle}[mr]
\bauthor{\bsnm{Athreya},~\bfnm{Krishna~B.}\binits{K.~B.}} \AND
\bauthor{\bsnm{Karlin},~\bfnm{Samuel}\binits{S.}}
(\byear{1967}).
\btitle{Limit theorems for the split times of branching processes}.
\bjournal{J. Math. Mech.}
\bvolume{17}
\bpages{257--277}.
\bid{mr={0216592}}
\end{barticle}
%
\endbibitem

%b3 ###
\bibitem{athreya}
%
\begin{bbook}[mr]
\bauthor{\bsnm{Athreya},~\bfnm{K.~B.}\binits{K.~B.}} \AND
\bauthor{\bsnm{Ney},~\bfnm{P.~E.}\binits{P.~E.}}
(\byear{2004}).
\btitle{Branching Processes}.
\bpublisher{Dover}, \baddress{Mineola, NY}.
\bid{mr={2047480}}
\end{bbook}
%
\endbibitem

%b4 ###
\bibitem{BerSid09}
%
\begin{barticle}[mr]
\bauthor{\bsnm{Bertoin},~\bfnm{Jean}\binits{J.}} \AND
\bauthor{\bsnm{Sidoravicius},~\bfnm{Vladas}\binits{V.}}
(\byear{2009}).
\btitle{The structure of typical clusters in large sparse random
configurations}.
\bjournal{J. Stat. Phys.}
\bvolume{135}
\bpages{87--105}.
\bid{doi={10.1007/s10955-009-9728-y}, mr={2505727}}
\end{barticle}
%
\endbibitem

%b5 ###
\bibitem{vcg-random-shanky}
%
\begin{barticle}[mr]
\bauthor{\bsnm{Bhamidi},~\bfnm{Shankar}\binits{S.}}
(\byear{2008}).
\btitle{First passage percolation on locally treelike networks. {I}. {D}ense
random graphs}.
\bjournal{J. Math. Phys.}
\bvolume{49}
\bpages{125218}.
\bid{doi={10.1063/1.3039876}, mr={2484349}}
\end{barticle}
%
\endbibitem

%b6 ###
\bibitem{BhaHofHoo09a}
%
\begin{bmisc}[auto:SpringerTagBib|2009-01-14|16:51:27]
\bauthor{\bsnm{Bhamidi},~\bfnm{S.}\binits{S.}}, \bauthor
{\bparticle{van~der
}\bsnm{Hofstad},~\bfnm{R.}\binits{R.}} \AND
\bauthor{\bsnm{Hooghiemstra},~\bfnm{G.}\binits{G.}}
(\byear{2009}).
\bhowpublished{Extreme value theory, {P}oisson--{D}irichlet
distributions and
first passage percolation on random networks. Preprint}.
\end{bmisc}
%
\endbibitem

%b7 ###
\bibitem{BinGolTeu89}
%
\begin{bbook}[mr]
\bauthor{\bsnm{Bingham},~\bfnm{N.~H.}\binits{N.~H.}},
\bauthor{\bsnm{Goldie},~\bfnm{C.~M.}\binits{C.~M.}} \AND
\bauthor{\bsnm{Teugels},~\bfnm{J.~L.}\binits{J.~L.}}
(\byear{1989}).
\btitle{Regular Variation}.
\bseries{Encyclopedia of Mathematics and Its Applications}
\bvolume{27}.
\bpublisher{Cambridge Univ. Press}, \baddress{Cambridge}.
\bid{mr={1015093}}
\end{bbook}
%
\endbibitem

%b8 ###
\bibitem{Boll01}
%
\begin{bbook}[mr]
\bauthor{\bsnm{Bollob{\'a}s},~\bfnm{B{\'e}la}\binits{B.}}
(\byear{2001}).
\btitle{Random Graphs},
\bedition{2nd} ed.
\bseries{Cambridge Studies in Advanced Mathematics}
\bvolume{73}.
\bpublisher{Cambridge Univ. Press}, \baddress{Cambridge}.
\bid{mr={1864966}}
\end{bbook}
%
\endbibitem

%b9 ###
\bibitem{weak-strong-diso}
%
\begin{barticle}[auto:SpringerTagBib|2009-01-14|16:51:27]
\bauthor{\bsnm{Braunstein},~\bfnm{L.~A.}\binits{L.~A.}},
\bauthor{\bsnm{Buldyrev},~\bfnm{S.~V.}\binits{S.~V.}},
\bauthor{\bsnm{Cohen},~\bfnm{R.}\binits{R.}},
\bauthor{\bsnm{Havlin},~\bfnm{S.}\binits{S.}} \AND
\bauthor{\bsnm{Stanley},~\bfnm{H.~E.}\binits{H.~E.}}
(\byear{2003}).
\btitle{Optimal paths in disordered complex networks}.
\bjournal{Phys. Rev. Lett.}
\bvolume{91}
\bpages{168701}.
\end{barticle}
%
\endbibitem

%b10 ###
\bibitem{buhler}
%
\begin{barticle}[mr]
\bauthor{\bsnm{B{\"u}hler},~\bfnm{Wolfgang~J.}\binits{W.~J.}}
(\byear{1971}).
\btitle{Generations and degree of relationship in supercritical {M}arkov
branching processes}.
\bjournal{Z. Wahrsch. Verw. Gebiete}
\bvolume{18}
\bpages{141--152}.
\bid{mr={0295444}}
\end{barticle}
%
\endbibitem

%b11 ###
\bibitem{ChuLu03}
%
\begin{barticle}[mr]
\bauthor{\bsnm{Chung},~\bfnm{Fan}\binits{F.}} \AND
\bauthor{\bsnm{Lu},~\bfnm{Linyuan}\binits{L.}}
(\byear{2003}).
\btitle{The average distance in a random graph with given expected degrees}.
\bjournal{Internet Math.}
\bvolume{1}
\bpages{91--113}.
\bid{mr={2076728}}
\end{barticle}
%
\endbibitem

%b12 ###
\bibitem{ChuLu06c}
%
\begin{bbook}[mr]
\bauthor{\bsnm{Chung},~\bfnm{Fan}\binits{F.}} \AND
\bauthor{\bsnm{Lu},~\bfnm{Linyuan}\binits{L.}}
(\byear{2006}).
\btitle{Complex Graphs and Networks}.
\bseries{CBMS Regional Conference Series in Mathematics}
\bvolume{107}.
\bpublisher{Amer. Math. Soc.},
\baddress{Providence, RI}.
\bid{mr={2248695}}
\end{bbook}
%
\endbibitem

%b13 ###
\bibitem{CohHav03}
%
\begin{bmisc}[auto:SpringerTagBib|2009-01-14|16:51:27]
\bauthor{\bsnm{Cohen},~\bfnm{R.}\binits{R.}} \AND
\bauthor{\bsnm{Havlin},~\bfnm{S.}\binits{S.}}
(\byear{2003}).
\bhowpublished{Scale-free networks are ultrasmall. \textit{Phys. Rev. Lett.}
\textbf{90} 058701}.
\end{bmisc}
%
\endbibitem

%b14 ###
\bibitem{DeiEskHofHoo09}
%
\begin{barticle}[mr]
\bauthor{\bsnm{Deijfen},~\bfnm{Maria}\binits{M.}}, \bauthor
{\bparticle{van~den
}\bsnm{Esker},~\bfnm{Henri}\binits{H.}}, \bauthor{\bparticle{van~der
}\bsnm{Hofstad},~\bfnm{Remco}\binits{R.}} \AND
\bauthor{\bsnm{Hooghiemstra},~\bfnm{Gerard}\binits{G.}}
(\byear{2009}).
\btitle{A preferential attachment model with random initial degrees}.
\bjournal{Ark. Mat.}
\bvolume{47}
\bpages{41--72}.
\bid{doi={10.1007/s11512-007-0067-4}, mr={2480915}}
\end{barticle}
%
\endbibitem

%b15 ###
\bibitem{DinKimLubPer09}
%
\begin{bmisc}[auto:SpringerTagBib|2009-01-14|16:51:27]
\bauthor{\bsnm{Ding},~\bfnm{J.}\binits{J.}},
\bauthor{\bsnm{Kim},~\bfnm{J.~H.}\binits{J.~H.}},
\bauthor{\bsnm{Lubetzky},~\bfnm{E.}\binits{E.}} \AND
\bauthor{\bsnm{Peres},~\bfnm{Y.}\binits{Y.}}
(\byear{2009}).
\bhowpublished{Diameters in supercritical random graphs via first passage
percolation. Preprint.
Available at}
\url{http://arxiv.org/abs/0906.1840}.
\end{bmisc}
%
\endbibitem

%b16 ###
\bibitem{DomHofHoo08}
%
\begin{barticle}[auto:SpringerTagBib|2009-01-14|16:51:27]
\bauthor{\bsnm{Dommers},~\bfnm{S.}\binits{S.}}, \bauthor
{\bparticle{van~der
}\bsnm{Hofstad},~\bfnm{R.}\binits{R.}} \AND
\bauthor{\bsnm{Hooghiemstra},~\bfnm{G.}\binits{G.}}
(\byear{2010}).
\btitle{Diameters in preferential attachment graphs}.
\bjournal{J. Stat. Phys.}
\bvolume{139}
\bpages{72--107}.
\end{barticle}
%
\endbibitem

%b17 ###
\bibitem{Durr06}
%
\begin{bbook}[mr]
\bauthor{\bsnm{Durrett},~\bfnm{Rick}\binits{R.}}
(\byear{2007}).
\btitle{Random Graph Dynamics}.
\bpublisher{Cambridge Univ. Press}, \baddress{Cambridge}.
\bid{mr={2271734}}
\end{bbook}
%
\endbibitem

%b18 ###
\bibitem{FerRam04}
%
\begin{barticle}[mr]
\bauthor{\bsnm{Fernholz},~\bfnm{Daniel}\binits{D.}} \AND
\bauthor{\bsnm{Ramachandran},~\bfnm{Vijaya}\binits{V.}}
(\byear{2007}).
\btitle{The diameter of sparse random graphs}.
\bjournal{Random Structures Algorithms}
\bvolume{31}
\bpages{482--516}.
\bid{doi={10.1002/rsa.20197}, mr={2362640}}
\end{barticle}
%
\endbibitem

%b19 ###
\bibitem{Gut05}
%
\begin{bbook}[mr]
\bauthor{\bsnm{Gut},~\bfnm{Allan}\binits{A.}}
(\byear{2005}).
\btitle{Probability: A Graduate Course}.
\bpublisher{Springer}, \baddress{New York}.
\bid{mr={2125120}}
\end{bbook}
%
\endbibitem

%b20 ###
\bibitem{ham-wel}
%
\begin{bincollection}[vtex]
\bauthor{\bsnm{Hammersley},~\bfnm{J.~M.}\binits{J.~M.}} \AND
\bauthor{\bsnm{Welsh},~\bfnm{D.~J.~A.}\binits{D.~J.~A.}}
(\byear{1965}).
\btitle{First-passage percolation, subadditive processes, stochastic networks,
and generalized renewal theory}.
In
\bbooktitle{Bernouilli--Bayes--Laplace: Anniversary Volume}
\bpages{61--110}.
\bpublisher{Springer}, \baddress{New York}.
\bid{mr={0198576}}
\end{bincollection}
%
\endbibitem

%b21 ###
\bibitem{douglas-fpp}
%
\begin{bincollection}[mr]
\bauthor{\bsnm{Howard},~\bfnm{C.~Douglas}\binits{C.~D.}}
(\byear{2004}).
\btitle{Models of first-passage percolation}.
In \bbooktitle{Probability on Discrete Structures}.
\bseries{Encyclopaedia Math. Sci.}
\bvolume{110}
\bpages{125--173}.
\bpublisher{Springer}, \baddress{Berlin}.
\bid{mr={2023652}}
\end{bincollection}
%
\endbibitem

%b22 ###
\bibitem{janson123}
%
\begin{barticle}[mr]
\bauthor{\bsnm{Janson},~\bfnm{Svante}\binits{S.}}
(\byear{1999}).
\btitle{One, two and three times {$\log n/n$} for paths in a complete graph
with random weights}.
\bjournal{Combin. Probab. Comput.}
\bvolume{8}
\bpages{347--361}.
\bid{doi={10.1017/S0963548399003892}, mr={1723648}}
\end{barticle}
%
\endbibitem

%b23 ###
\bibitem{Jans06b}
%
\begin{barticle}[mr]
\bauthor{\bsnm{Janson},~\bfnm{Svante}\binits{S.}}
(\byear{2009}).
\btitle{The probability that a random multigraph is simple}.
\bjournal{Combin. Probab. Comput.}
\bvolume{18}
\bpages{205--225}.
\bid{doi={10.1017/S0963548308009644}, mr={2497380}}
\end{barticle}
%
\endbibitem

%b24 ###
\bibitem{JanLuc07}
%
\begin{barticle}[mr]
\bauthor{\bsnm{Janson},~\bfnm{Svante}\binits{S.}} \AND
\bauthor{\bsnm{Luczak},~\bfnm{Malwina~J.}\binits{M.~J.}}
(\byear{2009}).
\btitle{A new approach to the giant component problem}.
\bjournal{Random Structures Algorithms}
\bvolume{34}
\bpages{197--216}.
\bid{doi={10.1002/rsa.20231}, mr={2490288}}
\end{barticle}
%
\endbibitem

%b25 ###
\bibitem{JanLucRuc00}
%
\begin{bbook}[mr]
\bauthor{\bsnm{Janson},~\bfnm{Svante}\binits{S.}},
\bauthor{\bsnm{{\L}uczak},~\bfnm{Tomasz}\binits{T.}} \AND
\bauthor{\bsnm{Rucinski},~\bfnm{Andrzej}\binits{A.}}
(\byear{2000}).
\btitle{Random Graphs}.
\bpublisher{Wiley}, \baddress{New York}.
\bid{mr={1782847}}
\end{bbook}
%
\endbibitem

%b26 ###
\bibitem{MolRee95}
%
\begin{binproceedings}[mr]
\bauthor{\bsnm{Molloy},~\bfnm{Michael}\binits{M.}} \AND
\bauthor{\bsnm{Reed},~\bfnm{Bruce}\binits{B.}}
(\byear{1995}).
\btitle{A critical point for random graphs with a given degree sequence}.
In \bbooktitle{Proceedings of the {S}ixth {I}nternational {S}eminar on {R}andom
{G}raphs and {P}robabilistic {M}ethods in {C}ombinatorics and {C}omputer
{S}cience, ``{R}andom {G}raphs'\textit{93}'' ({P}ozna\'n, \textit{1993})}
\bvolume{6}
\bpages{161--179}.
\bpublisher{Wiley},
\baddress{New York}.
\bid{mr={1370952}}
\end{binproceedings}
%
\endbibitem

%b27 ###
\bibitem{MolRee98}
%
\begin{barticle}[mr]
\bauthor{\bsnm{Molloy},~\bfnm{Michael}\binits{M.}} \AND
\bauthor{\bsnm{Reed},~\bfnm{Bruce}\binits{B.}}
(\byear{1998}).
\btitle{The size of the giant component of a random graph with a given degree
sequence}.
\bjournal{Combin. Probab. Comput.}
\bvolume{7}
\bpages{295--305}.
\bid{doi={10.1017/S0963548398003526}, mr={1664335}}
\end{barticle}
%
\endbibitem

%b28 ###
\bibitem{NorRei06}
%
\begin{barticle}[mr]
\bauthor{\bsnm{Norros},~\bfnm{Ilkka}\binits{I.}} \AND
\bauthor{\bsnm{Reittu},~\bfnm{Hannu}\binits{H.}}
(\byear{2006}).
\btitle{On a conditionally {P}oissonian graph process}.
\bjournal{Adv. in Appl. Probab.}
\bvolume{38}
\bpages{59--75}.
\bid{doi={10.1239/aap/1143936140}, mr={2213964}}
\end{barticle}
%
\endbibitem

%b29 ###
\bibitem{ReiNor04}
%
\begin{barticle}[auto:SpringerTagBib|2009-01-14|16:51:27]
\bauthor{\bsnm{Reittu},~\bfnm{H.}\binits{H.}} \AND
\bauthor{\bsnm{Norros},~\bfnm{I.}\binits{I.}}
(\byear{2004}).
\btitle{On the power law random graph model of massive data networks}.
\bjournal{Performance Evaluation}
\bvolume{55}
\bpages{3--23}.
\end{barticle}
%
\endbibitem

%b30 ###
\bibitem{smythe-fpp}
%
\begin{bbook}[auto:SpringerTagBib|2009-01-14|16:51:27]
\bauthor{\bsnm{Smythe},~\bfnm{R.~T.}\binits{R.~T.}} \AND
\bauthor{\bsnm{Wierman},~\bfnm{J.~C.}\binits{J.~C.}}
(\byear{1978}).
\btitle{First-Passage Percolation on the Square Lattice}.
\bseries{Lecture Notes in Math.}
 \bvolume{671}.
 \bpublisher{Springer},
 \baddress{Berlin}.
\end{bbook}
%
\endbibitem

%b31 ###
\bibitem{hofs2}
%
\begin{barticle}[mr]
\bauthor{\bparticle{van~den }\bsnm{Esker},~\bfnm{Henri}\binits{H.}},
\bauthor{\bparticle{van~der }\bsnm{Hofstad},~\bfnm{Remco}\binits{R.}},
\bauthor{\bsnm{Hooghiemstra},~\bfnm{Gerard}\binits{G.}} \AND
\bauthor{\bsnm{Znamenski},~\bfnm{Dmitri}\binits{D.}}
(\byear{2005}).
\btitle{Distances in random graphs with infinite mean degrees}.
\bjournal{Extremes}
\bvolume{8}
\bpages{111--141}.
\bid{doi={10.1007/s10687-006-7963-z}, mr={2275914}}
\end{barticle}
%
\endbibitem

%b32 ###
\bibitem{hofs-erdos-fpp}
%
\begin{barticle}[mr]
\bauthor{\bparticle{van~der }\bsnm{Hofstad},~\bfnm{Remco}\binits{R.}},
\bauthor{\bsnm{Hooghiemstra},~\bfnm{Gerard}\binits{G.}} \AND
\bauthor{\bsnm{Van~Mieghem},~\bfnm{Piet}\binits{P.}}
(\byear{2001}).
\btitle{First-passage percolation on the random graph}.
\bjournal{Probab. Engrg. Inform. Sci.}
\bvolume{15}
\bpages{225--237}.
\bid{doi={10.1017/S026996480115206X}, mr={1828576}}
\end{barticle}
%
\endbibitem

%b33 ###
\bibitem{hofs-flood}
%
\begin{barticle}[mr]
\bauthor{\bparticle{van~der }\bsnm{Hofstad},~\bfnm{Remco}\binits{R.}},
\bauthor{\bsnm{Hooghiemstra},~\bfnm{Gerard}\binits{G.}} \AND
\bauthor{\bsnm{Van~Mieghem},~\bfnm{Piet}\binits{P.}}
(\byear{2002}).
\btitle{The flooding time in random graphs}.
\bjournal{Extremes}
\bvolume{5}
\bpages{111--129}.
\bid{doi={10.1023/A:1022175620150}, mr={1965974}}
\end{barticle}
%
\endbibitem

%b34 ###
\bibitem{hofs3}
%
\begin{barticle}[mr]
\bauthor{\bparticle{van~der }\bsnm{Hofstad},~\bfnm{Remco}\binits{R.}},
\bauthor{\bsnm{Hooghiemstra},~\bfnm{Gerard}\binits{G.}} \AND
\bauthor{\bsnm{Van~Mieghem},~\bfnm{Piet}\binits{P.}}
(\byear{2005}).
\btitle{Distances in random graphs with finite variance degrees}.
\bjournal{Random Structures Algorithms}
\bvolume{27}
\bpages{76--123}.
\bid{doi={10.1002/rsa.20063}, mr={2150017}}
\end{barticle}
%
\endbibitem

%b35 ###
\bibitem{hofs1}
%
\begin{barticle}[mr]
\bauthor{\bparticle{van~der }\bsnm{Hofstad},~\bfnm{Remco}\binits{R.}},
\bauthor{\bsnm{Hooghiemstra},~\bfnm{Gerard}\binits{G.}} \AND
\bauthor{\bsnm{Znamenski},~\bfnm{Dmitri}\binits{D.}}
(\byear{2007}).
\btitle{Distances in random graphs with finite mean and infinite variance
degrees}.
\bjournal{Electron. J. Probab.}
\bvolume{12}
\bpages{703--766}
\bmisc{(electronic)}.
\bid{mr={2318408}}
\end{barticle}
%
\endbibitem

%b36 ###
\bibitem{VanHooHof00}
%
\begin{bmisc}[auto:SpringerTagBib|2009-01-14|16:51:27]
\bauthor{\bsnm{Van~Mieghem},~\bfnm{P.}\binits{P.}},
\bauthor{\bsnm{Hooghiemstra},~\bfnm{G.}\binits{G.}} \AND
\bauthor{\bparticle{van~der }\bsnm{Hofstad},~\bfnm{R.}\binits{R.}}
(\byear{2000}).
\bhowpublished{A scaling law for the hopcount. Technical
report~2000125, Delft
Univ. Technology. Available at}
\url{http://www.nas.ewi.tudelft.nl/people/Piet/}.
\end{bmisc}
%
\endbibitem

%b37 ###
\bibitem{wastlund}
%
\begin{bmisc}[auto:SpringerTagBib|2009-01-14|16:51:27]
\bauthor{\bsnm{W{\"a}stlund},~\bfnm{J.}\binits{J.}}
(\byear{2006}).
\bhowpublished{Random assignment and shortest path problems}.
\bhowpublished{In
\textit{Proceedings of the Fourth Colloquium on Mathematics and
Computer Science, Institut. \'Elie Cartan, Nancy France. DMTCS Proc.}
\textbf{AG}
31--38.}
\end{bmisc}
%
\endbibitem

\end{thebibliography}
\end{document}